\documentclass[10pt,a4paper]{article}

\usepackage[utf8]{inputenc}
\usepackage{authblk}
\usepackage{amsthm,amsmath,amsfonts,amssymb,mathtools}
\usepackage[numbers]{natbib}
\usepackage{tikz}
\usetikzlibrary{arrows}
\usepackage{enumitem}
\usepackage{bbm}

\usepackage{xcolor} 
\usepackage[colorlinks,citecolor=blue, linkcolor = blue,urlcolor=blue]{hyperref}\usepackage{hyperref}
\usepackage{graphicx}
\usepackage{titlesec}

\titleformat{\section}
  {\normalfont\fontsize{12}{15}\bfseries}{\thesection}{1em}{}

\newif\ifarxiv\arxivtrue

\usetikzlibrary{3d}
\makeatletter
\tikzoption{canvas is xy plane at z}[]{%
  \def\tikz@plane@origin{\pgfpointxyz{0}{0}{#1}}%
  \def\tikz@plane@x{\pgfpointxyz{1}{0}{#1}}%
  \def\tikz@plane@y{\pgfpointxyz{0}{1}{#1}}%
  \tikz@canvas@is@plane
}
\makeatother

\tikzset{xyp/.style={canvas is xy plane at z=5}}
\tikzset{xzp/.style={canvas is xz plane at y=0}}
\tikzset{yzp/.style={canvas is yz plane at x=6}}

\theoremstyle{plain}
\newtheorem{theorem}{Theorem}[section]
\newtheorem{proposition}[theorem]{Proposition}
\newtheorem{lemma}[theorem]{Lemma}
\newtheorem{corollary}[theorem]{Corollary}

\newtheorem{remark}{Remark}[section]
\newtheorem{assumption}{Assumption}[section]
\newtheorem{example}{Example}[section]

\newcommand{\supnorm}[1]{{\left\lVert#1\right\rVert}_{\infty}}
\newcommand{\nnorm}[1]{\left\lVert#1\right\rVert}
\newcommand{\norm}[1]{\left\lvert#1\right\rvert}

\newcommand{\specnorm}[1]{{\left\lVert#1\right\rVert}_{2}}

\newcommand{\E}[1]{\mathbb{E}\left[#1\right]}
\newcommand{\var}[1]{\mathrm{Var}\left(#1\right)}

\newcommand{\toinp}[0]{\xrightarrow[n \to \infty]{\mathbb{P}}}
\newcommand{\toind}[0]{\xrightarrow[n \to \infty]{\mathrm{d}}}
\newcommand{\as}[0]{\mathbb{P}\text{ a.s.}}

\newcommand{\covmat}[0]{\Sigma_{n,\theta}}
\newcommand{\shiftcovmat}[0]{\Sigma_{h(n),\theta}}

\newcommand{\rcovmat}[0]{\Sigma_{\theta}(s_{(n)})}

\newcommand{\rrcovmat}[0]{\Sigma_{\theta}(S_{(n)}(\omega))}

\newcommand{\acovmat}[0]{\widetilde{\Sigma}_{n,\theta}}
\newcommand{\shiftacovmat}[0]{\widetilde{\Sigma}_{h(n),\theta}}

\newcommand{\racovmat}[0]{\widetilde{\Sigma}_{\theta}(s_{(n)})}

\newcommand{\tcovmat}[0]{\Sigma_{n,\theta_{0}}}
\newcommand{\shifttcovmat}[0]{\Sigma_{h(n),\theta_{0}}}

\newcommand{\rtcovmat}[0]{\Sigma_{\theta_{0}}(s_{(n)})}

\newcommand{\rracovmat}[0]{\widetilde{\Sigma}_{\theta}(S_{(n)}(\omega))}

\newcommand{\shifttacovmat}[0]{\widetilde{\Sigma}_{h(n),\theta_{0}}}

\newcommand{\ddet}[1]{\operatorname{det}\left(#1\right)}
\DeclareMathOperator{\tr}{tr}

\newcommand{\eeexp}[0]{\operatorname{e}}
\newcommand{\llog}[1]{\operatorname{log}\left(#1\right)}
\newcommand{\tonormal}[0]{\xrightarrow[]{n \to \infty}}
\newcommand{\tonormalm}[0]{\xrightarrow[]{m \to \infty}}

\DeclareMathOperator*{\argmin}{argmin}

\newcommand{\ii}[0]{\mathrm{i}}
\newcommand{\natnum}[0]{\mathbb{N}_{+}}

\newcommand{\orange}[1]{\textcolor{orange}{#1}}
\renewcommand{\orange}[1]{{#1}}

\renewenvironment{abstract}
 {\small
  \begin{center}
  \bfseries \abstractname\vspace{-.5em}\vspace{0pt}
  \end{center}
  \list{}{%
    \setlength{\leftmargin}{10mm}
    \setlength{\rightmargin}{\leftmargin}%
  }%
  \item\relax}
 {\endlist}

\newcommand\extrafootertext[1]{%
    \bgroup
    \renewcommand\thefootnote{\fnsymbol{footnote}}%
    \renewcommand\thempfootnote{\fnsymbol{mpfootnote}}%
    \footnotetext[0]{#1}%
    \egroup
}

\title{Asymptotic analysis of ML-covariance parameter estimators based on covariance approximations}

\author{Reinhard Furrer}
\author{Michael Hediger}

\affil{Department of Mathematics, Winterthurerstrasse 190, 8057 Zurich, Switzerland}

\begin{document}

\maketitle

\extrafootertext{\textit{Email addresses}: \texttt{reinhard.furrer@math.uzh.ch}, \texttt{michael.hediger@math.uzh.ch}}

\begin{abstract}
Given a zero-mean Gaussian random field with a covariance function that belongs to a parametric family of covariance functions, we introduce a new notion of likelihood approximations, termed truncated-likelihood functions. Truncated-likelihood functions are based on direct functional approximations of the presumed family of covariance functions. For compactly supported covariance functions, within an increasing-domain asymptotic framework, we provide sufficient conditions under which consistency and asymptotic normality of estimators based on truncated-likelihood functions are preserved. We apply our result to the family of generalized Wendland covariance functions and discuss several examples of Wendland approximations. For families of covariance functions that are not compactly supported, we combine our results with the covariance tapering approach and show that ML estimators, based on truncated-tapered likelihood functions, asymptotically minimize the Kullback-Leibler divergence, when the taper range is fixed.\\[2mm]
\emph{Keywords:} Gaussian random fields, compactly supported covariance functions, likelihood approximations, consistency, asymptotic normality, covariance tapering.
\end{abstract}

\section{Introduction}

\subsection{\orange{On infill- and increasing-domain asymptotics}}

Maximum likelihood (ML) estimators for covariance parameters are highly popular in inference for random fields. Aiming towards asymptotic properties of such estimators, one needs to specify how the observation points and the associated sampling domain behave as the number of observation points increases. Two well-studied asymptotic frameworks are referred to as infill-domain asymptotics (also termed fixed-domain asymptotics) and increasing-domain asymptotics (see \cite{alma991014407499705251}, p.\ $100$ for an introduction of terms). In infill-domain asymptotics, observation points are sampled within a bounded sampling domain, whereas in increasing-domain asymptotics, the sampling domain grows as the number of observation points increases. When referring to infill- and increasing-domain asymptotics, one often places additional assumptions on the minimum distance between any two distinct observation points. In increasing-domain asymptotics, the latter distance is often assumed to be bounded away from zero, while in infill-domain asymptotics, one frequently assumes that distinct observation points can be sampled arbitrarily close to each other (see for example \cite{zhang2005towards}). There is a fair amount of literature which demonstrates that asymptotic properties of ML estimators for covariance parameters can be quite different under the two mentioned asymptotic frameworks (see \cite{zhang2005towards} or more lately \cite{bachoc2020asymptoticarXiv}). For example, it is known that some covariance parameters can not be estimated consistently under an infill-domain asymptotic framework (\cite{YING1991280}, \cite{zhang2004inconsistent}), whereas they can be estimated consistently, under given regularity conditions, within an increasing-domain asymptotic framework (\cite{mardia1984maximum}, \cite{bachoc2014asymptotic}). It is worth noting that in infill-domain asymptotics, these results can depend on the dimension $d$ of the Euclidean space $\mathbb{R}^{d}$, where the random field is assumed to be observed. For example, when the true covariance function belongs to the Mat\'{e}rn family (\cite{Matern}), and smoothness parameters are given, it is shown in \cite{zhang2004inconsistent}, that for $d = 1,2,3$, the scale and variance parameters can not be estimated consistently via an ML approach in an infill-domain asymptotic framework. The case where $d=4$ is still open, but for $d \geq 5$, it is shown in \cite{anderes2010consistent} that under infill-domain asymptotics, all covariance parameters of the Mat\'{e}rn family can be estimated consistently using an ML approach.

\subsection{\orange{Compactly supported covariance functions}}

In recent years, the dataset sizes have steadily increased such that statistical analyses on random fields can become quite expensive in terms of computational resources (see for example \cite{Flur:Furr:21} for a recent discussion). One prominent issue with large datasets is the large size of covariance matrices, constructed upon applying an underlying covariance function to given data. However, in certain fields of application, observed correlations are assumed to vanish beyond a certain cut-off distance (see \cite{https://doi.org/10.1002/qj.49712555417} p.\ $750-751$, and references therein, for an example in meteorology or also \cite{MR4439327} and \cite{Gerb:Moes:Furr:17}). On the other hand, in the context of real valued random fields, it is common practice to multiply a presumed covariance function with a known positive-definite and compactly supported covariance function, called the covariance taper. The resulting compactly supported covariance function is referred to as the tapered covariance function. For an introduction to covariance tapering we refer to \cite{MR2291261}. The use of compactly supported covariance functions can thus be of great importance for some fields of application. Not only do they potentially reflect the nature of the underlying covariance structure, but also, their application can lead to sparse covariance matrices. The latter are helpful in terms of the high computational costs in the context of large datasets. An excellent introduction to the construction of compactly supported covariance functions, associated to stationary and isotropic Gaussian random fields, is given in \cite{Gneiting2002}. Additional results are available in \cite{Zastavnyi2006}, \cite{b50d71882155d7fde924ff08f9ecdb1b} and \cite{CHERNIH201417}.

\subsection{\orange{Motivation}} \label{sec:motivation}
\orange{The parametric family of generalized Wendland covariance functions represents one example of a family of compactly supported covariance functions which allows, similar to the Mat\'{e}rn family, for a continuous parametrization of smoothness (in the mean square sense) of the underlying random field. Its origin is due to Wendland (\cite{wendland1995piecewise}) and an early adaptation for statistical applications was given by Gneiting (\cite{https://doi.org/10.1002/qj.49712555906}). In its general form (see \cite{Gneiting2002} and \cite{b50d71882155d7fde924ff08f9ecdb1b} for special cases) the generalized Wendland covariance function with smoothness parameters $\nu$ and $\kappa$, variance parameter $\sigma^{2}$ and range parameter $\beta$ is given by
\begin{equation} \label{StartWendland}
    \phi(t) \coloneqq \frac{\sigma^{2}}{\operatorname{B}(2\kappa,\nu + 1)\beta^{2\kappa + \nu}}\int_{t}^{\beta}w(w^{2} - t^{2})^{\kappa-1}(\beta-w)^{\nu}dw,
\end{equation}
if $t \in [0,\beta)$ and is zero otherwise. In the above display, $\operatorname{B}$ is the beta function. For technical details about valid parameter values, we refer to \cite{bevilacqua2019estimation} or Section~\ref{sec_application} of the present article. Clearly, in comparison with closed-form covariance functions, computing \eqref{StartWendland} is cumbersome, as it involves numerical integration. Depending on the support $\beta$ and a set of locations $s_{1}, \dotsc, s_{n} \in \mathbb{R}^{d}$, the $n \times n$ covariance matrix $\Sigma_{i,j} = \phi(\lVert s_{i} - s_{j} \rVert)$ requires at most $n(n-1)/2$ calculations of \eqref{StartWendland}. One strategy, which facilitates computing $\Sigma$, is to reduce the number of times \eqref{StartWendland} must be calculated. As an illustration, we give three examples which involve approximations $\tilde{\phi_{i}}$, $i = 1,2,3$, of $\phi$ (respectively approximations $\widetilde{\Sigma}_{i}$ of $\Sigma$):}  

\begin{enumerate}[label=(\alph*)]
    
    \item[\orange{($\tilde{\phi_{1}}$)}] \orange{Truncation of the support}

    \item[\orange{($\tilde{\phi_{2}}$)}] \orange{linear interpolation}

    \item[\orange{($\tilde{\phi_{3}}$)}] \orange{addition of a nugget effect }

\end{enumerate}

\orange{For $\tilde{\phi_{1}}$, we truncate $\phi$ to obtain $\tilde{\phi_{1}}$ which has a smaller support compared to $\phi$. This becomes especially interesting, when the original function $\phi$ tails off slowly (high degree of differentiability at the origin). As a result, $\widetilde{\Sigma}_{1}$ will be more sparse compared to $\Sigma$. Example $\tilde{\phi_{2}}$ is to predefine the numbers at which \eqref{StartWendland} is calculated. This is achieved by introducing a partition $0 < t_{1} < \dotsc < t_{N} = \beta$ of the support of $\phi$. Then, $\tilde{\phi_{2}}$ results in $N$ calculations of $\phi$. This defines a closed form approximation of $\phi$. Notice that $t_{1}, \dotsc, t_{N}$ do not need to be equispaced. Finally, $\tilde{\phi_{3}}$ can be interpreted as a tuning option for a given approximation $\tilde{\phi}_{*}$ of $\phi$:
\begin{equation*}
    \tilde{\phi_{3}}(t) \coloneqq \begin{cases}
    \tilde{\phi}_{*}(t) + \delta, & t=0, \\ \tilde{\phi}_{*}(t), & t \neq 0, \end{cases} \quad \delta \geq 0.
\end{equation*}
With regard to practical usage, this form of approximation increases numerical stability. Further, it allows for more flexibility in practice, where the number of observations $n$ is given and $\widetilde{\Sigma}_{*}$ based on $\tilde{\phi}_{*}$ might not be positive-definite.}

\orange{Following up the above examples, we picture an approximation $\tilde{\phi}$ of $\phi$ (respectively approximation $\widetilde{\Sigma}$ of $\Sigma$). Several questions arise:}

\begin{itemize}
    
    \item \orange{What are conditions on $\tilde{\phi}$ to ensure that $\widetilde{\Sigma}$ is asymptotically (as $n \xrightarrow[]{} \infty$) equivalent to $\Sigma$ and eventually (for $n$ large enough) remains positive-definite? }

    \item \orange{In terms of ML estimators for covariance parameters, how shall a log-likelihood approximation based on $\tilde{\phi}$ be defined?}

    \item \orange{Under which conditions on $\tilde{\phi}$ are ML estimators based on $\tilde{\phi}$ consistent and asymptotically normal?}

\end{itemize}

\orange{In the more general setting of a given parametric family of covariance functions, the present study gives a concrete context, where the latter questions are answered by introducing the notion of truncated-ML estimators.}

\subsection{\orange{Framework and contribution}}

Truncated-ML estimators for covariance parameters are based on truncated-likelihood functions. The latter are defined upon parametric families of sequences of functions, which approximate a presumed family of covariance functions on a common domain. Colloquially we will call these parametric sequences of functions covariance approximations. The respective matrices, constructed upon applying covariance approximations to a given collection of observation points, will be termed covariance matrix approximations. We will allow for covariance matrix approximations that are not necessarily positive semi-definite. Therefore, truncated-likelihood functions are more general than existing likelihood approximations methods such as low-rank, Vecchia, or covariance tapering approaches (see \cite{MR3996451} for a summary of commonly used methods). 

We work in an increasing-domain asymptotic framework, where collections of observation points are realizations of finite collections of a randomly perturbed regular grid (see also \cite{bachoc2014asymptotic}). We consider a stationary Gaussian random field, with a zero-mean function and a true unknown covariance function that belongs to a given parametric family of covariance functions. If the presumed family of covariance functions is compactly supported, we provide sufficient conditions under which truncated-ML estimators and (regular) ML estimators for covariance parameters are consistent and asymptotically normal. Some conditions imposed on families of covariance functions are identical to the conditions that were already considered in \cite{bachoc2014asymptotic}. The main difference is that we work with compactly supported covariance functions. Therefore, it is possible to simplify some of the conditions that were set up in \cite{bachoc2014asymptotic}. As for statistical applications, we apply these results to the family of generalized Wendland covariance functions. In contrast to the infill-domain asymptotic framework considered in \cite{bevilacqua2019estimation}, we show that under the studied increasing-domain asymptotic framework, under some conditions on the parameter space, (regular) ML estimators for variance and range parameters are consistent and asymptotically normal. Further, we show that the same asymptotic results are recovered for truncated-ML estimators, based on various generalized Wendland approximations, such as truncations, linear interpolations and added nugget effects. 

Additionally, we provide an extension to families of covariance functions which are not compactly supported. We combine our results with the covariance tapering approach. That is, we study covariance taper approximations and their asymptotic influence on the conditional Kullback-Leibler divergence of the misspecified distribution from the true distribution (see also \cite{bachoc2018asymptotic}). We show that the latter divergence is minimized by truncated-tapered ML estimators.

\subsection{\orange{Structure of the article}}

The rest of the article is organized as follows. Section~\ref{sec:context} establishes the context. We introduce some primary notation, define the sampling domain and the random field itself. In Section~\ref{sec:regularity} we introduce regularity conditions on covariance functions and approximations. In Section~\ref{sec:uniformequivof mat} we present intermediate asymptotic results on covariance matrices and approximations. Section~\ref{MLE_results} contains our main results: We introduce truncated-ML estimators and present results on consistency and asymptotic normality. In Section~\ref{sec_application}, we apply our results to the family of generalized Wendland covariance functions and discuss several examples of generalized Wendland approximations. Then, in the context of non-compactly supported covariance functions, Section~\ref{sec:app2} contains results on the asymptotic influence of taper approximations on the Kullback-Leibler divergence. Section~\ref{sec:conc} gives an outlook and some final comments. The \hyperref[appn]{Appendix} is split into three parts. Covariance approximations for isotropic random fields are discussed in Appendix~\ref{sec_application1}. Appendix~\ref{appn} contains additional supporting results, whereas all the proofs are left for Appendix~\ref{appnB}.

\section{Context} \label{sec:context}

\subsection{Primary notation}
The set $\natnum$ and $\mathbb{R}_{+}$ shall represent the set of positive integers and non-negative real numbers, respectively. For $d \in \natnum$, we use the notation $\mathrm{B}(x;r)$ ($\mathrm{B}[x;r]$) for the open (closed) ball of radius $r>0$ with center $x \in \mathbb{R}^{d}$. Given $n \in \natnum$, for some set $A \subset \mathbb{R}^{n}$, we write $\mathfrak{B}(A)$ for the Borel $\sigma$-algebra on $A$. 

For a vector $(w_{1}, \dotsc, w_{d}) = w \in \mathbb{R}^{d}$, we write $\nnorm{w} =\big(w_{1}^{2} + \cdots + w_{d}^{2}\big)^{1/2}$ for the Euclidean norm of $w$ on $\mathbb{R}^{d}$. In the case of $d=1$ we use the notation $\norm{\cdot}$ for the Euclidean norm. For two vectors $w,w^{\prime} \in \mathbb{R}^{d}$, $\langle w, w^{\prime} \rangle = w^{\mathrm{t}}w^{\prime} = \sum_{i = 1}^{d}w_{i}w^{\prime}_{i}$ represents the inner product that induces $\nnorm{\cdot}$ on $\mathbb{R}^{d}$. Given $D \subset \mathbb{R}^{d}$, we write $\mathcal{B}_{\text{C}}(D;S)$ for the space of real valued, uniformly bounded functions on $D$, having compact support $S \subset D$. If $f \in \mathcal{B}_{\text{C}}(D;S)$ and $f$ is also continuous, we use the notation $\mathcal{C}_{\text{C}}(D;S)$ instead of $\mathcal{B}_{\text{C}}(D;S)$. For $f \in \mathcal{C}_{\text{C}}(D;S)$ we write $\supnorm{f} = \sup\{\norm{f(h)} \colon h \in D\}$ for the uniform norm on $\mathcal{C}_{\text{C}}(D;S)$. For vectors $w \in \mathbb{R}^{d}$, $\norm{w}_{\infty} = \max_{i = 1, \dotsc, d}\norm{w_{i}}$ denotes the uniform norm on $\mathbb{R}^{d}$. 

For a real $n \times n$ matrix $A$, $\specnorm{A} = \max_{\{z \colon z^{\mathrm{t}}z=1\}}\langle z, A^{\mathrm{t}}A z\rangle^{1/2}$ denotes the spectral norm of $A$. We write $A \succ 0$ ($A \prec 0$) to indicate that $A$ is positive-definite (negative-definite). Further, $\lambda_{1}(A) \geq \cdots \geq \lambda_{n}(A)$ denote the $n$ real eigenvalues of a matrix $A \in \mathrm{S}_{n\times n}(\mathbb{R})$, where $\mathrm{S}_{n\times n}(\mathbb{R})$ represents the space of real symmetric $n \times n$ matrices. 

We use the notation $\nabla f(x) = \big(\frac{\partial f}{\partial x_{1}}(x), \dotsc, \frac{\partial f}{\partial x_{p}}(x)\big)$ for the gradient of $f$ at $x$, where $x \mapsto f(x)$ is any differentiable, real valued function, defined on some $E \subset \mathbb{R}^{p}$. Further, for a vector valued, differentiable function $g(x) = (g_{1}(x), \dotsc, g_{m}(x))$, with values in $\mathbb{R}^{m}$, defined on some $U \subset \mathbb{R}^{p}$, we write $J_{g}(x)_{l,k} = \frac{\partial g_{l}}{\partial x_{k}}(x)$, $1 \leq l \leq m$, $1 \leq k \leq p$, for the Jacobi-matrix of $g$ at~$x$. 

A mapping $Y$ from a probability space $(\Omega, \mathcal{F}, \mathbb{P})$ to a measure space $(E, \mathcal{A})$ will be called a random element if it is $\mathcal{F}/\mathcal{A}$ measurable. If we write that $Y \colon (\Omega, \mathcal{F}) \to (E, \mathcal{A})$ is measurable, we mean that it is $\mathcal{F}/\mathcal{A}$ measurable.  If $(Y_{n}){}_{n \in \natnum}$ denotes a sequence of random elements, where for any $n \in \natnum$, $Y_{n}$ is a mapping form a probability space $(\Omega, \mathcal{F}, \mathbb{P})$ to a measure space $(E, \mathcal{A})$, we use the notation 
\begin{gather*}
    Y_{n} \toinp Y \text{ and }Y_{n} \toind \mathcal{L},
\end{gather*}
to indicate convergence of $(Y_{n}){}_{n \in \natnum}$ to a random element $Y$ in probability and in distribution, respectively. Note that for convergence in distribution, the introduced notation indicates that the limit $Y$ has law $\mathcal{L}$ on $(E, \mathcal{A})$. \orange{A sequence of estimators $(\hat{\theta}_{n})_{n \in \natnum}$ for $\theta_{0} \in \mathbb{R}^{p}$ will be referred to as consistent if it converges in probability to $\theta_{0}$. Finally, $\mathcal{N}(\mu, \Sigma)$ indicates a multivariate normal distribution with mean vector $\mu$ and covariance matrix $\Sigma$.}   

\subsection{Random sampling scheme} \label{sec:context1}

On a probability space $(\Omega, \mathcal{F}, \mathbb{P})$, we consider a real valued Gaussian random function $Z$, which has sample functions on $\mathbb{R}^{d}$. We assume that $Z$ is stationary (homogeneous) with zero-mean function and covariance function $c_{\theta_{0}}(s)$, $s \in \mathbb{R}^{d}$, where $\theta_{0} \in \Theta$, with $\Theta \subset \mathbb{R}^{p}$, compact and convex. Thus, we consider a real valued random field $\{Z_{s} \colon s \in \mathbb{R}^{d}\}$, which has true and unknown covariance function $c_{\theta_{0}}$ that belongs to a family of covariance functions $\{c_{\theta} \colon \theta \in \Theta\}$.

\begin{figure}[!t]
    \centering
\includegraphics[width=0.9\textwidth]{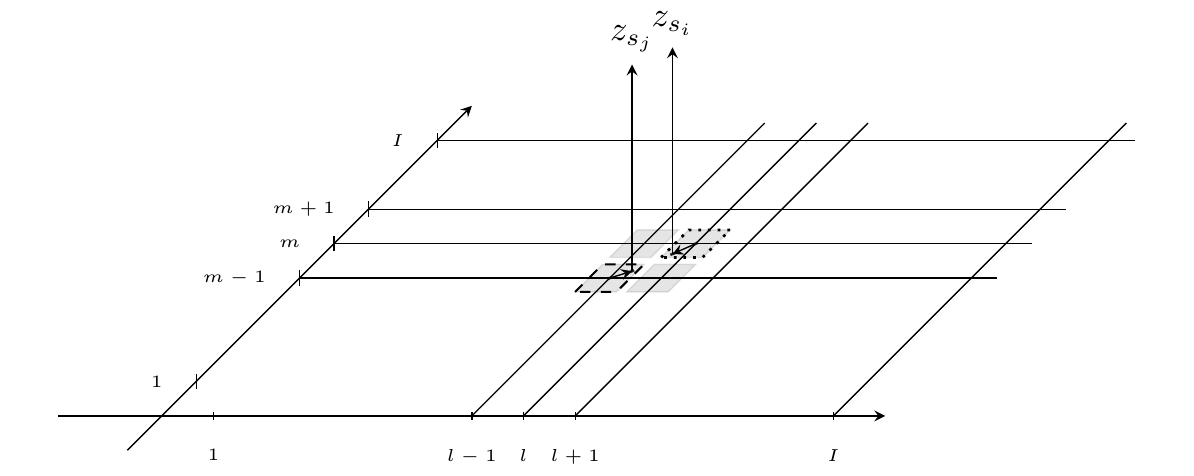}
\caption{For $\tau = 0.4$ and $I \in \natnum$, a random field $Z$ is observed at two realizations $s_{i}$ and $s_{j}$ of $S_{i} = v_{i} + \tau X_{i}$ and $S_{j} = v_{j} + \tau X_{j}$, $i \neq j$, $1 \leq i,j \leq I^{2}$. Dotted and dashed lines mark the borders of the ranges of $S_{i}$ and $S_{j}$, respectively.}
\label{fig:RandomPerturbedRegularGrid}
\end{figure}

Let $\mathcal{Q} \coloneqq [-1,1]^{d}$ and $X \colon \Omega \to \mathcal{Q}^{\natnum}$ be a stochastic process, defined on the same probability space $(\Omega, \mathcal{F}, \mathbb{P})$, but independent of $Z$. We assume that the sequence $(X_{i})_{i \in \natnum}$ is a sequence of independent random vectors with common law on $\mathcal{Q}$, which has a  strictly positive probability density function on $\mathcal{Q}$ (see also Remark~\ref{remark:conditionaldisti}). Given $\tau \in [0,1/2)$ and a sequence of deterministic points $(v_{i}){}_{i \in \natnum}$, with $v_{i} \in \natnum^{d}$, we define a randomly perturbed regular grid $S$ as the process
\begin{equation} \label{randomgrid}
    \big\{S_{i} \coloneqq v_{i} + \tau X_{i} \colon i \in \natnum\big\},
\end{equation}
where we assume that for all $I \in \natnum$, $\big\{v_{i}, 1 \leq i \leq I^{d}\big\} = \big\{1, \dotsc, I\big\}^{d}$. Therefore, for any $\omega \in \Omega$, $S(\omega)$ is a sequence on $\natnum$, with image $S[\natnum](\omega) \subset \prod_{i=1}^{\infty}\big(v_{i} + \tau \mathcal{Q}\big) \eqqcolon \mathcal{G}$ and any first $I^{d}$ coordinates are in $\{1,\dotsc,I\}^{d} + \tau \mathcal{Q}$ (see also Figure~\ref{fig:RandomPerturbedRegularGrid}). At this point we remark that if nothing is mentioned, the parameter $\tau \in [0,1/2)$ and the sequence $(v_{i}){}_{i \in \natnum}$ shall be fixed. Let $X_{(n)} \coloneqq (X_{1}, \dotsc, X_{n})$ and $S_{(n)} \coloneqq (S_{1}, \dotsc, S_{n})$ denote finite collections of $X$ and $S$, respectively. We use the notation $x_{(n)} \coloneqq (x_{1}, \dotsc, x_{n})$ for a vector that contains the first $n$ entries of a given sequence in $\mathcal{Q}^{\natnum}$. Correspondingly, given $\tau \in [0,1/2)$, $v_{1}, \dotsc, v_{n}$ and $x_{(n)} \in \mathcal{Q}^{n}$, we write $s_{(n)} \coloneqq (s_{1}, \dotsc, s_{n})$, $s_{i} = v_{i} + \tau x_{i}$, for $n$ perturbed grid locations in $\mathcal{G}_{n} \coloneqq \prod_{i=1}^{n}(v_{i} + \tau\mathcal{Q})$. 

On $(\Omega, \mathcal{F}, \mathbb{P})$, we define the random vector 
\begin{align} \label{composition}
    \omega \mapsto Z_{(n)}(\omega) \coloneqq \big(Z_{S_{1}(\omega)}(\omega), \dotsc, Z_{S_{n}(\omega)}(\omega)\big) = (z_{s_{1}}, \dotsc, z_{s_{n}})\eqqcolon z_{(n)},
\end{align}
which denotes $Z$ observed at a finite collection of $S$. The situation, where a Gaussian random field is assumed to be observed at a randomly perturbed regular grid, with parameter $\tau$ and deterministic points $(v_{i}){}_{i \in \natnum}$, as introduced above, is also considered in \cite{bachoc2014asymptotic}. 
Given $\theta \in \Theta$, we let $\rcovmat \coloneqq [c_{\theta}(s_{i}- s_{j})]_{1 \leq i,j \leq n}$ denote the non-random $n \times n$ covariance matrix based on an arbitrary $s_{(n)} \in \mathcal{G}_{n}$. On $(\Omega, \mathcal{F}, \mathbb{P})$, we write
\begin{gather*} 
    \omega \mapsto \covmat(\omega) \coloneqq \rrcovmat, \quad \theta \in \Theta,
\end{gather*}
for the $n \times n$ random covariance matrix based on a finite collection $S_{(n)}$ of $S$.

\begin{remark} \label{remark:conditionaldisti}
\normalfont Some technical remarks are worth pointing out. \orange{We assume that the random function $Z(s,\omega) \coloneqq Z_{s}(\omega)$, $s \in \mathbb{R}^{d}$, is measurable as a function from the measure space $(\mathbb{R}^{d} \times \Omega, \mathfrak{B}(\mathbb{R}^{d}) \otimes \mathcal{F})$ to $(\mathbb{R}, \mathfrak{B}(\mathbb{R}))$. That is to say that $Z$ is (jointly) measurable. This condition makes sure that the components $\omega \mapsto Z_{S_{i}(\omega)}(\omega) = Z(S_{i}(\omega),\omega)$, $i = 1, \dotsc, n$, of \eqref{composition} are $\mathcal{F}/\mathfrak{B}(\mathbb{R})$ measurable as the composition of the measurable functions $\omega \mapsto (S_{i}(\omega), \omega)$ and $(s,\omega) \mapsto Z(s,\omega)$. Thus, the random vector $Z_{(n)}$ is well defined}. Since $Z$ and $S$ are independent, it is readily seen that the conditional distribution of $Z_{(n)}$ given $S_{(n)} = s_{(n)}$ is Gaussian, with characteristic function $\operatorname{exp}(-(1/2)a^{\mathrm{t}}\tcovmat(\omega)a)$, $a \in \mathbb{R}^{n}$. In addition, we note that for fixed $\omega \in \Omega$, $S[\natnum](\omega)$ is not bounded and if we define $\Delta_{\tau} \coloneqq 1-2\tau$, we are given some fixed $\Delta_{\tau} > 0$, which is independent of $n \in \natnum$ and $\theta \in \Theta$, such that
\begin{equation} \label{IncreasingDomain}
    \inf_{n\in \natnum}\inf_{\substack{1 \leq i, j \leq n \\ i \neq j}}\nnorm{s_{i} - s_{j}} \geq \Delta_{\tau}.
\end{equation} 
Hence, we are in an increasing-domain asymptotic framework where the minimum distance between any two distinct observation points is bounded away from zero. The assumption that for any given $i \in \natnum$, $X_{i}$ has strictly positive probability density function on $\mathcal{Q}$, is purely technical (see also the proof of Theorem~\ref{thm1}). As it can be seen from the mentioned proof, if $\tau = 0$, the assumption becomes redundant.   
\end{remark}

\section{Regularity Conditions on covariance functions and covariance approximations} \label{sec:regularity}

\subsection{Regularity conditions on the family of covariance functions}

\begin{assumption}[Regularity conditions on $c_{\theta}$]\label{cov_a2} 
~\begin{enumerate}[label=(\arabic*)]

\item \label{3.2.1} There exist real constants $C$, $L < \infty$, which are independent of $\theta \in \Theta$, such that $c_{\theta} \in \mathcal{B}_{\text{C}}(\mathbb{R}^{d};S_{\theta})$, with $S_{\theta} \subset \mathrm{B}[0;C]$ and $\supnorm{c_{\theta}} \leq L$.  

\item \label{3.2.2} For any $s \in \mathbb{R}^{d}$, the first, second and third order partial derivatives of $\theta \mapsto c_{\theta}(s)$ exist. In addition, for any $q = 1,2,3$, $i_{1}, \dotsc, i_{q} \in \{1, \dotsc,p\}$, 
$\frac{\partial^{q} c_{\theta}}{\partial \theta_{i_1}\cdots\partial \theta_{i_q}} \in \mathcal{B}_{\text{C}}\big(\mathbb{R}^{d};S_{\theta}(i_{1}, \dotsc, i_{q})\big)$, where there exist constants $C^{\prime} \text{, }L^{\prime} < \infty$, which are independent of $\theta \in \Theta$, such that $S_{\theta}(i_{1}, \dotsc, i_{q}) \subset \mathrm{B}[0;C^{\prime}]$ and $\big\lVert\frac{\partial^{q} c_{\theta}}{\partial \theta_{i_1}\cdots\partial \theta_{i_q}}\big\rVert_{\infty} \leq~L^{\prime}$.   
 
\item \label{3.2.3} Fourier inversion holds, that is for any $\theta \in \Theta$
\begin{gather*}
    c_{\theta}(s) = \int_{\mathbb{R}^{d}}\!\hat{c}_{\theta}(f)\eeexp^{\ii \langle f, s \rangle}\mathrm{d}\!f,
\end{gather*}
with $\Theta \times \mathbb{R}^{d} \ni (\theta, f) \mapsto \hat{c}_{\theta}(f)$ continuous and strictly positive. 
\end{enumerate}
\end{assumption}

\begin{remark} \label{remark31}
Note that \ref{3.2.1} and \ref{3.2.2} of Assumption~\ref{cov_a2} are different to the conditions assumed in \cite{bachoc2014asymptotic} (compare also to Condition $3.2$ imposed in \cite{bachoc2018asymptotic}, or Condition~$4$ stated in \cite{bachoc2020asymptotic}). In \cite{bachoc2014asymptotic} it is assumed that a given covariance function $k_{\theta}$ is not only bounded on $\mathbb{R}^{d}$, but also it decays sufficiently fast in the Euclidean norm on $\mathbb{R}^{d}$. Explicitly, it is assumed in Condition~2.1 of \cite{bachoc2014asymptotic} that there exists a finite constant $A$, which is independent of $\theta \in \Theta$, such that for any $s \in \mathbb{R}^{d}$, $\norm{k_{\theta}(s)} \leq A/(1+\lVert s \rVert^{d+1})$. This polynomial decay condition on $k_{\theta}$ can be interpreted as a summability condition on the entries of the respective covariance matrices $K_{\theta}(s_{(n)})_{i,j} \coloneqq k_{\theta}(s_{i} - s_{j})$, which guaranties that the maximal eigenvalues of $K_{\theta}(s_{(n)})$ are uniformly bounded in $n \in \natnum$, $s_{(n)}\in \mathcal{G}_{n}$ and $\theta \in \Theta$ (see Lemmas D.1 and D.5 in \cite{bachoc2014asymptotic}). Note that the exponent $d+1$ can be replaced by $d + \alpha$, with $\alpha > 0$ some fixed constant (see also (6) in \cite{bachoc2020asymptoticarXiv}). In the present study we show that under the assumption of a minimal spacing between any two distinct observation points, if $c_{\theta}$ has compact support on $\mathbb{R}^{d}$, the number of possible observation points, which are covered by the support of $c_{\theta}$, must be bounded uniformly in $n \in \natnum$, $s_{(n)}\in \mathcal{G}_{n}$ and $\theta \in \Theta$ (see Lemma~\ref{lemma1}). This, together with the condition that $c_{\theta}$ is also uniformly bounded on $\Theta$ and $\mathbb{R}^{d}$, will be sufficient to conclude that the maximal eigenvalues of $\rcovmat$ are uniformly bounded in $n \in \natnum$, $s_{(n)}\in \mathcal{G}_{n}$ and $\theta \in \Theta$ (see Lemmas~\ref{eigenbehave} and~\ref{eigenbehave1x1}). Similar remarks can be made with regard to the conditions imposed on the partial derivatives of $c_{\theta}$ with respect to $\theta$ (see Lemma~\ref{lemmaEigenbehavePartials}). In addition, \ref{3.2.3} of Assumption~\ref{cov_a2} is also imposed in \cite{bachoc2014asymptotic} (compare also to \cite{Furrer_Bachoc_2016_smallest_eigenvalue} and \cite{bachoc2020asymptotic}). It guarantees that the minimal eigenvalues of $\rcovmat$ are bounded from below, uniformly in $n \in \natnum$, $s_{(n)}\in \mathcal{G}_{n}$ and $\theta \in \Theta$ (see Lemmas~\ref{eigenbehave} and~\ref{eigenbehave1x1}). Finally, we remark that within the framework of compactly supported covariance functions, the given conditions are very minimal and can be considered as classical in the context of ML estimation. Especially, if one is not interested in the asymptotic distribution, and rather seeks conditions under which ML estimators are consistent (with regard to a concrete example, we refer to Remark~\ref{remark:kappa}).             
\end{remark}

\subsection{Regularity conditions on the family covariance approximations}  \label{sec_cov_approx}

Given $\theta \in \Theta$, we let $(\tilde{c}_{m,\theta})_{m \in \natnum}$ denote a sequence of real valued functions defined on $\mathbb{R}^{d}$. The families $\big\{(\tilde{c}_{m,\theta})_{m \in \natnum} \colon \theta \in \Theta\big\}$ can be put under the following assumption.  

\begin{assumption}[Regularity conditions on $\tilde{c}_{m,\theta}$] \label{uniformly_approx_covfun} 
~\begin{enumerate}[label=(\arabic*)]

\item \label{3.4.0} For any $\theta \in \Theta$ and $m \in \natnum$, the function $\tilde{c}_{m,\theta}\colon \big(\mathbb{R}^{d}, \mathfrak{B}(\mathbb{R}^{d})\big) \to \big(\mathbb{R}, \mathfrak{B}(\mathbb{R})\big)$ is measurable and such that $\tilde{c}_{m,\theta}(s) = \tilde{c}_{m,\theta}(-s)$ for any $s \in \mathbb{R}^{d}$.

\item \label{3.4.1} For any $m \in \natnum$, $\tilde{c}_{m,\theta}$ satisfies \ref{3.2.1} of Assumption~\ref{cov_a2}, where respective constants $\widetilde{C}$ and $\widetilde{L}$ can be further chosen independently of $m \in \natnum$.

\item \label{3.4.2} $\sup_{\theta \in \Theta}\supnorm{\tilde{c}_{m,\theta}-c_{\theta}} \tonormalm 0$.

\item \label{3.4.3} For any $m \in \natnum$, $\tilde{c}_{m,\theta}$ satisfies \ref{3.2.2} of Assumption~\ref{cov_a2}, where respective constants $\widetilde{C}^{\prime}$ and $\widetilde{L}^{\prime}$ can be further chosen independently of $m \in \natnum$.   

\item \label{3.4.4} For any $q = 1,2,3$, $i_{1}, \dotsc, i_{q} \in \{1, \dotsc,p\}$, we have that
\begin{gather*}
    \sup_{\theta \in \Theta}\supnorm{\frac{\partial^{q} \tilde{c}_{m,\theta}}{\partial \theta_{i_1}\cdots\partial \theta_{i_q}} - \frac{\partial^{q} c_{\theta}}{\partial \theta_{i_1}\cdots\partial \theta_{i_q}}} \tonormalm 0.
\end{gather*}
\end{enumerate}

\end{assumption}

To make the notation easier, we write $(\tilde{c}_{m,\theta}) \coloneqq (\tilde{c}_{m,\theta}){}_{m \in \natnum}$. In the following, we formally introduce covariance matrix approximations (random and non-random versions). To do so, let $r \colon \natnum \to \natnum$ be such that $r(n) \xrightarrow{} \infty$ as $n \xrightarrow{} \infty$. Given $s_{(n)} \in \mathcal{G}_{n}$, we let $\racovmat \coloneqq [\tilde{c}_{r(n), \theta}(s_{i}- s_{j})]_{1 \leq i,j \leq n}$ denote the non-random $n \times n$ matrix based on a given family $\big\{(\tilde{c}_{m,\theta}) \colon \theta \in \Theta\big\}$. Then, on $(\Omega, \mathcal{F}, \mathbb{P})$, if $\big\{(\tilde{c}_{m,\theta}) \colon \theta \in \Theta\big\}$ is a family of Borel measurable sequences of functions, we write
\begin{gather*} 
    \omega \mapsto \acovmat(\omega)  \coloneqq \rracovmat,
\end{gather*}
for the $n \times n$ random matrix based on a finite collection $S_{(n)}$ of $S$. Colloquially we will use the term covariance approximation when we refer to a given family $\big\{(\tilde{c}_{m,\theta}) \colon \theta \in \Theta\big\}$, which can approximate a family of covariance functions $\{c_{\theta} \colon \theta \in \Theta\}$ in the sense of Assumption~\ref{uniformly_approx_covfun}. In these terms $\{c_{\theta} \colon \theta \in \Theta\}$ itself is a covariance approximation. The expression covariance matrix approximation will be used for both, $\racovmat$ and its random version $\acovmat$. Similar, we use the expression covariance matrix for both, $\rcovmat$ and $\covmat$.

\begin{remark} \label{remark:unifassump}
\normalfont \ref{3.4.0}, \ref{3.4.1} and \ref{3.4.3} of Assumption~\ref{uniformly_approx_covfun} are natural extensions of \ref{3.2.1} and \ref{3.2.2} of Assumption~\ref{cov_a2}. Notice that the measurability condition imposed in \ref{3.4.0} of Assumption~\ref{uniformly_approx_covfun} makes sure that $\acovmat$ is $\mathcal{F}/\mathfrak{B}(\mathbb{R}^{n^2})$ measurable. Condition \ref{3.4.2} of Assumption~\ref{uniformly_approx_covfun} specifies in which sense a family $\big\{(\tilde{c}_{m,\theta}) \colon \theta \in \Theta\big\}$ approximates the family $\{c_{\theta} \colon \theta \in \Theta\}$. We require that $(\tilde{c}_{m,\theta})$ converges uniformly on $\mathbb{R}^{d}$ to $c_{\theta}$, where the convergence is also uniform on the parameter space $\Theta$. In fact, we will show (see Lemmas~\ref{eigenbehave1x1} and~\ref{eigenbehave}) that the uniform convergence of $(\tilde{c}_{m,\theta})$ to $c_{\theta}$, together with the condition that the families $\{c_{\theta} \colon \theta \in \Theta\}$ and $\big\{(\tilde{c}_{m,\theta}) \colon \theta \in \Theta\big\}$ have uniformly bounded compact support, are, among others, sufficient criteria to proof that the matrices $\rcovmat$ and $\racovmat$ are asymptotically (as $n \xrightarrow[]{} \infty$) equivalent, uniformly on $\Theta$ and $\mathcal{G}$. Condition \ref{3.4.4} of Assumption~\ref{uniformly_approx_covfun} will allow us to conclude that a similar result holds true for the first, second and third order partial derivatives (with respect to $\theta$) of $\racovmat$ and $\rcovmat$. For concrete examples of covariance approximations, where the conditions of Assumption~\ref{uniformly_approx_covfun} are verified, we refer to Section~\ref{sec_application}.   
\end{remark}

\section{Uniform asymptotic equivalence of covariance matrices and covariance matrix approximations} \label{sec:uniformequivof mat}

This section presents intermediate results on covariance matrices and approximations. In particular, Lemma~\ref{eigenbehave} gives precise conditions under which $\acovmat$ eventually (for $n$ large enough) remains positive-definite with $\mathbb{P}$ probability one. 

\begin{lemma} \label{eigenbehave}
Assume that the family $\{c_{\theta} \colon \theta \in \Theta\}$ satisfies \ref{3.2.1} and \ref{3.2.3} of Assumption~\ref{cov_a2}. Consider $\big\{(\tilde{c}_{m,\theta}) \colon \theta \in \Theta\big\}$ that satisfies \ref{3.4.0}, \ref{3.4.1} and \ref{3.4.2} of Assumption~\ref{uniformly_approx_covfun}.
Then, we have that $\as$
\begin{equation*}
        \sup_{n\in \natnum}\sup_{\theta \in \Theta}\big\lVert\covmat\big\rVert_{2} < \infty\text{ and }\sup_{n\in \natnum}\sup_{\theta \in \Theta}\big\lVert\acovmat\big\rVert_{2}  < \infty.
\end{equation*}
In particular we can conclude that $\as$
\begin{equation*} 
        \sup_{\theta \in \Theta}\big\lVert\covmat - \acovmat\big\rVert_{2} \tonormal 0.
\end{equation*}
Further, it is true that $\as$
\begin{equation*}
       \inf_{n\in \natnum}\inf_{\theta \in \Theta}\lambda_{n}\big(\covmat\big) > 0, 
\end{equation*}
and there exists $N \in \natnum$ such that $\as$ 
\begin{equation*}
       \inf_{n\geq N}\inf_{\theta \in \Theta}\lambda_{n}\big(\acovmat\big) > 0.
\end{equation*}
\end{lemma}

\section{Truncated-ML estimators} \label{MLE_results}

\orange{Given a square matrix $A$, we define $\operatorname{det}_{+}(A)$ to be the product of the strictly positive eigenvalues of $A$. If all of the eigenvalues are less or equal to zero, $\operatorname{det}_{+}(A) = 1$. Further, we use the notation $A^{+}$ for the pseudoinverse of $A$ (sometimes called Moore-Penrose inverse)}. For the given collection $\{c_{\theta} \colon \theta \in \Theta\}$, we define, on $(\Omega, \mathcal{F}, \mathbb{P})$, for any $n \in \natnum$ and $\theta \in \Theta$, the random variable
\orange{\begin{align} \label{pseudo1}
   l_{n}(\theta) &\coloneqq \frac{1}{n}\operatorname{log}\big(\operatorname{det}_{+}(\covmat)\big) + \frac{1}{n}\big\langle Z_{(n)}, \covmat^{+}Z_{(n)}\big\rangle.
\end{align}
}Given $\omega \in \Omega$, $\theta \mapsto l_{n}(\theta)(\omega)$ shall be called the \textit{truncated-modified log-likelihood function} based on $\{c_{\theta} \colon \theta \in \Theta\}$. A sequence of estimators $\big(\orange{\hat{\theta}_{n}(c)}\big){}_{n \in \natnum}$, defined on $(\Omega, \mathcal{F}, \mathbb{P})$, will be called a sequence of \textit{truncated-ML estimators} for $\theta_{0}$ based on $\{c_{\theta} \colon \theta \in \Theta\}$, if for any $n \in \natnum$,
\begin{gather*}
    \hat{\theta}_{n}(c) \in \argmin_{\theta\in \Theta}{l}_{n}(\theta).
\end{gather*}  
Similarly, on $(\Omega, \mathcal{F}, \mathbb{P})$, for a given collection of sequences of real valued functions $\big\{(\tilde{c}_{m,\theta}) \colon \theta \in \Theta\big\}$, we introduce, for any $n \in \natnum$ and $\theta \in \Theta$, the random variable
\begin{align} \label{pseudo2}
   \tilde{l}_{n}(\theta) &\coloneqq \frac{1}{n}\operatorname{log}\big(\operatorname{det}_{+}(\acovmat)\big) + \frac{1}{n}\big\langle Z_{(n)}, \acovmat^{+}Z_{(n)}\big\rangle.
\end{align}
Then, for $\omega \in \Omega$, the function $\theta \mapsto \tilde{l}_{n}(\theta)(\omega)$ denotes the truncated-modified log-likelihood function based on $\big\{(\tilde{c}_{m,\theta}) \colon \theta \in \Theta\big\}$. A sequence of estimators $\big(\hat{\theta}_{n}(\tilde{c})\big){}_{n \in \natnum}$, defined on $(\Omega, \mathcal{F}, \mathbb{P})$, will be called a sequence of truncated-ML estimators for $\theta_{0}$ based on $\big\{(\tilde{c}_{m,\theta}) \colon \theta \in \Theta\big\}$, if for any $n \in \natnum$
\begin{equation} \label{pseudocriterion}
   \hat{\theta}_{n}(\tilde{c}) \in \argmin_{\theta\in \Theta}{\tilde{l}}_{n}(\theta).
\end{equation}   
At this point is important to note that for a given $\omega \in \Omega$, it is in general not true that $l_{n}(\theta)(\omega)$ and $\tilde{l}_{n}(\theta)(\omega)$ are continuous in $\theta$ for any $n \in \natnum$. Nevertheless, a consequence of Lemma~\ref{eigenbehave} is the following proposition: 

\begin{proposition} \label{orderN}
Assume that the family $\{c_{\theta} \colon \theta \in \Theta\}$ satisfies \ref{3.2.1} and \ref{3.2.3} of Assumption~\ref{cov_a2}. Consider $\big\{(\tilde{c}_{m,\theta}) \colon \theta \in \Theta\big\}$ that satisfies \ref{3.4.0}, \ref{3.4.1} and \ref{3.4.2} of Assumption~\ref{uniformly_approx_covfun}. Then, we have that for any $n \in \natnum$, $\as$,
\begin{gather*}
    l_{n}(\theta) = \frac{1}{n}\operatorname{log}\big(\operatorname{det}(\covmat)\big) + \frac{1}{n}\big\langle Z_{(n)}, \covmat^{-1}Z_{(n)}\big\rangle.
\end{gather*}
Further there exists $N \in \natnum$ such that for any $n \geq N$, $\as$,
\begin{gather*}
    \tilde{l}_{n}(\theta) = \frac{1}{n}\operatorname{log}\big(\operatorname{det}(\acovmat)\big) + \frac{1}{n}\big\langle Z_{(n)}, \acovmat^{-1}Z_{(n)}\big\rangle,
\end{gather*}
and we have that
\begin{gather*}
    \sup_{\theta \in \Theta}\big \lvert l_{n}(\theta) - \tilde{l}_{n}(\theta) \big \rvert \toinp 0.
\end{gather*}
\end{proposition}
Using Proposition~\ref{orderN}, we notice that if, for any $s \in \mathbb{R}^{d}$ and $m \in \natnum$, both $\theta \mapsto c_{\theta}(s)$ and $\theta \mapsto \tilde{c}_{m,\theta}(s)$ are $k$ times differentiable, we have that $\theta \mapsto l_{n}(\theta)(\omega)$ and $\theta \mapsto\tilde{l}_{n}(\theta)(\omega)$ are $k$ times differentiable for $n$ large enough, respectively. 

For the rest of the article, if we refer to truncated-ML estimators (without mentioning further whether estimators are based on families of covariance functions or approximations), we refer to both, truncated-ML estimators based on families of covariance functions and approximations. The same is applied for the notion of truncated-modified log-likelihood functions based on either covariance functions or approximations. However, if $\{c_{\theta} \colon \theta \in \Theta\}$ satisfies the assumptions of Proposition~\ref{orderN}, a sequence of truncated-ML estimators 
$\big(\hat{\theta}_{n}(c)\big){}_{n \in \natnum}$ shall be simply called a sequence of ML estimators for $\theta_{0}$. Similarly, we will simply refer to a modified log-likelihood function when the given family $\{c_{\theta} \colon \theta \in \Theta\}$ is under the assumptions of Proposition~\ref{orderN}.
        
\begin{remark}
\normalfont The introduction of truncated-modified log-likelihood functions is not standard. Modified refers to the fact that the log-likelihood for the Gaussian density function of a random vector $(Z_{s_{1}}, \dotsc, Z_{s_{n}})$ is scaled by $-2/n$. This is common practice in the literature about ML estimators for covariance parameters under an increasing-domain asymptotic framework (see for instance \cite{bachoc2014asymptotic}, \cite{bachoc2018asymptotic} and also \cite{bachoc2020asymptotic}). The matrices $\covmat(\omega)$ and $\acovmat(\omega)$ are not necessarily positive-definite. In particular, $\acovmat(\omega)$ can be negative-definite. If the matrices $\covmat(\omega)$ and $\acovmat(\omega)$ are not positive-definite, we truncate the log-likelihood by a pseudo-determinant and -inverse to obtain the functions $\theta \mapsto l_{n}(\theta)(\omega)$ and $\theta \mapsto \tilde{l}_{n}(\theta)(\omega)$. Hence, the use of the  expression ``truncated''.       
\end{remark}

\begin{remark}
\normalfont As it was mentioned in Remark $2.2$ of \cite{bachoc2018asymptotic}, for given $\omega \in \Omega$, we allow the functions $\theta \mapsto l_{n}(\theta)(\omega)$ and $\theta \mapsto \tilde{l}_{n}(\theta)(\omega)$ to have more than one minimizer. In which case the asymptotic results given in Section~\ref{Consasym} hold true for any given sequence of truncated-ML estimators. With regard to the existence of a minimizer we refer to Remark $2.1$ in \cite{bachoc2014asymptotic}.       
\end{remark}

\subsection{Consistency and asymptotic normality of truncated-ML estimators} \label{Consasym} 

The main results of this section are that under suitable conditions on the families of covariance functions and approximations, truncated-ML estimators for covariance parameters are not only consistent (Theorem~\ref{thm1} and Corollary~\ref{corollarythm1}) but also asymptotically normal (Theorem~\ref{thm2} and Corollary~\ref{corollarythm2}). In particular, we will make use of the conditions presented in Assumptions~\ref{cov_a2} and~\ref{uniformly_approx_covfun}. However, in the context of random fields that are observed at randomly perturbed regular grid locations as defined in (\ref{randomgrid}), we will further make use of the following two technical conditions that were also imposed in \cite{bachoc2014asymptotic}. Associated to the common range $\mathcal{Q}$, of the process $X$, we define the set $D_{\tau} \coloneqq \bigcup_{z\in \mathbb{Z}^{d}\setminus\{0\}}(z + \tau U_{\mathcal{Q}})$, where $U_{\mathcal{Q}} \coloneqq \{u_{1} - u_{2} \colon u_{1} \in \mathcal{Q}, u_{2} \in \mathcal{Q}\}$ denotes the set of differences between two points in $\mathcal{Q}$.  

\begin{assumption}[Asymptotic identifiability around $\theta_{0}$]\label{cov_a3} 
For $\tau = 0$, there does not exists $\theta \neq \theta_{0}$ such that $c_{\theta}(z) = c_{\theta_{0}}(z)$ for all $z \in \mathbb{Z}^{d}$. If $\tau \neq 0$, there does not exists $\theta \neq \theta_{0}$ such that $s \mapsto c_{\theta}(s) - c_{\theta_{0}}(s)$ is zero a.e.\ with respect to the Lebesgue measure on $D_{\tau}$ and $c_{\theta}(0) = c_{\theta_{0}}(0)$.  
\end{assumption}

\begin{assumption}[Local identifiability around $\theta_{0}$]\label{cov_a4} 
For $\tau = 0$, there does not exists $\mathbb{R}^{p}\setminus\{0\} \ni \alpha = (\alpha_{1}, \dotsc, \alpha_{p})$ such that  
$\sum_{k=1}^{p}\alpha_{k}\frac{\partial c_{\theta_{0}}}{\partial\theta_{k}}(z) = 0$ for all $z \in \mathbb{Z}^{d}$. For $\tau \neq 0$, there does not exists $\mathbb{R}^{p}\setminus\{0\} \ni \alpha = (\alpha_{1}, \dotsc, \alpha_{p})$ such that $s \mapsto \sum_{k=1}^{p}\alpha_{k}\frac{\partial c_{\theta_{0}}}{\partial\theta_{k}}(s)$ is zero a.e.\ with respect to the Lebesgue measure on $D_{\tau}$ and $\sum_{k=1}^{p}\alpha_{k}\frac{\partial c_{\theta_{0}}}{\partial\theta_{k}}(0) = 0$.   
\end{assumption}

\begin{theorem} \label{thm1}
Let $\big(\hat{\theta}_{n}(\tilde{c})\big){}_{n \in \natnum}$ be a sequence of truncated-ML estimators for $\theta_{0}$ based on $\big\{(\tilde{c}_{m,\theta}) \colon \theta \in \Theta\big\}$. Assume that $\{c_{\theta} \colon \theta \in \Theta\}$ satisfies Assumption~\ref{cov_a2} (regarding~\ref{3.2.2}, $q=1$ and the continuity of first order partial derivatives is sufficient) and Assumption~\ref{cov_a3}. Suppose further that $\big\{(\tilde{c}_{m,\theta}) \colon \theta \in \Theta\big\}$ satisfies Assumption~\ref{uniformly_approx_covfun} (regarding \ref{3.4.3} and \ref{3.4.4}, $q=1$ and the continuity of first order partial derivatives is sufficient).
Then, we have that
\begin{gather*}
    \hat{\theta}_{n}(\tilde{c}) \toinp \theta_{0}.
\end{gather*}
\end{theorem}

The following corollary is immediate.

\begin{corollary} \label{corollarythm1}
Suppose that $\{c_{\theta} \colon \theta \in \Theta\}$ satisfies Assumption~\ref{cov_a2} (regarding~\ref{3.2.2}, $q=1$ and the continuity of first order partial derivatives is sufficient) and Assumption~\ref{cov_a3}. Then, we can conclude that a sequence of ML estimators $\big(\hat{\theta}_{n}(c)\big){}_{n \in \natnum}$ for $\theta_{0}$ is consistent.
\end{corollary}

Before we present the results about asymptotic normality, it is helpful to consider some additional notation. Let $K \in \natnum$, such that for any $\omega \in \Omega$, the sequences of functions 
\begin{equation} \label{shiftloglik1}
    \big(l_{n,K}\!(\theta)(\omega)\big)_{n \in \natnum} \coloneqq \big(l_{n+K-1}(\theta)(\omega)\big)_{n \in \natnum}
\end{equation}
and
\begin{equation} \label{shiftloglik2}
   \big(\tilde{l}_{n,K}\!(\theta)(\omega)\big)_{n \in \natnum} \coloneqq \big(\tilde{l}_{n+K-1}(\theta)(\omega)\big)_{n \in \natnum}
\end{equation}
are differentiable with respect to $\theta$. Note that if $\{c_{\theta} \colon \theta \in \Theta\}$ satisfies Assumption~\ref{cov_a2} and the collection $\big\{(\tilde{c}_{m,\theta}) \colon \theta \in \Theta\big\}$ satisfies Assumption~\ref{uniformly_approx_covfun}, then we know about the existence of such a $K$ under application of Proposition~\ref{orderN}. For the given $K \in \natnum$, on $(\Omega, \mathcal{F}, \mathbb{P})$, we introduce the sequence of random functions
\begin{gather*}
    \big\{\big(\underbrace{\omega \mapsto G_{n,K}(\omega, \theta)}_{\eqqcolon G_{n,K}(\theta)}\big)_{n\in \natnum} \colon \theta \in \Theta\big\},
\end{gather*}
where for $n \in \natnum$ and $\theta \in \Theta$, the random vector $G_{n,K}(\theta)$ has components $G_{j,n,K}(\theta)$, $j = 1, \dotsc, p$, with
\begin{gather*}
    G_{j,n,K}(\theta) = \frac{\partial l_{n,K}}{\partial \theta_{j}}(\theta) - \mathbb{E}\bigg[\frac{\partial l_{n,K}}{\partial \theta_{j}}(\theta) \;\bigg|\; S_{(n)}\bigg],
\end{gather*}
and thus 
\begin{equation} \label{cFun2}
    G_{n,K}(\theta) = \nabla l_{n,K}(\theta) - \mathbb{E}\big[\nabla l_{n,K}(\theta)\;\big|\; S_{(n)}\big]. 
\end{equation}
Similarly, on $(\Omega, \mathcal{F}, \mathbb{P})$, we introduce the sequence of random functions
\begin{gather*}
    \big\{\big(\underbrace{\omega \mapsto \widetilde{G}_{n,K}(\omega, \theta)}_{\eqqcolon \widetilde{G}_{n,K}(\theta)}\big)_{n\in \natnum} \colon \theta \in \Theta\big\},
\end{gather*}
where for any $n \in \natnum$ and $\theta \in \Theta$, the components of $\widetilde{G}_{n,K}(\theta)$ are given by 
\begin{gather*}
    \widetilde{G}_{j,n,K}(\theta) = \frac{\partial \tilde{l}_{n,K}}{\partial \theta_{j}}(\theta) - \mathbb{E}\bigg[\frac{\partial \tilde{l}_{n,K}}{\partial \theta_{j}}(\theta)\;\bigg|\; S_{(n)}\bigg], \quad j = 1, \dotsc, p,
\end{gather*}
and thus
\begin{equation} \label{cFun}
    \widetilde{G}_{n,K}(\theta) = \nabla \tilde{l}_{n,K}(\theta) - \mathbb{E}\big[\nabla \tilde{l}_{n,K}(\theta)\;\big|\; S_{(n)}\big].
\end{equation}
If the collection $\{c_{\theta} \colon \theta \in \Theta\}$ satisfies Assumption~\ref{cov_a2}, we simply write, for any $n \in \natnum$,
\begin{gather*}
    J_{G_{n}}(\theta_{0}) \coloneqq J_{G_{n,1}}(\theta_{0}),
\end{gather*}
for the random Jacobi-matrix of $\theta \mapsto G_{n,1}(\theta)$ evaluated at $\theta_{0}$. 

\begin{theorem} \label{thm2}
Let $\big(\hat{\theta}_{n}(\tilde{c})\big){}_{n \in \natnum}$ be an sequence of truncated-ML estimators for $\theta_{0}$ based on $\big\{(\tilde{c}_{m,\theta}) \colon \theta \in \Theta\big\}$. Suppose that $\{c_{\theta} \colon \theta \in \Theta\}$ satisfies Assumptions~\ref{cov_a2}, \ref{cov_a3} and \ref{cov_a4}. Suppose further that $\big\{(\tilde{c}_{m,\theta}) \colon \theta \in \Theta\big\}$ satisfies Assumption~\ref{uniformly_approx_covfun}.
Then, we have that
\begin{equation}\label{asymNorm}
    n^{1/2}\big(\hat{\theta}_{n}(\tilde{c}) - \theta_{0}\big) \toind \mathcal{N}(0,\Lambda^{-1}),
\end{equation}
where $\mathrm{S}_{p \times p} \ni \Lambda \succ 0$ is deterministic and such that
\begin{gather*}
    \frac{1}{2}J_{G_{n}}(\theta_{0}) \toinp \Lambda \xleftarrow[n \to \infty]{\mathbb{P}}\frac{1}{2}J_{\widetilde{G}_{n,N}}(\theta_{0}),
\end{gather*}
with $N \in \natnum$ as in Proposition~\ref{orderN}. 
\end{theorem}

\begin{corollary} \label{corollarythm2}
Suppose that $\{c_{\theta} \colon \theta \in \Theta\}$ satisfies Assumptions~\ref{cov_a2}, \ref{cov_a3} and \ref{cov_a4}. Then, we can conclude that a sequence of ML estimators $\big(\hat{\theta}_{n}(c)\big){}_{n \in \natnum}$ for $\theta_{0}$ is such that 
\begin{gather*}
    n^{1/2}(\hat{\theta}_{n}(c) - \theta_{0}) \toind \mathcal{N}(0,\Lambda^{-1}),
\end{gather*}
with $\Lambda$ as in Theorem~\ref{thm2}. 
\end{corollary}

\begin{remark} \label{normalremark}
\normalfont Under Assumptions~\ref{cov_a2}, for any $K \in \natnum$, $\E{\nabla l_{n,K}(\theta)\;|\; S_{(n)}} = 0$ with $\mathbb{P}$ probability one. However, even if $\{c_{\theta} \colon \theta \in \Theta\}$ is under Assumptions~\ref{cov_a2} and $\big\{(\tilde{c}_{m,\theta}) \colon \theta \in \Theta\big\}$ is under Assumptions~\ref{uniformly_approx_covfun}, it is not in general true that $\as$ $\mathbb{E}\big[\nabla \tilde{l}_{n,N}(\theta)\;|\; S_{(n)}\big] = 0$, where $N$ is as in Proposition~\ref{orderN}. Notice further that under Assumptions~\ref{cov_a2}, for $\omega \in \Omega$, $-(n/2)J_{G_{n}}(\theta_{0})(\omega)$ represents the second derivative of the log-likelihood $\theta \mapsto -(n/2)l_{n}(\theta)(\omega)$ based on $\{c_{\theta} \colon \theta \in \Theta\}$.    
\end{remark}

\section{Example of application: Generalized Wendland functions} \label{sec_application} 

In this section we work in the same setting as in Section~\ref{sec:context1}, but we additionally assume that $Z$ is isotropic. Explicitly, for the given family of covariance functions $\{c_{\theta} \colon \theta \in \Theta\}$, we assume that there exists a parametric family $\{\varphi_{\theta} \colon \theta \in \Theta\}$ such that for any $\theta \in \Theta$, $s \in \mathbb{R}^{d}$, $c_{\theta}(s) = \varphi_{\theta}(\nnorm{s})$. The family $\{\varphi_{\theta} \colon \theta \in \Theta\}$ is called the radial version of $\{c_{\theta} \colon \theta \in \Theta\}$. We can recycle the notation of Section~\ref{sec:regularity} and easily translate Assumptions~\ref{cov_a2} and \ref{uniformly_approx_covfun} by considering families of approximations $\big\{(\tilde{\varphi}_{m, \theta}) \colon \theta \in \Theta\big\}$ for $\{\varphi_{\theta} \colon \theta \in \Theta\}$ on $\mathbb{R}_{+}$. This allows us to readily recover the results of Sections~\ref{sec:uniformequivof mat} and \ref{MLE_results} for isotropic random fields. For the details we refer to Assumptions~\ref{cov_a2_iso} and~\ref{uniformly_approx_covfun_iso}, as well as Theorems~\ref{consistencyradial} and~\ref{asymnormalradial} in Appendix~\ref{sec_application1}.

In terms of an explicit family of radial covariance functions, we reconsider the generalized Wendland covariance function which we have already introduced in 
\eqref{StartWendland} of Section~\ref{sec:motivation}. Let $\Theta \ni \theta \coloneqq (\sigma^{2},\beta)$, where $\Theta \coloneqq [\sigma^{2}_{\min},\sigma^{2}_{\max}]\times[\beta_{\min}, \beta_{\max}]$, with $0 < \sigma^{2}_{\min} < \sigma^{2}_{\max} < \infty$ and $1-2\tau < \beta_{\min} < \beta_{\max} < \infty$. We assume that the covariance function of the random field $Z$ is given by $\phi_{\theta_{0}}(\nnorm{s})$, $s \in \mathbb{R}^{d}$, $\theta_{0} \in \Theta$, where $\phi_{\theta_{0}}$ belongs to the family $\{\phi_{\theta} \colon \theta \in \Theta\}$ which is defined by
\begin{equation} \label{wendland}
    \phi_{\theta}(t) \coloneqq \sigma^{2}\phi_{\nu, \kappa}\bigg(\frac{t}{\beta}\bigg), \quad t \in [0, \infty),
\end{equation}
where 
\begin{equation*}
    \phi_{\nu,\kappa}(r) \coloneqq \begin{cases} \frac{1}{\operatorname{B}(2\kappa,\nu + 1)}\int_{r}^{1}u(u^{2} - r^{2})^{\kappa-1}(1-u)^{\nu}du, & r \in [0,1), \\
                  0, & r \in [1, \infty),\end{cases}
\end{equation*}
compare to \eqref{StartWendland} of 
Section~\ref{sec:motivation}. We treat $\kappa$ and $\nu$ as given but such that $\kappa > 0$ and $\nu \geq (d+1)/2 + \kappa$. Notice that the latter restriction on $\kappa$ and $\nu$ makes sure that for any $\theta \in \Theta$, $\phi_{\theta}$ belongs to the class $\Phi_{d}$, the class of real valued and continuous functions, defined on $\mathbb{R}_{+}$, which are strictly positive at the origin and such that for any finite collection of points in $\mathbb{R}^{d}$, evaluation at the Euclidean norm of pairwise differences between points of the collection results in a non-negative definite matrix (see for example \cite{Gneiting2002}). Actually, in the latter reference it is argued that for $\kappa > 0$, $\phi_{\nu,\kappa} \in \Phi_{d}$ if and only if $\nu \geq (d+1)/2 + \kappa$. For the respective family defined on $\mathbb{R}^{d}$, we use the notation $w_{\theta}(s) \coloneqq \phi_{\theta}(\lVert s \rVert)$.

\begin{remark} \label{remarkrange}
\normalfont The restriction $\beta_{\min} > 1-2\tau$ is imposed to proof that the family $\{w_{\theta} \colon \theta \in \Theta\}$ satisfies Assumptions~\ref{cov_a3} and~\ref{cov_a4} (see the proof of Propositions~\ref{check2}). This is not surprising, as $1-2\tau$ defines the minimal spacing between pairs of distinct observation points of the randomly perturbed regular grid, defined in (\ref{randomgrid}) of Section~\ref{sec:context1}. Further, as we have noted that $\phi_{\nu,\kappa} \in \Phi_{d}$ if and only if $\nu \geq (d+1)/2 + \kappa$, the two smoothness parameters $\nu$ and $\kappa$ can not be estimated without further constraints. 
\end{remark}

\begin{proposition} \label{check1}
Let $\kappa > 4$. Then, the family $\{\phi_{\theta} \colon \theta \in \Theta\}$ satisfies Assumption~\ref{cov_a2_iso}, where for any $\theta \in \Theta$ and for any $q = 1,2,3$, $i_{1}, \dotsc, i_{q} \in \{1, \dotsc,p\}$, the functions $t \mapsto \phi_{\theta}(t)$ and $t \mapsto \frac{\partial^{q} \phi_{\theta}}{\partial \theta_{i_1}\cdots\partial \theta_{i_q}}(t)$ are continuous on $\mathbb{R}_{+}$. 
\end{proposition}

\begin{proposition} \label{check2}
Let $\kappa > 2$. Then, the family $\{w_{\theta} \colon \theta \in \Theta\}$ satisfies Assumptions~\ref{cov_a3} and \ref{cov_a4}.
\end{proposition}

Using Propositions~\ref{check1} and \ref{check2}, under application of Theorems~\ref{consistencyradial} and \ref{asymnormalradial} (recall also Corollaries~\ref{corollarythm1} and \ref{corollarythm2}), we obtain the following result:

\begin{proposition} \label{wendlandConsistent}
Let $\kappa > 4$. A sequence $\big(\hat{\theta}_{n}(\phi)\big){}_{n \in \natnum}$ of ML estimators for $\theta_{0}$ based on $\{\phi_{\theta} \colon \theta \in \Theta\}$ is consistent. Further there exists a non-random symmetric $p \times p$ matrix $\Lambda \succ 0$ such that 
\begin{gather*}
    n^{1/2}\big(\hat{\theta}_{n}(\phi) - \theta_{0}\big) \toind \mathcal{N}(0,\Lambda^{-1}).
\end{gather*}
\end{proposition}

\begin{remark} \label{remark:kappa}
It is worth to note that the restriction $\kappa > 4$ is only needed for the asymptotic distribution of ML estimators, respectively truncated-ML estimators. In particular, in Proposition~\ref{check1}, if one only demands conditions involving first order partial derivatives of $\phi_{\theta}$, with respect to $\theta$, $\kappa > 2$ is sufficient. With regard to consistency of the estimator $\big(\hat{\theta}_{n}(\phi)\big){}_{n \in \natnum}$ in Proposition~\ref{wendlandConsistent}, $\kappa > 2$ is sufficient as well. The same applies for the truncated-ML estimators considered in Examples~\ref{ex:1},~\ref{ex:2},~\ref{ex:3} and~\ref{ex:4}. Keeping in mind the differentiability conditions imposed in Assumption~\ref{cov_a2_iso}, the given restrictions on $\kappa$ are not surprising (compare also to \cite{bevilacqua2019estimation}, within the infill-domain asymptotic framework). 
\end{remark} 

We discuss four examples of generalized Wendland approximations.

\begin{example}[Truncation of $\phi_{\theta}$] \label{ex:1}
\normalfont Let $\{\phi_{\theta} \colon \theta \in \Theta\}$ be as in Proposition~\ref{check1}. Let $\big\{(\mathfrak{T}_{m, \theta})\colon \theta \in \Theta\big\}$ be defined as follows: 
For $\theta \in \Theta$ and $m \in \natnum$, we set,
\begin{gather*}
    \orange{\mathfrak{T}_{m, \theta}(t)} \coloneqq \phi_{\theta}\mathbbm{1}_{[0,C_{m}]}(t), \quad t \in \mathbb{R}_{+}, \; \text{$C_{m} \xrightarrow[]{} \infty$ as $m \xrightarrow[]{} \infty$.}
\end{gather*}

\begin{proposition} \label{propositionExample0}
A sequence $\big(\hat{\theta}_{n}(\mathfrak{T})\big){}_{n \in \natnum}$ of truncated-ML estimators for $\theta_{0}$ based on $\big\{(\mathfrak{T}_{m, \theta})\colon \theta \in \Theta\big\}$ is consistent and we have that
\begin{equation*}
    n^{1/2}\big(\hat{\theta}_{n}(\mathfrak{T}) - \theta_{0}\big) \toind \mathcal{N}(0,\Lambda^{-1}),
\end{equation*}
where $\Lambda$ is defined as in Proposition~\ref{wendlandConsistent}.
\end{proposition}

\end{example}

In the following we let $M < \infty$ denote a real constant, which is independent of $\beta \in [\beta_{\min}, \beta_{\max}]$ such that $\beta_{\max} \leq M$.

\begin{example}[Trimmed Bernstein polynomials] \label{ex:2}
\normalfont  Let $\{\phi_{\theta} \colon \theta \in \Theta\}$ be as in Proposition~\ref{check1}. We consider a family $\big\{(\mathfrak{P}_{m, \theta})\colon \theta \in \Theta\big\}$ defined as follows: For $\theta \in \Theta$ and $m \in \natnum$, we set for $t \in \mathbb{R}_{+}$,
\begin{gather*}
    \orange{\mathfrak{P}_{m, \theta}(t)} \coloneqq \begin{cases} B_{m, \theta}(t;b_{m}), & t \leq M, \\
                  0, & t > M,  \end{cases}
\end{gather*}
with 
\begin{gather*}
    B_{m, \theta}(t;b_{m}) = \sum_{k=0}^{m}\phi_{\theta}\bigg(b_{m}\frac{k}{m}\bigg)\binom{m}{k}\bigg(\frac{t}{b_{m}}\bigg)^{k}\bigg(1-\frac{t}{b_{m}}\bigg)^{m-k},
\end{gather*}
the Bernstein polynomial of the function $\phi_{\theta}$ on $[0,b_{m})$, where $b_{m} \tonormalm \infty$ and we assume that $b_{m} = o(m)$. Thus, for any $0 \leq k \leq m$,
\begin{gather*}
    b_{m}\frac{k+1}{m} - b_{m}\frac{k}{m} \tonormalm 0,
\end{gather*}
the distance between adjacent points converge to zero as $m$ approaches infinity. See also \cite{Chlodovsky1937} for an introduction of Bernstein polynomials on unbounded intervals.

\begin{proposition} \label{propositionExample1}
The family $\big\{(\mathfrak{P}_{m, \theta})\colon \theta \in \Theta\big\}$ satisfies Assumption~\ref{uniformly_approx_covfun_iso}.
\end{proposition}

Using Propositions~\ref{check1}, \ref{check2} and \ref{propositionExample1}, under application of Theorems~\ref{consistencyradial} and \ref{asymnormalradial}, we have proven the following result:

\begin{proposition} \label{proposition2Example1}
A sequence $\big(\hat{\theta}_{n}(\mathfrak{P})\big){}_{n \in \natnum}$ of truncated-ML estimators for $\theta_{0}$ based on $\big\{(\mathfrak{P}_{m, \theta})\colon \theta \in \Theta\big\}$ is consistent and we have that
\begin{equation*}
    n^{1/2}\big(\hat{\theta}_{n}(\mathfrak{P}) - \theta_{0}\big) \toind \mathcal{N}(0,\Lambda^{-1}),
\end{equation*}
where $\Lambda$ is defined as in Proposition~\ref{wendlandConsistent}. 
\end{proposition}

\end{example}

\begin{example}[Linear interpolation] \label{ex:3}
\normalfont Let $\{\phi_{\theta} \colon \theta \in \Theta\}$ be as in Proposition~\ref{check1}. For a given $m \in \natnum$, we consider a partition of the interval $[0,M]$, $0 = t_{0} \leq t_{1} \leq \dotsc \leq t_{N_{m}} = M$, where $N_{m} \tonormalm \infty$ and for $0 \leq k \leq N_{m}$, $t_{k+1}^{m} - t_{k}^{m} \tonormalm 0$.
Then, we define the family $\big\{(\mathfrak{L}_{m, \theta})\colon \theta \in \Theta\big\}$ as follows: 
For $\theta \in \Theta$ and $m \in \natnum$, we set for $t \in \mathbb{R}_{+}$,
\begin{gather*}
    \orange{\mathfrak{L}_{m, \theta}(t)} \coloneqq \begin{cases} I_{m, \theta}(t;N_{m}), & t \leq M, \\
                  0, & t > M,  \end{cases}
\end{gather*}
where 
\begin{gather*}
    I_{m, \theta}(t;N_{m}) = \begin{cases} \phi_{\theta}(t_{k}^{m}) + \frac{\phi_{\theta}(t_{k+1}^{m})- \phi_{\theta}(t_{k}^{m})}{t_{k+1}^{m} - t_{k}^{m}}(t-t_{k}^{m}), & t \in [t_{k}^{m}, t_{k+1}^{m}], \\
                  0, & t \not\in [t_{k}^{m}, t_{k+1}^{m}].  \end{cases}
\end{gather*}
Thus, for a given $m \in \natnum$, $\mathfrak{L}_{m, \theta}$ represents a linear interpolation of the function $\phi_{\theta}$ on the interval $[0,M]$.

\begin{proposition} \label{propositionExample2}
The family $\big\{(\mathfrak{L}_{m, \theta})\colon \theta \in \Theta\big\}$ satisfies Assumption~\ref{uniformly_approx_covfun_iso}.
\end{proposition}

Using Propositions~\ref{check1}, \ref{check2} and \ref{propositionExample2}, under application of Theorems~\ref{consistencyradial} and \ref{asymnormalradial}, we have further proven the following result:

\begin{proposition} \label{proposition2Example2}
A sequence $\big(\hat{\theta}_{n}(\mathfrak{L})\big){}_{n \in \natnum}$ of truncated-ML estimators for $\theta_{0}$ based on $\big\{(\mathfrak{L}_{m, \theta})\colon \theta \in \Theta\big\}$ is consistent and we have that
\begin{equation*}
    n^{1/2}\big(\hat{\theta}_{n}(\mathfrak{L}) - \theta_{0}\big) \toind \mathcal{N}(0,\Lambda^{-1}),
\end{equation*}
where $\Lambda$ is defined as in Proposition~\ref{wendlandConsistent}. 
\end{proposition}

\end{example}

\begin{example}[Vanishing nugget effect] \label{ex:4}
Let $\{\phi_{\theta} \colon \theta \in \Theta\}$ be as in Proposition~\ref{check1} and consider a family $\big\{\big(\tilde{\phi}_{m, \theta}\big) \colon \theta \in \Theta\big\}$ that satisfies Assumption~\ref{uniformly_approx_covfun_iso}. Then, define for any $\theta \in \Theta$ and $m \in \natnum$, the function
\begin{equation} \label{originShift}
   \orange{\mathfrak{S}_{m, \theta}(t)} \coloneqq \begin{cases} \tilde{\phi}_{m, \theta}(t) + \delta(m), & t = 0, \\
                  \tilde{\phi}_{m, \theta}(t), & t \neq 0,  \end{cases}
\end{equation}
where $(\delta(m)){}_{m \in \natnum}$ is independent of $\theta \in \Theta$ and $t \in \mathbb{R}_{+}$ and such that $\delta(m) \xrightarrow{} 0$, as $m \xrightarrow{} \infty$. Note that since the family $\{\phi_{\theta} \colon \theta \in \Theta\}$ satisfies Assumption~\ref{uniformly_approx_covfun_iso}, we could also choose $\tilde{\phi}_{m, \theta} \equiv \phi_{\theta}$ in (\ref{originShift}).   

\begin{proposition} \label{propositionExample4}
A sequence $\big(\hat{\theta}_{n}\big(\mathfrak{S}\big)\big){}_{n \in \natnum}$ of truncated-ML estimators for $\theta_{0}$ based on $\big\{\big(\mathfrak{S}_{m, \theta}\big)\colon \theta \in \Theta\big\}$ is consistent and we have that
\begin{equation*}
    n^{1/2}\big(\hat{\theta}_{n}(\mathfrak{S}) - \theta_{0}\big) \toind \mathcal{N}(0,\Lambda^{-1}),
\end{equation*}
where $\Lambda$ is defined as in Proposition~\ref{wendlandConsistent}.
\end{proposition}

\end{example}

\begin{remark}\label{remark:askin}
As it was already mentioned in the introduction, computing \eqref{wendland} is costly. However, if $\kappa$ is a positive integer, closed form solutions of \eqref{wendland} exist. More specifically, if $\kappa = k \in \natnum$, then 
\begin{gather*}
    \phi_{\nu, k}(r) = A_{\nu + k}(r)P_{k}(r),
\end{gather*}
where $P_{k}$ is a polynomial of order $k$ and $A_{\nu + k}$ the Askey function (\cite{askey1973radial}) of order $\nu + k$,
\begin{gather*}
    A_{\nu + k}(r) = \begin{cases} (1-r)^{\nu + k}, & 0 \leq r < 1, \\
    0, & r \geq 1.  \end{cases}
\end{gather*}
In addition, if $\kappa \in (\natnum - 1/2)$, a positive half-integer, it is shown in \cite{b50d71882155d7fde924ff08f9ecdb1b} that further closed form solutions of \eqref{wendland}, involving polynomial, logarithmic and square root terms, exist. Thus, in the specific example of generalized Wendland covariance functions, covariance approximations will facilitate computing \eqref{wendland} when $\kappa \notin \natnum \cup (\natnum - 1/2)$.
\end{remark}

\section{Covariance taper approximations: Beyond compactly supported covariance functions} \label{sec:app2}

Asymptotic properties of (regular) tapered-ML estimators were addressed in both the infill- and increasing-domain asymptotic framework (see \cite{MR2504203}, \cite{MR2549562}, \cite{MR3173747} and \cite{Asymptotic_tapering}). The direct functional approximation approach studied here can be combined with covariance tapering. Given observations of $S$, it is known that under weak assumptions on the presumed covariance function, ML estimators based on tapered covariance functions (tapered-ML estimators) preserve consistency (see \cite{Asymptotic_tapering}, in particular Corollary~2 in the increasing-domain framework). However, this is the case for covariance tapers that have a compact support which is not fixed, but rather grows to the entire $\mathbb{R}^{d}$ as the number of observations from $S$ increases. Within an increasing-domain asymptotic framework, given a fixed compact support of the covariance taper, one can in general not expect tapered-ML estimators to be consistent. Still, under suitable conditions, tapered-ML estimators asymptotically minimize the Kullback-Leibler divergence (see for instance Theorem 3.3 in \cite{bachoc2018asymptotic}). Given the theory developed here, we can readily recover the same result for truncated-tapered ML estimators, ML estimators based on tapered covariance function, where the covariance taper is replaced with a functional approximation of it. To be more formal, let us remain in the setting of Section~\ref{sec:context}, but assume that $Z$ has true and unknown covariance function $k_{\theta_{0}}$, $\theta_{0} \in \Theta$, which belongs to a family $\{k_{\theta} \colon \theta \in \Theta\}$ which satisfies:
\begin{itemize}
    \item For any $s \in \mathbb{R}^{d}$, $\theta \mapsto k_{\theta}(s)$ is continuously differentiable
    \item There exist constants $A < \infty$ and $\alpha > 0$ such that for all $i = 1, \dotsc, p$, for all $s \in \mathbb{R}^{d}$ and for all $\theta \in \Theta$, $\lvert k_{\theta}(s) \rvert \leq A/(1+\lVert s \rVert^{d+\alpha})$ and $\big \lvert  \frac{\partial k_{\theta}}{\partial \theta_{i}}(s) \big \rvert \leq A/(1+\lVert s \rVert^{d+\alpha})$
    \item $\{k_{\theta} \colon \theta \in \Theta\}$ satisfies \ref{3.2.3} of Assumption~\ref{cov_a2}. 
\end{itemize}
The given assumptions are very weak and satisfied for instance for the Mat\'{e}rn family (see also Condition~2.1 in \cite{bachoc2014asymptotic} or Remark~\ref{remark31}). Then, we consider a fixed covariance taper $s \mapsto t_{\theta^{\prime}_{0}}(s)$, $\theta^{\prime}_{0} \in \Theta^{\prime}$, $\Theta^{\prime} \subset \mathbb{R}^{l}$, compact and convex. We assume that $t_{\theta^{\prime}_{0}}$ belongs to a family of tapers $\{t_{\theta^{\prime}} \colon \theta^{\prime} \in \Theta^{\prime}\}$ that satisfies Assumption~\ref{cov_a2} (regarding~\ref{3.2.2}, $q=1$ and the continuity of first order partial derivatives is sufficient). As we have seen in Section~\ref{sec_application} (Proposition~\ref{check1}), we may choose, with $\theta^{\prime}_{0} = (\beta_{0}, 1)$, $\kappa > 2$, $\nu \geq (d+1)/2 + \kappa$, a generalized Wendland taper (see also Remark~\ref{remark:kappa}). In the given context it is more convenient to write $t_{\beta_{0}} \coloneqq t_{\theta^{\prime}_{0}}$, where $\beta_{0}$ is the taper range, that is $t_{\beta_{0}}(s) = 0$ for $\lVert s \rVert \geq \beta_{0}$. Based on a finite collection $S_{(n)}$ of $S$, on $(\Omega, \mathcal{F}, \mathbb{P})$, we then define the tapered $n \times n$ covariance matrix ${R_{n,\theta}}_{i,j} \coloneqq k_{\theta}(S_{i} - S_{j})t_{\beta_{0}}(S_{i} - S_{j})$, $1 \leq i,j \leq n$. Additionally, we consider a covariance matrix approximation 
\begin{gather*}
    \widetilde{R}_{{n,\theta}_{i,j}} = k_{\theta}(S_{i} - S_{j})\tilde{t}_{r(n), \theta^{\prime}_{0}}(S_{i} - S_{j}), \quad  1 \leq i,j \leq n, \; \text{$r(n) \xrightarrow[]{} \infty$ as $n \xrightarrow[]{} \infty$,} 
\end{gather*}
of $R_{n,\theta}$, where $(\tilde{t}_{m,\theta^{\prime}_{0}})$ is a sequence of functions that belongs to a family of taper approximations 
$\{(\tilde{t}_{m,\theta^{\prime}}) \colon \theta^{\prime} \in \Theta^{\prime}\}$, for which Assumption~\ref{sec_cov_approx} applies. Again, we write $\tilde{t}_{m,\theta^{\prime}_{0}} \coloneqq \tilde{t}_{m,\beta_{0}}$, $m \in \natnum$, to highlight the fixed range parameter. We note that the results of Lemma~\ref{eigenbehave} and Proposition~\ref{orderN} remain true with $\covmat$ and $\acovmat$ replaced with $R_{n,\theta}$ and $\widetilde{R}_{n,\theta}$, respectively. We know (see Remark~\ref{remark:conditionaldisti}) that the conditional distribution of $Z_{(n)}$ given $S_{(n)}$ is given by the random variable $\omega \mapsto \mathcal{N}(0,K_{n,\theta_{0}}(\omega))$. On the other hand, we can assume a misspecified distribution $\omega \mapsto \mathcal{N}(0,R_{n,\theta}(\omega))$, where the true covariance matrix is replaced with the tapered covariance matrix $R_{n,\theta}(\omega)$, $\theta \in \Theta$. Then, we define the scaled (see \cite{bachoc2018asymptotic}) conditional Kullback-Leibler divergence of $\mathcal{N}(0,R_{n,\theta})$ from $\mathcal{N}(0,K_{n,\theta_{0}})$,  
\begin{gather*}
     d_{n,\theta} \coloneqq \frac{1}{n}\operatorname{log}\big(\operatorname{det}(R_{n,\theta}K_{n,\theta_{0}}^{-1})\big) + \frac{1}{n}\operatorname{tr}(K_{n,\theta_{0}}R_{n,\theta}^{-1})-1. 
\end{gather*}
The distribution $\mathcal{N}(0,R_{n,\theta})$ shall be called a regular taper miss-specified distribution. If we choose $n \geq N$ ($N$ as in Proposition~\ref{orderN}), we can even further misspecify the distribution of $Z_{(n)}$ given $S_{(n)}$ by replacing $R_{n,\theta}$ with $\widetilde{R}_{n,\theta}$ in $\mathcal{N}(0,R_{n,\theta})$. This gives rise to the scaled conditional Kullback-Leibler divergence of $\mathcal{N}(0,\widetilde{R}_{n,\theta})$ from $\mathcal{N}(0,K_{n,\theta_{0}})$, 
\begin{gather*}
     \tilde{d}_{n,\theta} \coloneqq \frac{1}{n}\operatorname{log}\big(\operatorname{det}(\widetilde{R}_{n,\theta}K_{n,\theta_{0}}^{-1})\big) + \frac{1}{n}\operatorname{tr}(K_{n,\theta_{0}}\widetilde{R}_{n,\theta}^{-1})-1.
\end{gather*}
We use the notation $\big(\hat{\theta}_{n}(kt)\big){}_{n \in \natnum}$ and $\big(\hat{\theta}_{n}(k\tilde{t})\big){}_{n \in \natnum}$ for ML and truncated-ML estimators for $\theta_{0}$ with respect to $\{k_{\theta}t_{\beta_{0}} \colon \theta \in \Theta\}$ and $\{(k_{\theta}\tilde{t}_{m,\beta_{0}}) \colon \theta \in \Theta\}$, respectively. In accordance with the literature about tapered-ML estimators, the estimators $\big(\hat{\theta}_{n}(kt)\big){}_{n \in \natnum}$ and 
$\big(\hat{\theta}_{n}(k\tilde{t})\big){}_{n \in \natnum}$ are then further referred to as tapered-ML estimators and truncated-tapered ML estimators, respectively. We can now state the following theorem:

\begin{theorem} \label{thm:taper}
We have that $\as$ 
\begin{equation} \label{eq:taper1}
 \sup_{\theta \in \Theta} \big\lvert d_{n,\theta} - \tilde{d}_{n,\theta} \big \rvert \tonormal 0, 
\end{equation}
and as $n \xrightarrow[]{} \infty$,
\begin{equation} \label{eq:taper2}
d_{n,\hat{\theta}_{n}(k\tilde{t})} = \inf_{\theta \in \Theta}d_{n,\theta} + \delta_{n},
\end{equation}
where $\delta_{n} \toinp 0$.
\end{theorem}

Therefore, in the given scenario, truncated-tapered ML estimators asymptotically minimize the conditional Kullback-Leibler divergence of taper misspecified distributions from the true distribution (compare also to Theorem 3.3 in \cite{bachoc2018asymptotic}). Thus, in terms of Kullback-Leibler divergence, truncated-tapered ML estimators and tapered-ML estimators perform asymptotically equally well.

\section{Discussion and outlook} \label{sec:conc}
With the introduction of truncated-likelihood functions, we allow for more far-reaching forms of covariance approximations, such as linear interpolations or polynomial approximations. Our approximation approach relates directly to the presumed covariance function. Thus, combinations with existing approximation methods such as low-rank or covariance tapering approaches are well possible. We studied the quality of truncated-ML estimators from an asymptotic point of view. For compactly supported covariance functions, the conditions imposed in Sections~\ref{sec:regularity} and \ref{MLE_results} permit us to obtain truncated-ML estimators that are asymptotically well-behaving. That is, we obtain estimators that are consistent and asymptotically normal. Our proof strategies were strongly influenced by \cite{bachoc2014asymptotic}. We have provided a comprehensive analysis for the family of generalized Wendland covariance functions. That is, we give precise conditions on smoothness, variance and range parameters, under which ML estimators for variance and range parameters are consistent and asymptotically normal. To our knowledge, this does not exists in the literature so far (compare also to \cite{bevilacqua2019estimation}, within the infill-domain asymptotic context). Further, we gave four examples of generalized Wendland approximations, for which truncated-ML estimators preserve consistency and asymptotic normality.

We now discuss some open questions. Our results on consistency and asymptotic normality depend on the condition that correlations vanish beyond a certain distance. It would be of interest to recover the consistency and asymptotic normality results for truncated-ML estimators, where the assumption of a compact support is dropped. To this end, we recall that the imposed conditions on covariance functions and approximations resulted in the uniform asymptotic equivalence of covariance matrices and approximations. Using this, we  established the existence of a positive integer $N$, after which covariance matrix approximations remain positive-definite. Expanding to non-compactly supported covariance function, this result  remains unchanged, as long as covariance matrices and approximations are uniformly asymptotically equivalent (uniformly on the parameter and sample space). Thus, in this case, consistency and asymptotic normality can be recovered, even when presumed covariance functions are no longer compactly supported. However, as a mere condition, the asymptotic equivalence of covariance matrices and approximations is of little practical importance. Thus, the case of non-compactly supported covariance functions deserves further attention.

From a more applied point of view, our results provide a strong theoretical basis for further research. It remains to test and extend the given examples of covariance approximations. The four examples of generalized Wendland approximations and their effect on parameter estimations were discussed from a theoretical point of view. An important next step is to provide numerical implementations and practical comparisons.

In conclusion, for large datasets built upon correlated data, the present work provides an essential missing piece in the area of covariance approximations.


\begin{appendix} 
\section{Covariance approximations for isotropic random fields} \label{sec_application1} 

We consider families of approximations $\big\{(\tilde{\varphi}_{m, \theta}) \colon \theta \in \Theta\big\}$ for $\{\varphi_{\theta} \colon \theta \in \Theta\}$ on $\mathbb{R}_{+}$ and translate (recycling the notation of Section~\ref{sec:regularity}) Assumptions~\ref{cov_a2} and \ref{uniformly_approx_covfun} as follows:

\begin{assumption}[Regularity conditions on $\varphi_{\theta}$]\label{cov_a2_iso} 
~\begin{enumerate}[label=(\arabic*)]

\item \label{6.1.1} There exist real constants $C$, $L < \infty$, which are independent of $\theta \in \Theta$, such that $\varphi_{\theta} \in \mathcal{B}_{\text{C}}(\mathbb{R}_{+};S_{\theta})$, with $S_{\theta} \subset [0,C]$ and $\supnorm{\varphi_{\theta}} \leq L$.

\item \label{6.1.2} For any $t \in \mathbb{R}_{+}$, the first, second and third order partial derivatives of $\theta \mapsto \varphi_{\theta}(t)$ exist. In addition, for any $q = 1,2,3$, $i_{1}, \dotsc, i_{q} \in \{1, \dotsc,p\}$, $\frac{\partial^{q} \varphi_{\theta}}{\partial \theta_{i_1}\cdots\partial \theta_{i_q}} \in \mathcal{B}_{\text{C}}(\mathbb{R}_{+};S_{\theta}(i_{1}, \dotsc, i_{q}))$, where there exist constants $C^{\prime}\text{, }L^{\prime} < \infty$, which are independent of $\theta \in \Theta$, such that $S_{\theta}(i_{1}, \dotsc, i_{q}) \subset [0,C^{\prime}]$ and $\big \lVert \frac{\partial^{q} \varphi_{\theta}}{\partial \theta_{i_1}\cdots\partial \theta_{i_q}}\big \rVert_{\infty} \leq~L^{\prime}$. 

\item \label{6.1.3} Fourier inversion holds, that is for any $\theta \in \Theta$,
\begin{gather*}
    \varphi_{\theta}(\nnorm{s}) = \int_{\mathbb{R}^{d}}\!\hat{c}_{\theta}(f)\eeexp^{\ii \langle f, s \rangle}\mathrm{d}\!f,
\end{gather*}
where $\Theta \times \mathbb{R}^{d} \ni (\theta, f) \mapsto \hat{c}_{\theta}(f)$ is continuous and strictly positive. 
\end{enumerate}
\end{assumption}

\begin{assumption}[Regularity conditions on $\tilde{\varphi}_{m,\theta}$] \label{uniformly_approx_covfun_iso} 
~\begin{enumerate}[label=(\arabic*)]

\item \label{6.2.1} For any $\theta \in \Theta$, for any $m \in \natnum$, the function $\tilde{\varphi}_{m,\theta}\colon (\mathbb{R}_{+}, \mathfrak{B}(\mathbb{R}_{+})) \to (\mathbb{R}, \mathfrak{B}(\mathbb{R}))$ is measurable. 

\item \label{6.2.2} For any $m \in \natnum$, $\tilde{\varphi}_{m,\theta}$ satisfies \ref{6.1.1} of Assumption~\ref{cov_a2_iso}, where respective constants $\widetilde{C}$ and $\widetilde{L}$ can be further chosen independently of $m \in \natnum$. 

\item \label{6.2.3} $\sup_{\theta \in \Theta}\supnorm{\tilde{\varphi}_{m,\theta}-\varphi_{\theta}} \tonormalm 0$.

\item \label{6.2.4} For any $m \in \natnum$, $\tilde{\varphi}_{m,\theta}$ satisfies \ref{6.2.2} of Assumption~\ref{cov_a2_iso}, where respective constants $\widetilde{C}^{\prime}$ and $\widetilde{L}^{\prime}$ can be further chosen independently of $m \in \natnum$. 

\item \label{6.2.5} For any $q = 1,2,3$, $i_{1}, \dotsc, i_{q} \in \{1, \dotsc,p\}$, we have that
\begin{gather*}
    \sup_{\theta \in \Theta}\supnorm{\frac{\partial^{q} \tilde{\varphi}_{m,\theta}}{\partial \theta_{i_1}\cdots\partial \theta_{i_q}} - \frac{\partial^{q} \varphi_{\theta}}{\partial \theta_{i_1}\cdots\partial \theta_{i_q}}} \tonormalm 0.
\end{gather*}
\end{enumerate}
\end{assumption}

Note that the family $\{\varphi_{\theta} \colon \theta \in \Theta\}$ satisfies Assumption~\ref{cov_a2_iso} if and only if $\{c_{\theta} \colon \theta \in \Theta\}$ satisfies Assumption~\ref{cov_a2}. Further, for any $n \in \natnum$ and $\theta \in \Theta$, we have that 
\begin{gather*}
    \covmat = \big[\varphi_{\theta}(\nnorm{S_{i}- S_{j}})\big]_{1 \leq i,j \leq n},
\end{gather*}
on $(\Omega, \mathcal{F}, \mathbb{P})$. Thus, a sequence of truncated-ML estimators for $\theta_{0}$ based on $\{c_{\theta} \colon \theta \in \Theta\}$ is a sequence of truncated-ML estimators for $\theta_{0}$ based on $\{\varphi_{\theta} \colon \theta \in \Theta\}$. If we define a sequence of truncated-ML estimators $\big(\hat{\theta}_{n}(\tilde{\varphi})\big){}_{n \in \natnum}$ for $\theta_{0}$ based on a given $\big\{(\tilde{\varphi}_{m, \theta}) \colon \theta \in \Theta\big\}$ upon replacing $\acovmat$ in (\ref{pseudocriterion}) with the random $n \times n$ matrix $\big[\tilde{\varphi}_{r(n), \theta}(\nnorm{S_{i}- S_{j}})\big]_{1 \leq i,j \leq n}$, we can recover the results of Sections~\ref{sec:uniformequivof mat} and \ref{MLE_results}:

\begin{theorem} \label{consistencyradial}
Let $\big(\hat{\theta}_{n}(\tilde{\varphi})\big){}_{n \in \natnum}$ be a sequence of truncated-ML estimators for $\theta_{0}$ based on $\big\{(\tilde{\varphi}_{m, \theta}) \colon \theta \in \Theta\big\}$. Assume that $\{\varphi_{\theta} \colon \theta \in \Theta\}$ satisfies Assumption~\ref{cov_a2_iso} (regarding~\ref{6.1.2}, $q=1$ and the continuity of first order partial derivatives is sufficient) and $\{c_{\theta} \colon \theta \in \Theta\}$ satisfies Assumption~\ref{cov_a3}. Suppose further that $\big\{(\tilde{\varphi}_{m, \theta}) \colon \theta \in \Theta\big\}$ satisfies Assumption~\ref{uniformly_approx_covfun_iso} (regarding \ref{6.2.4} and \ref{6.2.5}, $q=1$ and the continuity of first order partial derivatives is sufficient).
Then,
\begin{gather*}
    \hat{\theta}_{n}(\tilde{\varphi}) \toinp \theta_{0}.
\end{gather*}
\end{theorem}

\begin{theorem} \label{asymnormalradial}
Let $\big(\hat{\theta}_{n}(\tilde{\varphi})\big){}_{n \in \natnum}$ be an sequence of truncated-ML estimators for $\theta_{0}$ based on $\big\{(\tilde{\varphi}_{m, \theta}) \colon \theta \in \Theta\big\}$. Suppose that $\{\varphi_{\theta} \colon \theta \in \Theta\}$ satisfies Assumption~\ref{cov_a2_iso} and $\{c_{\theta} \colon \theta \in \Theta\}$ satisfies Assumptions~\ref{cov_a3} and \ref{cov_a4}. Assume further that $\big\{(\tilde{\varphi}_{m, \theta}) \colon \theta \in \Theta\big\}$ satisfies Assumption~\ref{uniformly_approx_covfun_iso}.
Then, we have that
\begin{equation*}
    n^{1/2}\big(\hat{\theta}_{n}(\tilde{\varphi}) - \theta_{0}\big) \toind \mathcal{N}(0,\Lambda^{-1}),
\end{equation*}
with $\Lambda$ as in Theorem~\ref{thm2}.
\end{theorem}

\section{Supporting results} \label{appn}
Let $r \colon \natnum \to \natnum$ be such that $r(n) \xrightarrow{} \infty$ as $n \xrightarrow{} \infty$. For the families $\{c_{\theta} \colon \theta \in \Theta\}$ and $\big\{(\tilde{c}_{m,\theta}) \colon \theta \in \Theta\big\}$, we introduce, for any $n \in \natnum$ and $\theta \in \Theta$, for an arbitrary $s_{(n)} \in \mathcal{G}_{n}$, for any $q = 1,2,3$, $i_{1}, \dotsc, i_{q} \in \{1, \dotsc,p\}$, the non-random $n \times n$ matrices
\begin{gather*}
    \frac{\partial^{q} \rcovmat}{\partial \theta_{i_1}\cdots\partial \theta_{i_q}} \coloneqq \bigg[\frac{\partial^{q} c_{\theta}}{\partial \theta_{i_1}\cdots\partial \theta_{i_q}}(s_{i}- s_{j})\bigg]_{1 \leq i,j \leq n},
\end{gather*}
and
\begin{gather*}
    \frac{\partial^{q} \racovmat}{\partial \theta_{i_1}\cdots\partial \theta_{i_q}} \coloneqq \bigg[\frac{\partial^{q} \tilde{c}_{r(n), \theta}}{\partial \theta_{i_1}\cdots\partial \theta_{i_q}}(s_{i}- s_{j})\bigg]_{1 \leq i,j \leq n},
\end{gather*}
whenever the above partial derivatives with respect to $\theta$ exist. Further, for Borel measurable sequences of functions $\big\{(\tilde{c}_{m,\theta}) \colon \theta \in \Theta\big\}$, we introduce, on $(\Omega, \mathcal{F}, \mathbb{P})$, the $n \times n$ random matrices
\begin{gather*} 
    \omega \mapsto \frac{\partial^{q} \covmat}{\partial \theta_{i_1}\cdots\partial \theta_{i_q}}(\omega) \coloneqq \frac{\partial^{q} \rrcovmat}{\partial \theta_{i_1}\cdots\partial \theta_{i_q}},
\end{gather*}
and
\begin{gather*} 
    \omega \mapsto \frac{\partial^{q} \acovmat}{\partial \theta_{i_1}\cdots\partial \theta_{i_q}}(\omega) \coloneqq \frac{\partial^{q} \rracovmat}{\partial \theta_{i_1}\cdots\partial \theta_{i_q}},
\end{gather*}
whenever the above partial derivatives with respect to $\theta$ exist.

\begin{lemma} \label{lemma1}
Let $C$, $L < \infty$ be some real constants. Consider $g \colon \mathbb{R}^{d} \to \mathbb{R}_{+}$ such that $g \in \mathcal{B}_{\text{C}}(\mathbb{R}^{d};S)$, with $S \subset \mathrm{B}[0; C]$ and $\supnorm{g} \leq L$. Then, for any $i \in \natnum$, for any sequence $(s_{j})_{j \in \natnum} \in \mathcal{G}$,
\begin{equation} \label{lemma11}
    \sum_{j \in \natnum}g(s_{i} - s_{j}) \leq L R(d,C,\tau),
\end{equation}
where $R(d,C,\tau) \coloneqq (2^{2d}dC^{d-1})/\Delta_{\tau}^{d}$, with $\Delta_{\tau} = 1-2\tau$. Further, we also have that
\begin{equation} \label{lemma12}
    \sum_{\substack{j \in \natnum \\ \norm{v_{i}-v_{j}}_{\infty} \geq C + 1}}g(s_{i} - s_{j}) =0. 
\end{equation}
\end{lemma}

\begin{remark} \label{remarkLemmacomp}
\normalfont We would like to point out that Lemma~\ref{lemma1} resembles Lemmas D.1 and D.3 of \cite{bachoc2014asymptotic}, where $f \colon \mathbb{R}^{d} \to \mathbb{R}_{+}$, which is such that $f(s) \leq 1/(1 + \norm{s}_{\infty}^{d+1})$, is replaced with a compactly supported function $g$, defined as in Lemma~\ref{lemma1}.
\end{remark}

\begin{lemma} \label{lemma3}
Let $\widetilde{C}$, $\widetilde{L} < \infty$ be some real constants. Consider a sequence of functions $(g_{m}){}_{m \in \natnum}$, with values in $\mathbb{R}_{+}$, where for any $m \in \natnum$, $g_{m} \in \mathcal{B}_{\text{C}}(\mathbb{R}^{d};\widetilde{S}_{m})$, with $\widetilde{S}_{m} \subset \mathrm{B}\big[0; \widetilde{C}\big]$ and $\supnorm{g_{m}} \leq \widetilde{L}$. Then, for any $i \in \natnum$, for any sequence $(s_{j})_{j \in \natnum} \in \mathcal{G}$,
\begin{equation} \label{lemma31}
    \sup_{m \in \natnum}\sum_{j \in \natnum}g_{m}(s_{i} - s_{j}) \leq \widetilde{L} R(d,\widetilde{C},\tau),
\end{equation}
where $R(d,\widetilde{C},\tau) \coloneqq (2^{2d}d\widetilde{C}^{d-1})/\Delta_{\tau}^{d}$, with $\Delta_{\tau} = 1-2\tau$. Further we also have that
\begin{equation} \label{lemma32}
    \sup_{m \in \natnum}\sum_{\substack{j \in \natnum \\ \norm{v_{i}-v_{j}}_{\infty} \geq \widetilde{C} + 1}}g_{m}(s_{i} - s_{j}) =0. 
\end{equation}
\end{lemma}

\begin{lemma} \label{eigenbehave1x1}
Assume that $\{c_{\theta} \colon \theta \in \Theta\}$ satisfies \ref{3.2.1} and \ref{3.2.3} of Assumption~\ref{cov_a2}. Consider $\big\{(\tilde{c}_{m, \theta}) \colon \theta \in \Theta\big\}$ that satisfies \ref{3.4.0}, \ref{3.4.1} and \ref{3.4.2} of Assumption~\ref{uniformly_approx_covfun}.
Then, we have that 
\begin{equation} \label{r1.22}
        \sup_{n\in \natnum}\sup_{s_{(n)} \in \mathcal{G}_{n}}\sup_{\theta \in \Theta}\big \lVert\rcovmat\big \rVert_{2} < \infty \text{,  }\sup_{n\in \natnum}\sup_{s_{(n)} \in \mathcal{G}_{n}}\sup_{\theta \in \Theta}\big \lVert\racovmat \big \rVert_{2} < \infty,
\end{equation}
and in particular
\begin{equation} \label{r1.2}
        \sup_{s_{(n)} \in \mathcal{G}_{n}}\sup_{\theta \in \Theta}\big \lVert\rcovmat - \racovmat\big \rVert_{2} \tonormal 0.
\end{equation}
Further, we have that
\begin{equation} \label{r1.3}
       \inf_{n\in \natnum}\inf_{s_{(n)} \in \mathcal{G}_{n}}\inf_{\theta \in \Theta}\lambda_{n}\big(\rcovmat\big) > 0, 
\end{equation}
and there exists $N \in \natnum$ such that  
\begin{equation} \label{r1.1}
       \inf_{n\geq N}\inf_{s_{(n)} \in \mathcal{G}_{n}}\inf_{\theta \in \Theta}\lambda_{n}\big(\racovmat\big) > 0.
\end{equation}
\end{lemma}

\begin{corollary} \label{corollary1}
Let $\{c_{\theta} \colon \theta \in \Theta\}$, $\big\{(\tilde{c}_{m, \theta}) \colon \theta \in \Theta\big\}$ and $N$ be as in Lemma~\ref{eigenbehave1x1}. Then, we have that
\begin{gather*}
    \sup_{n \in \natnum}\sup_{s_{(n)} \in \mathcal{G}_{n}}\sup_{\theta \in \Theta} \big\lVert\rcovmat^{-1}\big \rVert_{2} < \infty\text{, }\sup_{n\geq N}\sup_{s_{(n)} \in \mathcal{G}_{n}}\sup_{\theta \in \Theta}\big\lVert\racovmat^{-1}\big \rVert_{2} < \infty.
\end{gather*}
In addition we can conclude that
\begin{gather*}
    \sup_{s_{(n)} \in \mathcal{G}_{n}}\sup_{\theta \in \Theta}\big\lVert\rcovmat^{+} - \racovmat^{+}\big \rVert_{2} \tonormal 0.
\end{gather*} 
In particular we have that $\as$ 
\begin{gather*}
\sup_{n \in \natnum}\sup_{\theta \in \Theta}\big\lVert\covmat^{-1}\big \rVert_{2},\; \sup_{n\geq N}\sup_{\theta \in \Theta}\big\lVert\acovmat^{-1}\big \rVert_{2} < \infty\text{ and }\sup_{\theta \in \Theta}\big\lVert\covmat^{+} - \acovmat^{+}\big \rVert_{2} \tonormal 0.
\end{gather*}
\end{corollary}

\begin{lemma} \label{lemmaEigenbehavePartials}
Suppose that $\{c_{\theta} \colon \theta \in \Theta\}$ satisfies \ref{3.2.2} of Assumption~\ref{cov_a2}. Consider $\big\{(\tilde{c}_{m, \theta}) \colon \theta \in \Theta\big\}$ that satisfies \ref{3.4.0}, \ref{3.4.3} and \ref{3.4.4} of Assumption~\ref{uniformly_approx_covfun}. Then, for any $q = 1,2,3$, $i_{1}, \dotsc, i_{q} \in \{1, \dotsc,p\}$, we have that (\ref{r1.22}) and (\ref{r1.2}) of Lemma~\ref{eigenbehave1x1} are satisfied with $\rcovmat$ and $\racovmat$ replaced with the respective partial derivatives $\frac{\partial^{q} \rcovmat}{\partial \theta_{i_1}\cdots\partial \theta_{i_q}}$ and $\frac{\partial^{q} \racovmat}{\partial \theta_{i_1}\cdots\partial \theta_{i_q}}$. In particular, for any $q = 1,2,3$, $i_{1}, \dotsc, i_{q} \in \{1, \dotsc,p\}$, we have that $\as$ 
\begin{gather*}
    \sup_{n \in \natnum}\sup_{\theta \in \Theta}\bigg\lVert\frac{\partial^{q} \covmat}{\partial \theta_{i_1}\cdots\partial \theta_{i_q}}\bigg \rVert_{2} < \infty \text{, } \sup_{n \in \natnum}\sup_{\theta \in \Theta}\bigg\lVert\frac{\partial^{q} \acovmat}{\partial \theta_{i_1}\cdots\partial \theta_{i_q}}\bigg \rVert_{2} < \infty,
\end{gather*}
and in addition it is true that for any $q = 1,2,3$, $i_{1}, \dotsc, i_{q} \in \{1, \dotsc,p\}$, $\as$
\begin{gather*}
    \sup_{\theta \in \Theta}\bigg\lVert\frac{\partial^{q} \covmat}{\partial \theta_{i_1}\cdots\partial \theta_{i_q}}- \frac{\partial^{q} \acovmat}{\partial \theta_{i_1}\cdots\partial \theta_{i_q}}\bigg \rVert_{2} \tonormal 0.
\end{gather*}
\end{lemma}

\begin{lemma} \label{lemmalogdet}
Let $I \in \natnum$ be fixed. On $(\Omega, \mathcal{F}, \mathbb{P})$, for $k = 1, \dotsc, I$, we consider a sequence of $n \times n$ random symmetric matrices $\big(\widetilde{A}_{k,n,\theta}\big){}_{n \in \natnum}$, $\theta \in \Theta$, such that $\as$, for any $k = 1, \dotsc, I$, $\sup_{n \in \natnum}\sup_{\theta \in \Theta}\big\lVert\widetilde{A}_{k,n,\theta}\big\rVert < \infty$. Further we assume that there exists $N \in \natnum$ such that $\as$, for $k = 1, \dotsc, I$, $\inf_{n\geq N}\inf_{\theta \in \Theta}\lambda_{n}\big(\widetilde{A}_{k,n,\theta}\big) > 0$. Let $\big(A_{k,n,\theta}\big){}_{n \in \natnum}$, $\theta \in \Theta$, $k = 1, \dotsc, I$, be another sequence of $n \times n$ random symmetric matrices, defined on the same probability space, which is such that $\as$, for $k = 1, \dotsc, I$, 
\begin{gather*}
    \sup_{n \in \natnum}\sup_{\theta \in \Theta}\specnorm{A_{k,n,\theta}} < \infty\text{ and }\inf_{n\geq N}\inf_{\theta \in \Theta}\lambda_{n}(A_{k,n,\theta}) > 0.
\end{gather*}
Finally we also assume that $\as$, for any $k = 1, \dotsc, I$,  
\begin{gather*}
     \sup_{\theta \in \Theta}\big\lVert\widetilde{A}_{k,n,\theta}-A_{k,n,\theta}\big \rVert_{2} \tonormal 0.
\end{gather*}
Then, we have that $\as$
\begin{equation*} 
    \sup_{\theta \in \Theta}\norm{\frac{1}{n}\operatorname{log}
    \bigg(\operatorname{det}_{+}\bigg(\prod_{k=1}^{I}A_{k,n,\theta}\bigg)\bigg) - \frac{1}{n}
    \operatorname{log}\bigg(\operatorname{det}_{+}\bigg(\prod_{k=1}^{I}\widetilde{A}_{k,n,\theta}\bigg)\bigg)} \tonormal 0. 
\end{equation*}
\end{lemma}

\begin{lemma} \label{quadraticForm}
On $(\Omega, \mathcal{F}, \mathbb{P})$, consider two sequences of $n \times n$ random matrices $\big(A_{n,\theta}\big){}_{n \in \natnum}$ and $\big(\widetilde{A}_{n,\theta}\big){}_{n \in \natnum}$, $\theta \in \Theta$, such that $\as$
\begin{gather*}
    \sup_{\theta \in \Theta}\big \lVert A_{n,\theta} - \widetilde{A}_{n,\theta}\big \rVert_{2} \tonormal 0.
\end{gather*}
Then, we have that
\begin{equation} \label{ConvQuadratic}
    \sup_{\theta \in \Theta}\frac{1}{n}\big \lvert \langle Z_{(n)}, A_{n,\theta}Z_{(n)}\rangle - \langle Z_{(n)}, \widetilde{A}_{n,\theta}Z_{(n)}\rangle \big \rvert \toinp 0.
\end{equation}
\end{lemma}

\begin{lemma} \label{lemma2}
Suppose that $\{c_{\theta} \colon \theta \in \Theta\}$ satisfies Assumption~\ref{cov_a2} and \ref{cov_a4} (regularity conditions for partial derivatives up to order $q =2$ are sufficient). Suppose further that $\big\{(\tilde{c}_{m,\theta}) \colon \theta \in \Theta\big\}$ satisfies Assumption~\ref{uniformly_approx_covfun} (regularity conditions for partial derivatives up to order $q =2$ are sufficient). Let $N$ be as in Proposition~\ref{orderN} and define $\big\{(G_{n,N}(\theta)){}_{n \in \natnum}\colon \theta \in \Theta\big\}$ and $\big\{\big(\widetilde{G}_{n,N}(\theta)\big){}_{n \in \natnum}\colon \theta \in \Theta\big\}$ as in (\ref{cFun2}) and (\ref{cFun}), respectively. We then have that
 \begin{equation} \label{l2.1}
    \big \lVert J_{\widetilde{G}_{n,N}}(\theta_{0}) - J_{G_{n,N}}(\theta_{0})\big \rVert_{2} \toinp 0.
\end{equation}
Further, we conclude that the random $p\times p$ matrix $J_{\widetilde{G}_{n,N}}(\theta_{0})$ converges in probability $\mathbb{P}$ to a non-random matrix $2\Lambda$, where $\mathrm{S}_{p \times p} \ni \Lambda \succ 0.$ \end{lemma}


\section{Proofs} \label{appnB} 

\subsection{Proof of results in Appendix~\ref{appn}} \label{sec_support}

\begin{proof}[Proof of Lemma~\ref{lemma1}]
Let $\left(s_{j}\right)_{j \in \natnum} \in \mathcal{G}$. For $j \in \natnum$ such that $\norm{v_{i}-v_{j}}_{\infty} \geq C + 1$ we have that $\norm{s_{i}-s_{j}}_{\infty} \geq C$ and thus $\nnorm{s_{i}-s_{j}} \geq C$ as well (since $\norm{w}_{\infty}\leq \nnorm{w}$ for any $w \in \mathbb{R}^{d}$). Therefore, $(\ref{lemma12})$ follows since we have assumed that $g$ has compact support $S \subset \mathrm{B}\left[0; C\right]$. The proof of (\ref{lemma11}) depends on the fact that there exists a minimal spacing $\Delta_{\tau}> 0$ between any two distinct observation points (see \eqref{IncreasingDomain}). This allows us to show that for some arbitrary $i \in \natnum$, if $N_{s_{i}, C}$ denotes the cardinality of the set $\{j \in \natnum 
\colon \nnorm{s_{j} -s_{i}} \leq C\} \subset \{j \in \natnum 
\colon \norm{s_{j} -s_{i}}_{\infty} \leq C\}$, we have that $N_{s_{i}, C} \leq R\left(d,C,\tau\right)$. For a complete argument one could for example consider the proof of Lemma~4 in \cite{Asymptotic_tapering}. Using this we can estimate,
\begin{align*}
    \sum_{j\in \natnum}g\left(s_{i}-s_{j}\right)
    &= \sum_{j \in \natnum}g\left(s_{i}-s_{j}\right)\mathbbm{1}_{[0,C]}\left(\nnorm{s_{i}-s_{j}}\right) \\
    &\leq L\sum_{j \in \natnum}\mathbbm{1}_{[0,C+1]}\left(\nnorm{s_{i}-s_{j}}\right) \\
    &\leq L R\left(d,C,\tau\right),
\end{align*}
and thus also (\ref{lemma11}) is proven.    
\end{proof}

\begin{proof}[Proof of Lemma~\ref{lemma3}]
The proof is similar to the proof of Lemma~\ref{lemma1} and hence we consider the lemma as proven.
\end{proof}

\begin{proof}[Proof of Lemma~\ref{eigenbehave1x1}]
Let $C$, $L$ and $\widetilde{C}$, $\widetilde{L}$ be defined as in \ref{3.2.1} of Assumption~\ref{cov_a2} and \ref{3.4.1} of Assumption~\ref{uniformly_approx_covfun}, respectively. We use Lemma~\ref{lemma1} to show that there exists a real constant $M>0$, which does not depend on $n \in \natnum$, $s_{(n)} \in \mathcal{G}_{n}$ and $\theta \in \Theta$ such that for any $n \in \natnum$, $s_{(n)} \in \mathcal{G}_{n}$, 
\begin{equation} \label{UniformBoundMaxEigenvalue}
    \max\left\{\big \lVert\rcovmat\big \rVert_{2}, \big \lVert \racovmat \big \rVert_{2}\right\} \leq M.
\end{equation}
To see this, let $C_{*} \coloneqq \max\{C, \widetilde{C}\}$ and $L_{*} \coloneqq \max\{L, \widetilde{L}\}$. Using \ref{3.2.1} of Assumption~\ref{cov_a2}, we have that for any $\theta \in \Theta$, $c_{\theta} \in \mathcal{B}_{\text{C}}(\mathbb{R}^{d};S_{\theta})$, where now $S_{\theta} \subset \mathrm{B}\left[0;C_{*}\right]$ and $\supnorm{c_{\theta}} \leq L_{*}$, with $C_{*}$ and $L_{*}$ finite constants that are independent of $n \in \natnum$ and $\theta \in \Theta$. Thus we can write, for any $n \in \natnum$, $s_{(n)} \in \mathcal{G}_{n}$ and $\theta \in \Theta$, by Gershgorin circle theorem, 
\begin{align*}
    \specnorm{\rcovmat} &\leq \max_{i = 1, \dotsc,n}\sum_{j=1}^{n}\norm{c_{\theta}\left(s_{i}-s_{j}\right)} \\
    &\leq \sup_{i \in \natnum}\sum_{j\in\natnum}\norm{c_{\theta}\left(s_{i}-s_{j}\right)} \leq L_{*} R\left(d,C_{*},\tau\right) \eqqcolon M,
\end{align*}
under application of Lemma~\ref{lemma1}, with $R\left(d,C_{*},\tau\right) = (2^{2d}d C_{*}^{d-1})/\Delta_{\tau}^{d}$, where $\Delta_{\tau} = 1-2\tau$. Note that $M$ is independent of $n \in \natnum$, $s_{(n)} \in \mathcal{G}_{n}$ and $\theta \in \Theta$. Similarly, by \ref{3.4.1} of Assumption~\ref{uniformly_approx_covfun} we then use Lemma~\ref{lemma3}, together with Gershgorin circle theorem, to show that for any $n \in \natnum$, $s_{(n)} \in \mathcal{G}_{n}$ and $\theta \in \Theta$, $\lVert\racovmat\rVert_{2}\leq M$ as well. This shows (\ref{UniformBoundMaxEigenvalue}). 
Thus, we have established that
\begin{gather*}
    \sup_{n \in \natnum}\sup_{s_{(n)} \in \mathcal{G}_{n}}\sup_{\theta \in \Theta}\specnorm{\racovmat} \leq M,
\end{gather*}
and \begin{gather*}
    \sup_{n \in \natnum}\sup_{s_{(n)} \in \mathcal{G}_{n}}\sup_{\theta \in \Theta}\specnorm{\rcovmat} \leq M,
\end{gather*} 
and therefore (\ref{r1.22}) of Lemma~\ref{eigenbehave1x1} is verified. It is shown in \cite{bachoc2014asymptotic} (Proposition~D.4) that because of the increasing-domain setting, where there exists a minimal distance between any two observation points (see \eqref{IncreasingDomain}), and since \ref{3.2.3} of Assumption~\ref{cov_a2} is satisfied,
\begin{gather*}
    \inf_{n\in \natnum}\inf_{x_{(n)} \in \mathcal{Q}^{n}}\inf_{\theta \in \Theta}\lambda_{n}\left(\rcovmat\right) > 0.
\end{gather*}
This shows (\ref{r1.3}) of Lemma~\ref{eigenbehave1x1}. Using this result, we can fix some $\delta > 0$ (small enough, independent of $n\in \natnum$, $s_{(n)} \in \mathcal{G}_{n}$ and $\theta \in \Theta$), such that for any $s_{(n)} \in \mathcal{G}_{n}$, 
\begin{gather*}
0 < \varepsilon \coloneqq \frac{\delta}{R\left(d,C_{*},\tau\right)} < \min_{\nnorm{a} = 1}\langle a, \rcovmat a\rangle. 
\end{gather*}
For the above $\delta > 0$, we can then find $N \in \natnum$ such that, 
\begin{equation} \label{conf_frob1}
    \sup_{n \geq N}\sup_{s_{(n)} \in \mathcal{G}_{n}}\sup_{\theta \in \Theta}\specnorm{\rcovmat - \racovmat} \leq \delta.
\end{equation}
This is valid since for the given $\varepsilon > 0$, by the uniform convergence of $(\tilde{c}_{r(n), \theta})$ to $c_{\theta}$ (see \ref{3.4.2} of Assumption~\ref{uniformly_approx_covfun}), we find $N \in \natnum$ such that for any $n \geq N$, for any $s_{(n)} \in \mathcal{G}_{n}$ and $1 \leq i,j \leq n$, 
\begin{gather*}
    \sup_{\theta \in \Theta}\norm{\left[\rcovmat\right]_{i,j} - \left[\racovmat\right]_{i,j}} < \varepsilon.
\end{gather*}
Then, if we define
\begin{gather*}
    \mathbb{R}^{d} \ni s \mapsto g_{r(n)}(s) \coloneqq \sup_{\theta \in \Theta}\norm{\left(c_{\theta} - \tilde{c}_{r(n), \theta}\right)(s)}, \quad n \geq N,
\end{gather*}
since we have assumed that the families $\{c_{\theta} \colon \theta \in \Theta\}$ and $\big\{(\tilde{c}_{m,\theta}) \colon \theta \in \Theta\big\}$ have compact supports, which belong to $\mathrm{B}\left[0;C_{*}\right]$, we have that $g_{r(n)}(s) = 0$ for $\lVert s \rVert \geq C_{*}$. Thus, by Gershgorin circle theorem, under application of Lemma~\ref{lemma3}, for $n \geq N$ and $s_{(n)} \in \mathcal{G}_{n}$,  
\begin{align*}
    \specnorm{\racovmat - \rcovmat} & \leq \max_{i = 1, \dotsc,n}\sum_{j=1}^{n}g_{r(n)}\left(s_{i}-s_{j}\right)  \\
    &\leq \sup_{i \in \natnum}\sum_{j \in \natnum}g_{r(n)}\left(s_{i}-s_{j}\right) \leq  \varepsilon R\left(d,C_{*},\tau\right).
\end{align*}
Since $\varepsilon R\left(d,C_{*},\tau\right)$ is independent of $n\in \natnum$, $s_{(n)} \in \mathcal{G}_{n}$ and $\theta \in \Theta$, we can conclude that (\ref{conf_frob1}) must be satisfied. Using (\ref{conf_frob1}), we have, for $n \geq N$, $s_{(n)} \in \mathcal{G}_{n}$ and $\theta \in \Theta$, and for vectors $a$ such that $\nnorm{a} = 1$, that  
\begin{align*}
\norm{\langle a, \rcovmat a\rangle - \langle a, \racovmat a\rangle} &= \norm{\langle a, \big(\rcovmat-\racovmat\big)a\rangle}\\
              &\leq \specnorm{\rcovmat - \racovmat} \leq \delta, 
\end{align*}
under application of the Cauchy--Schwarz inequality. In conclusion we have for vectors $a$ such that $\nnorm{a} = 1$, for $n \geq N$, $s_{(n)} \in \mathcal{G}_{n}$ and $\theta \in \Theta$, 
\begin{align*}
    \min_{\nnorm{a} = 1}\langle a, \rcovmat a\rangle - \delta &\leq \min_{\nnorm{a} = 1}\langle a, \racovmat a\rangle.
\end{align*}
But we know that $\inf_{n \geq N}\inf_{s_{(n)} \in \mathcal{G}_{n}}\inf_{\theta \in \Theta}\min_{\nnorm{a} = 1}\langle a, \rcovmat a\rangle > 0$ and $\delta > 0$ was chosen small enough (but otherwise arbitrary). Thus, we have also proven (\ref{r1.1}) of Lemma~\ref{eigenbehave1x1}. Notice that (\ref{r1.2}) is proven with (\ref{conf_frob1}), hence the proof of Lemma~\ref{eigenbehave1x1} is complete. 
\end{proof}

\begin{proof}[Proof of Corollary~\ref{corollary1}]
This follows from Lemma~\ref{eigenbehave1x1}. 
\end{proof}

\begin{proof}[Proof of Lemma~\ref{lemmaEigenbehavePartials}]
We omit a formal argument and argue that one can proof Lemma~\ref{lemmaEigenbehavePartials} using the same way of reasoning as in the proof of Lemma~\ref{eigenbehave1x1}.
\end{proof}

\begin{proof}[Proof of Lemma~\ref{lemmalogdet}]
For $n \geq N$ ($N$ as in the statement) and $\theta \in \Theta$, we can write $\as$
\begin{gather*}
    \operatorname{det}_{+}\bigg(\prod_{k=1}^{I}\widetilde{A}_{k,n,\theta}\bigg) = \operatorname{det}\big(\underbrace{\widetilde{A}_{I,n,\theta}^{1/2} \cdots \widetilde{A}_{2,n,\theta}^{1/2}\widetilde{A}_{1,n,\theta}\widetilde{A}_{2,n,\theta}^{1/2}\cdots\widetilde{A}_{I,n,\theta}^{1/2}}_{\eqqcolon \widetilde{B}_{n,\theta}}\big), 
\end{gather*}
and 
\begin{gather*}
    \operatorname{det}_{+}\bigg(\prod_{k=1}^{I}A_{k,n,\theta}\bigg) = \operatorname{det}\big(\underbrace{A_{I,n,\theta}^{1/2} \cdots A_{2,n,\theta}^{1/2}A_{1,n,\theta}A_{2,n,\theta}^{1/2}\cdots A_{I,n,\theta}^{1/2}}_{\eqqcolon B_{n,\theta}}\big). 
\end{gather*}
Note that $\widetilde{B}_{n,\theta}$ and $B_{n,\theta}$ are random symmetric matrices. Further, for each of the random symmetric matrices
\begin{gather*}
  \widetilde{A}_{I,n,\theta}^{1/2}, \dotsc, \widetilde{A}_{2,n,\theta}^{1/2}, \; \widetilde{A}_{1,n,\theta}, \; A_{I,n,\theta}^{1/2}, \dotsc, A_{2,n,\theta}^{1/2},\; A_{1,n,\theta},  
\end{gather*}
the smallest eigenvalue is strictly greater that zero, $\as$, uniformly in $n \geq N$ and $\theta \in \Theta$ and hence we have that
\begin{equation} \label{cb0}
    \inf_{n\geq N}\inf_{\theta \in \Theta}\lambda_{n}\left(B_{n,\theta}\right) > 0 \text{ and }\inf_{n\geq N}\inf_{\theta \in \Theta}\lambda_{n}\left(\widetilde{B}_{n,\theta}\right) > 0, \; \; \as
\end{equation}
In addition, since $\as$ for $k = 1, \dotsc, I$, by assumption
\begin{gather*}
    \sup_{n \geq N}\sup_{\theta \in \Theta}\specnorm{\widetilde{A}_{k,n,\theta}}, \; \sup_{n \geq N}\sup_{\theta \in \Theta}\specnorm{\widetilde{A}_{k,n,\theta}} < \infty,
\end{gather*}
and
\begin{gather*} 
    \sup_{\theta \in \Theta}\specnorm{\widetilde{A}_{k,n,\theta} - A_{k,n,\theta}} \tonormal 0,
\end{gather*}
we also have that
\begin{equation} \label{cb1}
    \sup_{n \geq N}\sup_{\theta \in \Theta}\specnorm{\widetilde{B}_{n,\theta}}, \; \sup_{n \geq N}\sup_{\theta \in \Theta}\specnorm{B_{n,\theta}} < \infty,
\end{equation}
and
\begin{equation} \label{cb2}
    \sup_{\theta \in \Theta}\specnorm{B_{n,\theta} - \widetilde{B}_{n,\theta}} \tonormal 0.
\end{equation}
Using (\ref{cb0}), (\ref{cb1}) and (\ref{cb2}), we pick $\delta > 0$ arbitrary and define 
\begin{gather*}
    0 < \varepsilon \coloneqq \frac{\delta}{\sup_{n \geq N}\sup_{\theta \in \Theta}\specnorm{B_{n,\theta}}}
\end{gather*}
such that for some given integer $N^{*} \geq N$, $\as$,
\begin{gather*}
    \sup_{n \geq N^{*}}\sup_{\theta \in \Theta}\specnorm{\widetilde{B}_{n,\theta}^{-1} - B_{n,\theta}^{-1}}<\varepsilon.
\end{gather*}
Now write
\begin{align} \label{logdet2}
    \frac{1}{n}\operatorname{log}\left(\frac{\operatorname{det}(\prod_{k=1}^{I}A_{k,n,\theta})}{\operatorname{det}(\prod_{k=1}^{I}\widetilde{A}_{k,n,\theta})}\right) &= \frac{1}{n}\llog{\ddet{B_{n,\theta}\widetilde{B}_{n,\theta}^{-1}}} \nonumber \\
    &= \frac{1}{n}\tr\left(\llog{B_{n,\theta}\widetilde{B}_{n,\theta}^{-1}}\right) \nonumber \\
    &= \frac{1}{n}\sum_{i=1}^{n}\llog{\lambda_{i}\left(B_{n,\theta}\widetilde{B}_{n,\theta}^{-1}\right)}.
\end{align}
We can then estimate (\ref{logdet2}) from above and below as
\begin{gather*}
    \llog{\lambda_{n}\left(B_{n,\theta}\widetilde{B}_{n,\theta}^{-1}\right)} \leq \frac{1}{n}\sum_{i=1}^{n}\llog{\lambda_{i}\left(B_{n,\theta}\widetilde{B}_{n,\theta}^{-1}\right)} \leq \llog{\lambda_{1}\left(B_{n,\theta}\widetilde{B}_{n,\theta}^{-1}\right)}.
\end{gather*}
But for the given $\varepsilon > 0$, for $n \geq N^{*}$, we have that $\as$ \begin{align*}
    \lambda_{1}\left(B_{n,\theta}\widetilde{B}_{n,\theta}^{-1}\right) &\leq \lambda_{1}\left(B_{n,\theta}\right)\lambda_{1}\left(\widetilde{B}_{n,\theta}^{-1}\right) \\
    &= \big\lVert B_{n,\theta}\big\rVert_{2}\big\lVert\widetilde{B}_{n,\theta}^{-1}\big\rVert_{2} \\
    &\leq \big\lVert B_{n,\theta}\big\rVert_{2}\big\lVert \widetilde{B}_{n,\theta}^{-1} - B_{n,\theta}^{-1}\big\rVert_{2} + \big\lVert B_{n,\theta}\big\rVert_{2}\big\lVert B_{n,\theta}^{-1}\big\rVert_{2} \\
    &\leq 1 + \delta.
\end{align*}
On the other hand, by (\ref{cb2}), we also have that for $n \geq N^{*}$, $\as$
\begin{align*}
\lambda_{n}\left(B_{n,\theta}\widetilde{B}_{n,\theta}^{-1}\right) &\geq \lambda_{n}\left(B_{n,\theta}\right)\lambda_{n}\left(\widetilde{B}_{n,\theta}^{-1}\right) \\
&= \left(\min_{\{a\colon \nnorm{a} = 1\}}\left\langle a, B_{n,\theta}a\right\rangle\right)\left(\min_{\{a\colon \nnorm{a} = 1\}}\left\langle a, \widetilde{B}_{n,\theta}^{-1}a\right\rangle\right) \\
&\geq \left(\min_{\{a\colon \nnorm{a} = 1\}}\left\langle a, B_{n,\theta}a\right\rangle\right)\left(\min_{\{a\colon \nnorm{a} = 1\}}\left\langle a, B_{n,\theta}^{-1}a\right\rangle - \varepsilon\right) \\
&\geq 1-\delta.
\end{align*}
Since $\delta >0$ was arbitrary and independent of $\theta \in \Theta$, the lemma is proven.
\end{proof}

\begin{proof}[Proof of Lemma~\ref{quadraticForm}]
First, using the Cauchy--Schwarz inequality and the compatibility of the spectral norm with the Euclidean norm, we can estimate $\as$
\begin{gather*}
    \frac{1}{n}\norm{\langle Z_{(n)}, A_{n,\theta}Z_{(n)}\rangle - \langle Z_{(n)}, \widetilde{A}_{n,\theta}Z_{(n)}\rangle} \leq \sup_{\theta \in \Theta}\specnorm{A_{n,\theta}-\widetilde{A}_{n,\theta}}\frac{\nnorm{Z_{(n)}}^{2}}{n}
\end{gather*}
Let us fix some arbitrary $\varepsilon > 0$ such that for $n$ large enough we have that $\as$,
\begin{gather*}
    \sup_{\theta \in \Theta}\specnorm{A_{n,\theta}-\widetilde{A}_{n,\theta}} < \varepsilon.
\end{gather*}
Then, let $\delta > 0$ be arbitrary and notice that
\begin{gather*}
    \mathbb{P}\left(\varepsilon n^{-1}\nnorm{Z_{(n)}}^{2} > \delta \; \middle | \; S_{(n)} = s_{(n)}\right) = \mathbb{P}\left(\varepsilon n^{-1}\nnorm{\rtcovmat^{1/2}V_{n}}^{2} > \delta\right),
\end{gather*}
where $V_{n}$ is a Gauss vector, defined on $\left(\Omega, \mathcal{F}, \mathbb{P}\right)$, with zero-mean vector and identity covariance matrix. Then, we use Markov's inequality to estimate
\begin{align*}
    \mathbb{P}\left(\varepsilon n^{-1}\nnorm{\rtcovmat^{1/2}V_{n}}^{2} > \delta\right) &\leq \varepsilon n^{-1}\delta^{-1}\E{ \nnorm{\rtcovmat^{1/2}V_{n}}^{2}} \\
    &\leq \varepsilon \delta^{-1} \specnorm{\rtcovmat^{1/2}}^{2},
\end{align*}
where the latter term is bounded uniformly in $s_{(n)} \in \mathcal{G}_{n}$ and $n \in \natnum$ (see Lemma~\ref{eigenbehave1x1}). Thus we conclude that 
\begin{gather*}
    \sup_{s_{(n)} \in \mathcal{G}_{n}}\mathbb{P}\left( \sup_{\theta \in \theta}\specnorm{A_{n,\theta} - \widetilde{A}_{n,\theta}}\frac{\nnorm{Z_{(n)}}^{2}}{n} > \delta\; \middle | \; S_{(n)} = s_{(n)}\right) \tonormal 0,
\end{gather*}
which shows that
\begin{gather*}
    \sup_{\theta \in \theta}\specnorm{A_{n,\theta} - \widetilde{A}_{n,\theta}}\frac{\nnorm{Z_{(n)}}^{2}}{n} \toinp 0,
\end{gather*}
and thus the proof is complete.
\end{proof}

\begin{proof}[Proof of Lemma~\ref{lemma2}]
For $n \in \natnum$, let $h(n) = n + N-1$. Then, for $k = 1, \dotsc,p$, we have that $\as$
\begin{align*}
    \begin{split}
    \frac{\partial \tilde{l}_{n,N}}{\partial \theta_{k}}(\theta_{0}) &= \frac{1}{h(n)}\bigg(\tr\bigg(\shifttacovmat^{-1}\frac{\partial \shifttacovmat}{\partial \theta_{k}}\bigg)\\
    &\quad-\langle Z_{(h(n))} , \shifttacovmat^{-1}\frac{\partial \shifttacovmat}{\partial \theta_{k}}\shifttacovmat^{-1}Z_{(h(n))}\rangle\bigg), 
    \end{split}
\end{align*}   
and
\begin{equation*}
    \begin{split}
    \E{\frac{\partial \tilde{l}_{n,N}}{\partial \theta_{k}}(\theta_{0})\;\middle|\; S_{(h(n))}} &= \frac{1}{h(n)}\bigg(\tr\bigg(\shifttacovmat^{-1}\frac{\partial \shifttacovmat}{\partial \theta_{k}}\bigg) \\
    & \quad - \tr\bigg(\shifttacovmat^{-1}\frac{\partial \shifttacovmat}{\partial \theta_{k}}\shifttacovmat^{-1}\shifttcovmat\bigg)\bigg).
\end{split}
\end{equation*}
Similar expressions can then be calculated for $l_{n,N}$ based on $\shifttcovmat$. We can further calculate, for $n \in \natnum$, for $1 \leq k,l \leq p$, $\as$,
\begin{gather*}
   \frac{\partial^{2}\tilde{l}_{n,N}}{\partial \theta_{k}\partial \theta_{l}}(\theta_{0}) = \frac{1}{h(n)}\tr\left(\widetilde{A}_{1,h(n),\theta_{0}}^{kl}\right) + \frac{1}{h(n)}\langle Z_{(h(n))}, \widetilde{A}_{2,h(n),\theta_{0}}^{kl}Z_{(h(n))}\rangle,
\end{gather*}
where 
\begin{align} \label{eq:pl1}
    \widetilde{A}_{1,h(n),\theta_{0}}^{kl} &\coloneqq -\shifttacovmat^{-1}\frac{\partial \shifttacovmat}{\partial \theta_{k}}\shifttacovmat^{-1}\frac{\partial \shifttacovmat}{\partial \theta_{l}}+  \shifttacovmat^{-1}\frac{\partial^{2}\shifttacovmat}{\partial \theta_{k}\partial \theta_{l}},
\end{align}
and
\begin{align} \label{eq:pl2}
    \begin{split}
    \widetilde{A}_{2,h(n),\theta_{0}}^{kl} &\coloneqq 2\shifttacovmat^{-1}\frac{\partial \shifttacovmat}{\partial \theta_{k}}\shifttacovmat^{-1}\frac{\partial \shifttacovmat}{\partial \theta_{l}}\shifttacovmat^{-1}\\
    &\quad-\shifttacovmat^{-1}\frac{\partial^{2}\shifttacovmat}{\partial \theta_{k}\partial \theta_{l}}\shifttacovmat^{-1}.
    \end{split}
\end{align}
In addition, for $n \in \natnum$, we also have that $\as$,
\begin{align*}
   \frac{\partial \left(\E{\frac{\partial \tilde{l}_{n,N}}{\partial \theta_{l}}(\theta_{0})\; \middle| \; S_{(h(n))}}\right)}{\partial\theta_{k}}
    &= \frac{1}{h(n)}\bigg(\tr\left(\widetilde{A}_{1,h(n),\theta_{0}}^{kl}\right) + \tr\left(\widetilde{A}_{2,h(n),\theta_{0}}^{kl}\shifttcovmat\right)\bigg).
\end{align*}
Again, similar expressions can be obtained for $l_{n,N}$ based on $\shifttcovmat$, where for $n \in \natnum$, $1 \leq k,l \leq p$, the respective terms $A_{1,h(n),\theta_{0}}^{kl}$ and $A_{2,h(n),\theta_{0}}^{kl}$ are defined as in (\ref{eq:pl1}) and (\ref{eq:pl2}), respectively, but $\shifttacovmat$ is replaced with $\shifttcovmat$. Then, we have for $n \in \natnum$, for $k,l = 1, \dotsc ,p$, $\as$, 
\begin{align*}
    \begin{split}
    \norm{ \frac{\partial^{2}\tilde{l}_{n,N}}{\partial \theta_{k}\partial \theta_{l}}(\theta_{0}) -  \frac{\partial^{2}l_{n,N}}{\partial \theta_{k}\partial \theta_{l}}(\theta_{0})} &\leq \frac{1}{h(n)}\norm{\tr\left(\widetilde{A}_{1,h(n),\theta_{0}}^{kl} - A_{1,h(n),\theta_{0}}^{kl}\right)} \\
    &\quad+ \frac{1}{h(n)}\norm{\langle Z_{(h(n))}, \widetilde{A}_{2,h(n),\theta_{0}}^{kl} - A_{2,h(n),\theta_{0}}^{kl}Z_{(h(n))}\rangle} \\
    & \leq \frac{1}{h(n)}\norm{\tr\left(\widetilde{A}_{1,h(n),\theta_{0}}^{kl} - A_{1,h(n),\theta_{0}}^{kl}\right)} \\
    &\quad + \specnorm{\widetilde{A}_{2,h(n),\theta_{0}}^{kl} - A_{2,h(n),\theta_{0}}^{kl}}\frac{\nnorm{Z_{(h(n))}}^{2}}{h(n)}.
    \end{split}
\end{align*}
We can apply Lemma~\ref{quadraticForm} to the sequence of random matrices $\big(\widetilde{A}_{2,h(n),\theta_{0}}^{kl}\big){}_{n\in \natnum}$ and $\big(A_{2,h(n),\theta_{0}}^{kl}\big){}_{n\in \natnum}$ to conclude under application of Lemma~\ref{eigenbehave} (see also Corollary~\ref{corollary1} and Lemma~\ref{lemmaEigenbehavePartials}) that 
\begin{gather*}
    \specnorm{\widetilde{A}_{2,h(n),\theta_{0}}^{kl} - A_{2,h(n),\theta_{0}}^{kl}}\frac{\nnorm{Z_{(h(n))}}^{2}}{h(n)} \toinp 0.
\end{gather*}
We also have $\as$ 
\begin{gather*}
    \frac{1}{n}\norm{\tr\left(\widetilde{A}_{1,h(n),\theta_{0}}^{kl} - A_{1,h(n),\theta_{0}}^{kl}\right)} \tonormal 0,
\end{gather*}
using the triangular inequality, von Neumann's trace inequality and Lemma~\ref{eigenbehave} (see also Corollary~\ref{corollary1} and Lemma~\ref{lemmaEigenbehavePartials}). Hence, we have shown that for any $k,l = 1, \dotsc, p$
\begin{gather*}
    \norm{ \frac{\partial^{2}\tilde{l}_{n,N}}{\partial \theta_{k}\partial \theta_{l}}(\theta_{0}) -  \frac{\partial^{2}l_{n,N}}{\partial \theta_{k}\partial \theta_{l}}(\theta_{0})} \toinp 0.
\end{gather*}
In addition, we have that for any $k,l = 1, \dotsc, p$, $\as$, the expression
\begin{gather*}
    \norm{\frac{\partial \big(\operatorname{E}\big[\frac{\partial \tilde{l}_{n,N}}{\partial \theta_{l}}(\theta_{0})\; \big \lvert \; S_{(h(n))}\big]\big)}{\partial\theta_{k}} -\frac{\partial \big(\operatorname{E}\big[\frac{\partial l_{n,N}}{\partial \theta_{l}}(\theta_{0})\; \big \lvert \; S_{(h(n))}\big]\big)}{\partial\theta_{k}}}
\end{gather*}
is bounded from above by
\begin{gather*}
    \frac{\big\lvert\operatorname{tr}\big(\widetilde{A}_{1,h(n),\theta_{0}}^{kl}- A_{1,h(n),\theta_{0}}^{kl}\big)\big\rvert
    + \big\lvert\operatorname{tr}\big(\big(\widetilde{A}_{2,h(n),\theta_{0}}^{kl}- A_{2,h(n),\theta_{0}}^{kl}\big)\shifttcovmat\big)\big\rvert}{h(n)},
\end{gather*}
which again, under application of the triangular inequality, von Neumann's trace inequality and Lemma~\ref{eigenbehave} (see also Corollary~\ref{corollary1} and Lemma~\ref{lemmaEigenbehavePartials}), converges to zero $\as$ Hence, we have shown that for any $k,l = 1, \dotsc, p$,
\begin{gather*}
    \norm{\big(J_{\widetilde{G}_{n,N}}(\theta_{0})\big)_{kl} - \big(J_{G_{n,N}}(\theta_{0})\big)_{kl} } \toinp 0,
\end{gather*}
which concludes the proof of (\ref{l2.1}). Now it is shown in \cite{bachoc2014asymptotic} (see Propositions D.$7$ and D.$8$ and also consider the proofs of Propositions $3.2$ and $3.3$), under application of Lemmas~\ref{lemma1}, \ref{lemma3}, \ref{eigenbehave}, \ref{lemmaEigenbehavePartials} and Corollary~\ref{corollary1}, that  
\begin{gather*}
    J_{G_{n,N}}(\theta_{0}) \toinp 2\Lambda,
\end{gather*}
where $\Lambda$ is the $\as$ limit of a sequence of $p \times p$  matrices $\big(H_{h(n)}(\theta_{0})\big){}_{n \in \natnum}$ defined as
\begin{gather*}
     \left\{\left[\frac{1}{2h(n)}\tr\left(\shifttcovmat^{-1}\frac{\partial \shifttcovmat}{\partial \theta_{k}}\shifttcovmat^{-1}\frac{\partial \shifttcovmat}{\partial \theta_{l}}\right)\right]_{1 \leq k,l \leq p} \colon n \in \natnum\right\}.
\end{gather*}
Further, by Assumption~\ref{cov_a4}, it is concluded that the limit $\Lambda$ is such that $\Lambda \succ 0$. But then, we use (\ref{l2.1}) to show that
\begin{gather*}
    J_{\widetilde{G}_{n,N}}(\theta_{0}) \toinp 2\Lambda,
\end{gather*}
as well, which concludes the proof of Lemma~\ref{lemma2}.
\end{proof}

\subsection{Proof of results in Section~\ref{sec:uniformequivof mat}}

\begin{proof}[Proof of Lemma~\ref{eigenbehave}]
We rely on Lemma~\ref{eigenbehave1x1} and a proof is evident.
\end{proof}

\subsection{Proof of results in Section~\ref{MLE_results}}

To simplify the notation, we write $\big(\tilde{\theta}_{n}\big){}_{n \in \natnum} \coloneqq \big(\hat{\theta}_{n}(\tilde{c})\big){}_{n \in \natnum}$.

\begin{proof}[Proof of Proposition~\ref{orderN}]
The statement is verified as a consequence of Lemmas~\ref{eigenbehave},~\ref{lemmalogdet} and~\ref{quadraticForm}.
\end{proof}

\begin{proof}[Proof of Theorem~\ref{thm1}]
Let $N \in \natnum$ be as in Lemma~\ref{eigenbehave} (or Proposition~\ref{orderN}) and define, for any $\omega \in \Omega$, the sequence $\big(\tilde{l}_{n,N}\!\left(\theta\right)\left(\omega\right)\big){}_{n \in \natnum}$ as in (\ref{shiftloglik2}) of Section~\ref{Consasym}. We note that, under the given assumptions of Theorem~\ref{thm1}, the first order partial derivatives with respect to $\theta$ exist for the sequence $\big(\tilde{l}_{n,N}\!\left(\theta\right)\left(\omega\right)\big){}_{n \in \natnum}$. Then, we define the sequence of estimators $\big(\tilde{\theta}_{n,N}\big){}_{n \in \natnum} \coloneqq \big(\tilde{\theta}_{n + N-1}\big){}_{n \in \natnum}$. Therefore $\tilde{\theta}_{n,N}$ minimizes $\tilde{l}_{n,N}\!\left(\theta\right)$ $\as$ for any $n \in \natnum$. To prove that
\begin{gather*}
    \tilde{\theta}_{n} \toinp \theta_{0},
\end{gather*}
it is sufficient to show that
\begin{equation} \label{totoproof}
    \tilde{\theta}_{n,N} \toinp \theta_{0}.
\end{equation}
We consider a similar approach as given in \cite{bachoc2014asymptotic}. As $N$ is fixed, we write for $n \in \natnum$, $h(n) = n + N-1$. 
Under the assumptions of the theorem we have that $\as$
\begin{equation} \label{to1}
    \var{\tilde{l}_{n,N}(\theta)\;\vert\;S_{(h(n))}} \tonormal 0,
\end{equation}
and
\begin{equation} \label{to2}
    \max_{k=1,\dotsc,p}\sup_{\theta \in \Theta}\norm{\frac{\partial}{\partial \theta_{k}}\tilde{l}_{n,N}(\theta)} = O_{\mathbb{P}}(1)\text{ as } n \to \infty.
\end{equation}
To see it, we remark that $\as$ (using Proposition~\ref{orderN}),
\begin{gather*}
    \var{\tilde{l}_{n,N}(\theta)\;\vert\;S_{(h(n))}} = \frac{2}{h(n)^{2}}\tr\bigg(\underbrace{\shiftacovmat^{-1}\shifttcovmat\shiftacovmat^{-1}\shifttcovmat}_{\eqqcolon \widetilde{A}_{h(n),\theta}}\bigg).
\end{gather*}
From here, we can use von Neumann's trace inequality to show that $\as$
\begin{gather*}
    \norm{\tr\left(\widetilde{A}_{h(n),\theta}\right)} \leq h(n) \bigg\rVert\shiftacovmat^{-1}\bigg\rVert_{2}^{2}\bigg\rVert\shifttcovmat\bigg\rVert_{2}^{2}.
\end{gather*}
Now, by Lemma~\ref{eigenbehave} (and Corollary~\ref{corollary1}) we can conclude that there exists a real constant $M_{0} > 0$, such that for any $n \in \natnum$, $\as$, $\mathrm{Var}\big(\tilde{l}_{n,N}(\theta)\;\vert\;S_{(h(n))}\big) \leq M_{0}/h(n)$, which proofs~(\ref{to1}). For (\ref{to2}), we first notice that by Lemma~\ref{lemmaEigenbehavePartials}, there exist constants $M_{1}$, $M_{2} > 0$ (which are independent of $n \in \natnum$, $s_{(n)} \in \mathcal{G}_{n}$ and $\theta \in \Theta$) such that $\as$
\begin{gather*}
    \sup_{\theta \in \Theta}\specnorm{\acovmat^{-1}} < M_{1}\text{ and }\max_{k=1,\dotsc,p}\sup_{\theta \in \Theta}\specnorm{\frac{\partial}{\partial \theta_{k}}\acovmat} \leq M_{2}.
\end{gather*}
Using this result have that $\as$ 
\begin{align} \label{re0}
\begin{split}
\max_{k=1,\dotsc,p}\sup_{\theta \in \Theta}\norm{\frac{\partial}{\partial \theta_{k}}\tilde{l}_{n,N}(\theta)} &= \max_{k=1,\dotsc,p}\sup_{\theta \in \Theta}\bigg \vert \frac{1}{h(n)} \tr\left(\shiftacovmat^{-1}\frac{\partial}{\partial \theta_{k}}\shiftacovmat\right) \\
& \quad - \frac{1}{h(n)}\langle Z_{(h(n))}, \shiftacovmat^{-1}\frac{\partial}{\partial \theta_{k}}\shiftacovmat\shiftacovmat^{-1}Z_{(h(n))}\rangle \bigg \vert\\
                 &\leq M_{1}M_{2} + M_{1}^{2}M_{2}\frac{\nnorm{Z_{(h(n))}}^{2}}{h(n)}. 
\end{split}
\end{align} 
Let $V_{h(n)}$ be a Gauss vector on $(\Omega, \mathcal{F}, \mathbb{P})$, with zero-mean vector and $h(n) \times h(n)$ identity covariance matrix. Then (see also Remark~\ref{remark:conditionaldisti}), for any finite $M > 0$, we have that the probability
\begin{gather*}
    \mathbb{P}\big(M_{1}M_{2} + M_{1}^{2}M_{2}h(n)^{-1}\nnorm{Z_{(h(n))}}^{2} > M \mid S_{(h(n))} = s_{(h(n))}\big)
\end{gather*}
is bounded from above by 
\begin{gather*}
\mathbb{P}\big(M_{1}M_{2}\big(1+h(n)^{-1}\lVert V_{h(n)}\rVert^{2}\big) > M\big).
\end{gather*}
Therefore, $\as$,
\begin{gather*}
    \mathbb{P}\big(M_{1}M_{2} + M_{1}^{2}M_{2}h(n)^{-1}\nnorm{Z_{(h(n))}}^{2} > M \mid S_{(h(n))}\big)
\end{gather*}
is bounded from above by $\mathbb{P}\big(M_{1}M_{2}\big(1+h(n)^{-1}\lVert V_{h(n)}\rVert^{2}\big) > M\big)$ as well. Since $M_{1}M_{2}\big(1+h(n)^{-1}\lVert V_{h(n)}\rVert^{2}\big) = O_{\mathbb{P}}(1)\text{ as } n \to \infty$, (\ref{to2}) is shown. 

Notice further that $\Theta$ is convex, $\theta \mapsto \tilde{l}_{n,N}(\theta)$ is continuously differentiable and by (\ref{re0})
\begin{gather*}
    \sup_{n \in \natnum}\mathbb{E}\bigg[\max_{k=1,\dotsc,p}\sup_{\theta \in \Theta}\norm{\frac{\partial}{\partial \theta_{k}}\tilde{l}_{n,N}(\theta)}\bigg] < \infty.
\end{gather*}
Thus, under application of Corollary $2.2$ of \cite{NeweyWhitney_uniformConvinProb}, with (\ref{to1}) and (\ref{to2}), we can conclude that
\begin{equation} \label{consistency2}
    \sup_{\theta \in \Theta}\norm{\tilde{l}_{n,N}(\theta)- \mathbb{E}\left[\tilde{l}_{n,N}(\theta)\;\vert\;S_{(h(n))}\right]}\toinp 0.
\end{equation}
To continue, we define the sequences of random variables
\begin{gather*}
    \big(D_{h(n), \theta, \theta_{0}}\big)_{n \in \natnum} \coloneqq \big(\mathbb{E}\big[l_{n,N}(\theta)\;|\; S_{(h(n))}\big]- \mathbb{E}\big[l_{n,N}(\theta_{0})\;|\; S_{(h(n))}\big]\big)_{n \in \natnum},\\
    \big(\widetilde{D}_{h(n),\theta, \theta_{0}}\big)_{n \in \natnum} \coloneqq \big(\mathbb{E}\big[\tilde{l}_{n,N}(\theta) \;|\; S_{(h(n))}\big]- \mathbb{E}\big[\tilde{l}_{n,N}(\theta_{0})\;|\; S_{(h(n))}\big]\big)_{n \in \natnum}.
\end{gather*}
For any $n \in \natnum$, we have that $\as$,
\begin{equation*}
\begin{split}
D_{h(n), \theta, \theta_{0}} &= \frac{1}{h(n)}\llog{\ddet{\shiftcovmat}} + \frac{1}{h(n)}\tr\left(\shiftcovmat^{-1}\shifttcovmat\right) \\
& \quad - \frac{1}{h(n)}\llog{\ddet{\shifttcovmat}} - \frac{1}{h(n)}\tr\left(\shifttcovmat^{-1}\shifttcovmat\right).
\end{split}
\end{equation*}
Similarly, For any $n \in \natnum$, we have that $\as$
\begin{equation*}
    \begin{split}
    \widetilde{D}_{h(n),\theta, \theta_{0}} &= \frac{1}{h(n)}\llog{\ddet{\shiftacovmat}} + \frac{1}{h(n)}\tr\left(\shiftacovmat^{-1}\shifttcovmat\right) \\
    & \quad - \frac{1}{h(n)}\llog{\ddet{\shifttacovmat}} - \frac{1}{h(n)}\tr\left(\shifttacovmat^{-1}\shifttcovmat\right).
\end{split}
\end{equation*}
Notice that because of~(\ref{consistency2}) we have that 
\begin{equation*} 
    \sup_{\theta \in \Theta}\norm{\left(\tilde{l}_{n,N}(\theta) - \tilde{l}_{n,N}(\theta_{0})\right) - \widetilde{D}_{h(n),\theta, \theta_{0}}}\toinp 0.
\end{equation*}
Further, it is shown in \cite{bachoc2014asymptotic} (see the proof of Proposition $3.1$) that under application of Lemma~\ref{eigenbehave1x1} there exists some constant $B > 0$ (which does not depend on $n \in \natnum$) such that $\as$ 
\begin{gather*} 
    D_{h(n), \theta, \theta_{0}} \geq B\underbrace{\frac{1}{h(n)}\sum_{i,j=1}^{h(n)}\bigg(c_{\theta}(S_{i} - S_{j}) - c_{\theta_{0}}(S_{i} - S_{j})\bigg)^{2}}_{\eqqcolon D_{2, h(n),\theta, \theta_{0}}}.
\end{gather*}
Under application of Lemmas~\ref{lemma1}, \ref{lemma3}, \ref{eigenbehave1x1} and Corollary~\ref{corollary1}, it is then shown in the proof of Proposition $3.1$ of \cite{bachoc2014asymptotic} that either, if $\tau = 0$, $D_{2, h(n),\theta, \theta_{0}}$ is deterministic and we have that
\begin{gather*}
\sup_{\theta \in \Theta}\norm{D_{2, h(n),\theta, \theta_{0}} - D_{\infty, \theta, \theta_{0}}} \tonormal 0, 
\end{gather*}
where the limit is given by $D_{\infty, \theta, \theta_{0}} = \sum_{z\in\mathbb{Z}^{d}}\big(c_{\theta}(z) - c_{\theta_{0}}(z)\big)^{2}$. Or $\tau > 0$ and it is concluded that
\begin{gather*}
\sup_{\theta \in \Theta}\norm{D_{2, h(n),\theta, \theta_{0}} - D_{\infty, \theta, \theta_{0}}} \toinp 0, 
\end{gather*}
where in this case
\begin{gather*}
    D_{\infty, \theta, \theta_{0}} = \int_{D_{\tau}}\!\big(c_{\theta}(s) - c_{\theta_{0}}(s)\big)^{2}f(s)ds + \big(c_{\theta}(0) - c_{\theta_{0}}(0)\big)^{2}.
\end{gather*}
Notice that because of the assumption that $(X_{i})_{i \in \natnum}$ is independent with common law that has a strictly positive probability density function, the function $f$ is strictly positive almost everywhere with respect to the Lebesgue measure on $D_{\tau}$ (see the end of the proof of Proposition~$3.1$ in \cite{bachoc2014asymptotic}). In either case, we can thus conclude that 
\begin{gather*}
\sup_{\theta \in \Theta}\norm{D_{2, h(n),\theta, \theta_{0}} - D_{\infty, \theta, \theta_{0}}} \toinp 0, 
\end{gather*}
where for any $\alpha > 0$, because of Assumption~\ref{cov_a3}, $\inf_{\theta \colon \norm{\theta - \theta_{0}}\geq \alpha}D_{\infty, \theta, \theta_{0}} > 0$, and the limit $D_{\infty, \theta, \theta_{0}}$ is deterministic. We now want to show that there exists some $N_{2} \geq N$ such that for any $n \geq N_{2}$, for any $\theta \in \Theta$, $\as$,
\begin{equation} \label{consistency3}
    \widetilde{D}_{h(n),\theta,\theta_{0}} \geq B D_{2, h(n),\theta, \theta_{0}},
\end{equation}
as well. In this case, with $D_{2, h(n),\theta, \theta_{0}}$ a random function on $\Omega$ and $D_{\infty, \theta, \theta_{0}}$ a deterministic function of $\theta \in \Theta$, we would have for any fixed $\tau \geq 0$, and for any given $\alpha > 0$,   
\begin{gather*}
    \sup_{\theta \in \Theta}\norm{D_{2, h(n), \theta, \theta_{0}} - D_{\infty, \theta, \theta_{0}}} \toinp 0,\\
    \inf_{\theta \colon \norm{\theta - \theta_{0}} \geq \alpha}D_{\infty, \theta, \theta_{0}} > D_{\infty, \theta_{0}, \theta_{0}} = 0,
\end{gather*}
where the sequence of estimators $\big(\tilde{\theta}_{n,N}\big){}_{n \in \natnum}$, is such that for $n \geq N_{2}$, $\as$, 
\begin{align*}
    &D_{2, h(n), \tilde{\theta}_{n,N}, \theta_{0}} \\
    &= D_{2, h(n), \tilde{\theta}_{n,N}, \theta_{0}} - \frac{1}{B}\big(\tilde{l}_{n,N}\big(\tilde{\theta}_{n,N}\big) - \tilde{l}_{n,N}\big(\theta_{0}\big)\big) + \frac{1}{B}\big(\tilde{l}_{n,N}\big(\tilde{\theta}_{n,N}\big) - \tilde{l}_{n,N}\big(\theta_{0}\big)\big) \\
    &\leq \frac{\widetilde{D}_{h(n), \tilde{\theta}_{n,N}, \theta_{0}}}{B}  - \frac{1}{B}\big(\tilde{l}_{n,N}\big(\tilde{\theta}_{n,N}\big) - \tilde{l}_{n,N}\big(\theta_{0}\big)\big)  + \frac{1}{B}\big(\tilde{l}_{n,N}\big(\tilde{\theta}_{n,N}\big) - \tilde{l}_{n,N}\big(\theta_{0}\big)\big) \\
    &\leq D_{2, h(n), \theta_{0}, \theta_{0}}  + \underbrace{\frac{1}{B}\sup_{\theta \in \Theta}\norm{\big(\tilde{l}_{n,N}(\theta) - \tilde{l}_{n,N}(\theta_{0})\big) - \widetilde{D}_{h(n),\theta, \theta_{0}}}}_{\toinp \; 0},
\end{align*}
and we can conclude the proof of Theorem~\ref{thm1}, using Theorem~$5.7$ of \cite{van2000asymptotic}. Hence, it remains to show~(\ref{consistency3}). 
We write $\as$, 
\begin{align*}
    \begin{split}
    \norm{\widetilde{D}_{h(n), \theta, \theta_{0}} - D_{h(n), \theta, \theta_{0}}} &\leq \norm{\widetilde{A}_{1, h(n), \theta, \theta_{0}} - A_{1, h(n), \theta, \theta_{0}}} \\
    & \quad + \norm{\widetilde{A}_{2, h(n), \theta, \theta_{0}} - A_{2, h(n), \theta, \theta_{0}}}\\
    &\quad+ \norm{\widetilde{A}_{3, h(n), \theta, \theta_{0}} - A_{3, h(n), \theta, \theta_{0}}}, 
    \end{split}
\end{align*}
where
\begin{align*}
    \begin{split}
    \widetilde{A}_{1, h(n), \theta, \theta_{0}} - A_{1, h(n), \theta, \theta_{0}} &= \frac{1}{h(n)}\llog{\ddet{\shifttcovmat\shiftcovmat^{-1}}} \\
    &\quad- \frac{1}{h(n)}\llog{\ddet{\shifttacovmat\shiftacovmat^{-1}}},
    \end{split}
\end{align*}
\begin{align*}
    \widetilde{A}_{2, h(n), \theta, \theta_{0}} - A_{2, h(n), \theta, \theta_{0}} &= \frac{1}{h(n)}\tr\left(\left[\shiftacovmat^{-1}-\shiftcovmat^{-1}\right]\shifttcovmat\right),
\end{align*}
and
\begin{align*}
    \widetilde{A}_{3, h(n), \theta, \theta_{0}} - A_{3, h(n), \theta, \theta_{0}} &= \frac{1}{h(n)}\tr\left(\left[\shifttacovmat^{-1} - \shifttcovmat^{-1}\right]\shifttcovmat\right).
\end{align*}
By Lemma~\ref{eigenbehave}, Corollary~\ref{corollary1} and Lemma~\ref{lemmalogdet}, we already conclude that $\as$ 
\begin{gather*}
   \norm{\widetilde{A}_{1, h(n), \theta, \theta_{0}} - A_{1, h(n), \theta, \theta_{0}}} 
\end{gather*}
converges to zero uniformly in $\theta \in \Theta$ as $n \xrightarrow{} \infty$. Further, we can conclude that $\as$ 
\begin{align*}
    \big\lvert\widetilde{A}_{2, h(n), \theta, \theta_{0}} - A_{2, h(n), \theta, \theta_{0}}\big\rvert &\leq \big\lVert\shiftacovmat^{-1}-\shiftcovmat^{-1}\big\rVert_{2}\big\lVert\shifttcovmat\big\rVert_{2},
\end{align*}
and thus, since by Lemma~\ref{eigenbehave} $\as$ $\lVert\shifttcovmat\rVert_{2}$, $\lVert\shiftacovmat\rVert_{2}$ and $\lVert\shiftcovmat\rVert_{2}$ are finite, uniformly in $n \in \natnum$ and $\theta \in \Theta$, and $\as$
\begin{gather*}
    \specnorm{\shiftacovmat^{-1}-\shiftcovmat^{-1}} \tonormal 0,
\end{gather*}
uniformly in $\theta \in \Theta$, by application of Corollary~\ref{corollary1}, we can also see that $\as$ 
\begin{gather*}
    \norm{\widetilde{A}_{2, h(n), \theta, \theta_{0}} - A_{2, h(n), \theta, \theta_{0}}}
\end{gather*}
converges to zero as $n \xrightarrow{} \infty$, uniformly in $\theta \in \Theta$. Using a similar argument we can also show that $\as$ the term $\norm{\widetilde{A}_{3, h(n), \theta, \theta_{0}} - A_{3, h(n), \theta, \theta_{0}}}$ converges to zero as $n \to \infty$, uniformly in $\theta \in \Theta$. Hence, we have shown that $\as$
\begin{gather*}
    \norm{\widetilde{D}_{h(n), \theta, \theta_{0}} - D_{h(n), \theta, \theta_{0}}} \tonormal 0, 
\end{gather*}
uniformly in $\theta \in \Theta$ and we can argue that $\as$ there exists some $N_{2} \geq N$ such that for all $n \geq N_{2}$, $\widetilde{D}_{h(n), \theta, \theta_{0}} \geq B D_{2, h(n), \theta, \theta_{0}}$ on $\Omega$, which shows~(\ref{consistency3}). Therefore, we have that 
\begin{gather*}
    \tilde{\theta}_{n,N} \toinp \theta_{0},
\end{gather*}
which concludes the proof.
\end{proof}

\begin{proof}[Proof of Corollary~\ref{corollarythm1}]
This follows from Theorem~\ref{thm1} when we define, for any $\theta \in \Theta$ and $m \in \natnum$, $\tilde{c}_{m,\theta}(s) \coloneqq c_{\theta}(s)$, for all $s \in \mathbb{R}^{d}$.
\end{proof}

\begin{proof}[Proof of Theorem~\ref{thm2}]
Let $N \in \natnum$ be as in Proposition~\ref{orderN} and define, for any $\omega \in \Omega$, the sequences of functions $\big(l_{n,N}\!\left(\theta\right)\left(\omega\right) \big){}_{n \in \natnum}$ and $\big(\tilde{l}_{n,N}\!\left(\theta\right)\left(\omega\right)\big){}_{n \in \natnum}$ as in (\ref{shiftloglik1}) and (\ref{shiftloglik2}) of Section~\ref{Consasym}, respectively. From the proof of Theorem~\ref{thm1} we know that sequence of estimators $\big(\tilde{\theta}_{n,N}\big){}_{n \in \natnum} \coloneqq \big(\tilde{\theta}_{n + N-1}\big){}_{n \in \natnum}$, is such that 
\begin{gather*} 
    \tilde{\theta}_{n,N} \toinp \theta_{0}.
\end{gather*}
Define $\big\{(G_{n,N}\left(\theta\right)){}_{n \in \natnum}\colon \theta \in \Theta\big\}$ and $\big\{(\widetilde{G}_{n,N}\left(\theta\right)){}_{n \in \natnum}\colon \theta \in \Theta\big\}$ as in (\ref{cFun2}) and (\ref{cFun}) respectively. For $n \in \natnum$ we set $h(n) = n + N-1$. We have for $k = 1, \dotsc, p$, $\as$, for $n \in \natnum$, 
\begin{align*}
    \begin{split}
    \tilde{c}_{k,n,N}(\theta_{0}) &= \frac{\partial \tilde{l}_{n,N}}{\partial \theta_{k}}(\theta) - \E{\frac{\partial \tilde{l}_{n,N}}{\partial \theta_{k}}(\theta)\; \middle | \; S_{(h(n))}} \\
    &= \frac{1}{h(n)}\tr\big(\underbrace{\shifttacovmat^{-1}\frac{\partial \shifttacovmat}{\partial \theta_{k}}\shifttacovmat^{-1}\shifttcovmat}_{\eqqcolon \widetilde{M}_{k,h(n)}}\big) \\
    & \quad+ \frac{1}{h(n)}\langle Z_{(h(n))} ,\underbrace{-\shifttacovmat^{-1}\frac{\partial \shifttacovmat}{\partial \theta_{k}}\shifttacovmat^{-1}}_{\eqqcolon \widetilde{N}_{k,h(n)}}Z_{(h(n))} \rangle,
\end{split}
\end{align*}
where, by Lemma~\ref{eigenbehave} (see also Corollary~\ref{corollary1} and Lemma~\ref{lemmaEigenbehavePartials}) $\as$, for $n \in \natnum$, $\lVert\widetilde{M}_{k,h(n)}\rVert_{2}$ and $\lVert\widetilde{N}_{k,h(n)}\rVert_{2}$ are finite, uniformly in $\theta \in \Theta$. Further, notice that $\as$
\begin{gather*}
    \tr\left(\widetilde{M}_{k,h(n)} + \widetilde{N}_{k,h(n)}\shifttcovmat\right) = 0 \quad \forall\, k \in \{1, \dotsc,p\}.
\end{gather*}
From the proof of Lemma~\ref{lemma2}, we already know that $\as$, $H_{h(n)}(\theta_{0}) \tonormal \Lambda$, where for any $k,l = 1, \dotsc,p$, for any $n \in \natnum$,
\begin{gather*}
    \left[H_{h(n)}(\theta_{0})\right]_{kl} = \frac{1}{2h(n)}\tr\left(\shifttcovmat^{-1}\frac{\partial \shifttcovmat}{\partial \theta_{k}}\shifttcovmat^{-1}\frac{\partial \shifttcovmat}{\partial \theta_{l}}\right).
\end{gather*}
Now, if we define, on $\left(\Omega, \mathcal{F}, \mathbb{P}\right)$, the sequence of random $p \times p$ matrices
\begin{gather*}
    \bigg\{\underbrace{\frac{1}{2h(n)}\left[\tr\left(\widetilde{N}_{k,h(n)}\shifttcovmat \widetilde{N}_{\left\{l,h(n)\right\}}\shifttcovmat\right)\right]_{\{1 \leq k,l \leq p\}}}_{\eqqcolon \widetilde{H}_{h(n)}(\theta_{0})}\colon n \in \natnum\bigg\} 
\end{gather*}
we also have that, for any $k,l = 1, \dotsc p$, $\as$, $\big\lvert\big[\widetilde{H}_{h(n)}(\theta_{0})\big]_{kl}-\Sigma_{kl}\big\rvert \tonormal 0$. This follows from the fact that for any $k,l = 1, \dotsc p$, we have that $\as$
\begin{gather*}
    \norm{\left[\widetilde{H}_{h(n)}(\theta_{0})\right]_{kl}- \left[H_{h(n)}(\theta_{0})\right]_{kl}} \leq \frac{1}{2h(n)}\norm{\tr\left(\widetilde{B}_{h(n)}^{kl} - B_{h(n)}^{kl}\right)},
\end{gather*} 
where 
\begin{align*}
   \widetilde{B}_{h(n)}^{kl} &= \shifttacovmat^{-1}\frac{\partial \shifttacovmat}{\partial \theta_{k}}\shifttacovmat^{-1}\shifttcovmat\shifttacovmat^{-1}\frac{\partial \shifttacovmat}{\partial \theta_{l}}\shifttacovmat^{-1}\shifttcovmat,
\end{align*}
and
\begin{align*}
    B_{h(n)}^{kl} &= \shifttcovmat^{-1}\frac{\partial \shifttcovmat}{\partial \theta_{k}}\shifttcovmat^{-1}\shifttcovmat\shifttcovmat^{-1}\frac{\partial \shifttcovmat}{\partial \theta_{l}}\shifttcovmat^{-1}\shifttcovmat,
\end{align*}
and again, under application of the triangular inequality, von Neumann's trace inequality and Lemma~\ref{eigenbehave} (see also Corollary~\ref{corollary1} and Lemma~\ref{lemmaEigenbehavePartials}) we thus have that $\as$ $\lVert\widetilde{H}_{h(n)}(\theta_{0}) - H_{h(n)}(\theta_{0})\rVert_{2}\tonormal 0$. But $\Lambda$ is the $\as$ limit of $\big\{H_{h(n)}(\theta_{0}) \colon n \in \natnum\big\}$ and hence we conclude that $\Lambda$ is also the $\as$ limit of $\big\{\widetilde{H}_{h(n)}(\theta_{0}) \colon n \in \natnum\big\}$. Then, we can apply Proposition D.9 of \cite{bachoc2014asymptotic} to conclude that
\begin{gather*}
    h(n)^{1/2}\widetilde{G}_{n,N}(\theta_{0}) \toind \mathcal{N}\left(0, 4\Lambda\right).
\end{gather*}
Notice that because the family $\big\{(\tilde{c}_{m,\theta}) \colon \theta \in \Theta\big\}$ satisfies Assumption~\ref{uniformly_approx_covfun}, we have that for fixed $\omega \in \Omega$, $\theta \mapsto \widetilde{G}_{n,N}(\omega, \theta)$ is twice differentiable in $\theta$ and we can argue exactly as in the proof of Theorem~\ref{thm1} to conclude that the sequence 
\begin{gather*}
    \left(\sup_{\theta \in \Theta}\max_{\{1 \leq k,l,m \leq p\}}\norm{\frac{\partial \left(\frac{\partial \widetilde{G}_{m,n,N}}{\partial \theta_{l}}\right)}{\partial \theta_{k}}(\theta)}\right)_{n \in \natnum},
\end{gather*}
is bounded in probability $\mathbb{P}$.
In addition, by Lemma~\ref{lemma2}, we also have that
\begin{gather*}
   \specnorm{J_{\widetilde{G}_{n,N}}(\theta_{0}) - J_{G_{n,N}}(\theta_{0}) } \toinp 0.
\end{gather*}
Finally, the sequence of estimators $\big(\tilde{\theta}_{n,N}\big){}_{n \in \natnum}$ is consistent and such that
\begin{gather*}
    \mathbb{P}\left(\widetilde{G}_{n,N}\big(\tilde{\theta}_{n,N}\big) = 0\right) \tonormal 1.
\end{gather*}
Thus we conclude, using for example Proposition D.10 in \cite{bachoc2014asymptotic} that
\begin{gather*}
    h(n)^{1/2}\left(\tilde{\theta}_{n,N} - \theta_{0}\right) \toind \mathcal{N}(0,\Lambda^{-1}).
\end{gather*}
Since $N$ was fixed, we can conclude that 
\begin{gather*}
    n^{1/2}\left(\tilde{\theta}_{n} - \theta_{0}\right) \toind \mathcal{N}(0,\Lambda^{-1}),
\end{gather*}
as well. 
\end{proof}

\begin{proof}[Proof of Corollary~\ref{corollarythm2}]
The result follows from Theorem~\ref{thm2}, a proof is evident when we define the family $\big\{\big(\tilde{c}_{m,\theta}\big) \colon \theta \in \Theta\big\}$ as in the proof of Corollary~\ref{corollarythm1}.
\end{proof}

\subsection{Proof of results in Appendix~\ref{sec_application1}} \label{sec_applications}

\begin{proof}[Proof of Theorem~\ref{consistencyradial}]
The proof is similar to the proof of Theorem~\ref{thm1}.
\end{proof}

\begin{proof}[Proof of Theorem~\ref{asymnormalradial}]
The proof is similar to the proof of Theorem~\ref{thm2}.
\end{proof}

\subsection{Proof of results in Section~\ref{sec_application}}
Since $\nu$ and $\kappa$ are assumed to be known, we put $c_{\nu, \kappa} \coloneqq \operatorname{B}(2\kappa,\nu + 1)$. We define the function $f_{\nu, \kappa}(r,u) \coloneqq u(u^{2} - r^{2})^{\kappa-1}(1-u)^{\nu}$, $(r,u) \in [0,1] \times [0,1]$. We recall that for $r = 0$,  
\begin{gather*}
    c_{\nu,\kappa} = \int_{0}^{1}f_{\nu,\kappa}(0,u)du.
\end{gather*}

\begin{proof}[Proof of Proposition~\ref{check1}]
We have already seen that for any $\theta \in \Theta$, given known $\kappa > 0$ and $\nu \geq (d+1)/2 + \kappa$, $\phi_{\theta}$ is continuous on $\mathbb{R}_{+}$. Further, for any $\theta \in \Theta$, $\phi_{\theta}$ has compact support $S_{\theta} = [0,\beta] \subset [0,\beta_{\max}]$ and since $\kappa > 4$, we can also see that for any $\delta > 0$, for any $t \in \mathbb{R}_{+}$, $\phi_{\theta}(t+\delta) \leq \phi_{\theta}(t)$ and thus $\phi_{\theta}(t) \leq \phi_{\theta}(0) = \sigma^{2}$, which implies that $\supnorm{\phi_{\theta}} \leq \sigma^{2}_{\max}$. Hence, with $C \coloneqq \beta_{\max}$ and $L \coloneqq \sigma^{2}_{\max}$, $C$ and $L$ are independent of $\theta \in \Theta$ and hence we can conclude that \ref{6.1.1} of Assumption~\ref{cov_a2_iso} is satisfied with $\mathcal{B}_{\text{C}}(\mathbb{R}_{+};S_{\theta})$ replaced with $\mathcal{C}_{\text{C}}(\mathbb{R}_{+};S_{\theta})$. It is now sufficient to show that for any $\theta \in \Theta$, for any $q = 1,2,3$, $i_{1}, \dotsc, i_{q} \in \{1, \dotsc,p\}$, there exist constants $C_{\theta}(i_{1}, \dotsc, i_{q})$, $L_{\theta}(i_{1}, \dotsc, i_{q})< \infty$, such that
\begin{equation} \label{toshow1}
\frac{\partial^{q} \phi_{\theta}}{\partial \theta_{i_1}\cdots\partial \theta_{i_q}} \in \mathcal{C}_{\text{C}}(\mathbb{R}_{+};S_{\theta}(i_{1}, \dotsc, i_{q})),
\end{equation}  
where
\begin{align} \label{toshow2}
    \begin{split}
    S_{\theta}(i_{1}, \dotsc, i_{q}) \subset \left[0,C_{\theta}(i_{1}, \dotsc, i_{q})\right] \subset \left[0,C(i_{1}, \dotsc, i_{q})\right], \\
    \supnorm{\frac{\partial^{q} \phi_{\theta}}{\partial \theta_{i_1}\cdots\partial \theta_{i_q}}} \leq L_{\theta}(i_{1}, \dotsc, i_{q}) \leq L(i_{1}, \dotsc, i_{q}),
    \end{split}
\end{align}
with 
$C(i_{1}, \dotsc, i_{q})$, $L(i_{1}, \dotsc, i_{q}) < \infty$, independent of $\theta \in \Theta$. This means that in general we need to check the above condition for $2 + 2^2 + 2^3 = 14$ partial derivatives. Let us first focus on the partial derivatives with respect to the range parameter $\beta \in \left[\beta_{\min}, \beta_{\max}\right]$. 
For $r \in [0,1]$ we write
\begin{gather*}
    \phi_{\nu, \kappa}(r) = c_{\nu,\kappa}^{-1}\int_{a(r)}^{b(r)}f_{\nu,\kappa}(r,u)du,
\end{gather*}
where $[0,1] \ni r \mapsto a(r) = r$ and $[0,1] \ni r \mapsto b(r) \equiv 1$ are continuously differentiable on $[0,1]$. To simplify the notation we will put $f_{\nu,\kappa} \coloneqq f$. Then, since $f \colon [0,1] \times [0,1] \to \mathbb{R}$ is continuous and for any $u \in [0,1]$, since $\kappa > 2$,
\begin{gather*}
    \frac{\partial f}{\partial r}(r,u) = -2r(\kappa-1)u\left(u^{2}-r^{2}\right)^{\kappa-2}(1-u)^{\nu},
\end{gather*}
exists and is continuous on the rectangle $[0,1] \times [0,1]$, we can conclude, using the general Leibniz integral rule, that $[0,1] \ni r \mapsto \frac{d\phi_{\nu, \kappa}}{dr}\left(r\right)$ is continuous and given by
\begin{align*}
    \begin{split}
    \frac{d\phi_{\nu, \kappa}}{dr}\left(r\right) &= \underbrace{f(r,b(r))\frac{db}{dr}(r)}_{=0} - \underbrace{f(r,a(r))\frac{da}{dr}(r)}_{=0} \\ 
    &\quad -2r(\kappa-1)c_{\nu,\kappa}^{-1}\underbrace{\int_{a(r)}^{b(r)}u\left(u^{2}-r^{2}\right)^{\kappa-2}(1-u)^{\nu}du}_{= c_{\nu,\kappa-1}\phi_{\nu, \kappa -1}(r)}.
    \end{split}
\end{align*}
Hence, for $t \in [0,\beta]$
\begin{equation} \label{partialbeta}
    \frac{\partial \phi_{\theta}}{\partial \beta}(t) = -\sigma^{2}\frac{t}{\beta^{2}}\frac{d\phi_{\nu, \kappa}}{dr}\left(\frac{t}{\beta}\right) = \frac{2t^{2}(\kappa-1)}{\beta^{3}}\frac{c_{\nu,\kappa-1}}{c_{\nu,\kappa}}\sigma^{2}\phi_{\nu, \kappa -1}\left(\frac{t}{\beta}\right),
\end{equation}  
exists and is continuous as a function of $t$. But clearly, as for $r \in [1, \infty)$ $\phi_{\nu, \kappa}$ is zero, we have that for $t \in [\beta, \infty)$, $\frac{\partial \phi_{\theta}}{\partial \beta}(t)$ exists as well and is continuous as a function of $t$. Hence, for any $t \in \mathbb{R}_{+}$, $\frac{\partial \phi_{\theta}}{\partial \beta}(t)$ exists, is given by (\ref{partialbeta}) and is continuous as a function of $t$. Thus, by monotonicity of $\phi_{\nu, \kappa -1}$ we define
\begin{gather*}
    L_{\theta}(2) \coloneqq \frac{2\beta^{2}(\kappa-1)}{\beta^{3}}\frac{c_{\nu,\kappa-1}}{c_{\nu,\kappa}}\sigma^{2} = \frac{2(\kappa-1)}{\beta}\frac{c_{\nu,\kappa-1}}{c_{\nu,\kappa}}\sigma^{2}\text{ and }C_{\theta}(2) \coloneqq \beta,
\end{gather*}
and have that 
\begin{gather*}
    \frac{\partial \phi_{\theta}}{\partial \beta} \in \mathcal{C}_{\text{C}}(\mathbb{R}_{+};[0,C_{\theta}(2)])\text{, } \supnorm{\frac{\partial \phi_{\theta}}{\partial \beta}} \leq L_{\theta}(2).
\end{gather*}
Further, we find
\begin{gather*}
    \sup_{\theta \in \Theta}L_{\theta}(2) \leq \underbrace{\frac{2(\kappa-1)}{\beta_{\min}}\frac{c_{\nu,\kappa-1}}{c_{\nu,\kappa}}\sigma_{\sup}^{2}}_{\eqqcolon L(2)}\text{ and }\sup_{\theta \in \Theta}C_{\theta}(2) \leq \underbrace{\beta_{\max}}_{\eqqcolon C(2)},
\end{gather*}
where $L(2)$ and $C(2)$ do not depend on $\theta \in \Theta$. Since we have assumed that $\kappa > 4$, we can now repeat the arguments, which led to (\ref{partialbeta}), for another two times, and conclude that for any $t \in \mathbb{R}_{+}$, $\frac{\partial^{2} \phi_{\theta}}{\partial \beta^{2}}(t)$, and $\frac{\partial^{3} \phi_{\theta}}{\partial \beta^{3}}(t)$ exits as well, are given by
\begin{equation} \label{partialsbeta2}
    \frac{\partial^{2} \phi_{\theta}}{\partial \beta^{2}}(t) =  q_{1}(t),
\end{equation}
with
\begin{align*}
    \begin{split}
    q_{1}(t) &= \frac{4t^{4}(\kappa-1)(\kappa-2)}{\beta^{6}}\frac{c_{\nu,\kappa-2}}{c_{\nu,\kappa}}\sigma^{2}\phi_{\nu,\kappa-2}\left(\frac{t}{\beta}\right) \\
    & \quad - \frac{6t^{2}(\kappa-1)}{\beta^{4}}\frac{c_{\nu,\kappa-1}}{c_{\nu,\kappa}}\sigma^{2}\phi_{\nu,\kappa-1}\left(\frac{t}{\beta}\right)
    \end{split}
\end{align*}
and $\frac{\partial^{3} \phi_{\theta}}{\partial \beta^{3}}(t) = q_{2}(t)$, with 
\begin{equation*}
    \begin{split}
    q_{2}(t) &= \frac{8t^{6}(\kappa-1)(\kappa-2)(\kappa-3)}{\beta^{9}}\frac{c_{\nu,\kappa-3}}{c_{\nu,\kappa}}\sigma^{2}\phi_{\nu, \kappa-3}\left(\frac{t}{\beta}\right) \\
    & \quad+ \frac{24 t^{2}(\kappa-1)}{\beta^{5}}\frac{c_{\nu,\kappa-1}}{c_{\nu,\kappa}}\sigma^{2}\phi_{\nu,\kappa-1}\left(\frac{t}{\beta}\right) \\
    & \quad - \frac{36t^{4}(\kappa-1)(\kappa-2)}{\beta^{7}}\frac{c_{\nu,\kappa-2}}{c_{\nu,\kappa}}\sigma^{2}\phi_{\nu,\kappa-2}\left(\frac{t}{\beta}\right),
    \end{split}
\end{equation*}
and are both continuous as a function of $t \in \mathbb{R}_{+}$. Therefore, since $\phi_{\nu, \kappa -1}$, $\phi_{\nu, \kappa -2}$ and $\phi_{\nu, \kappa -3}$ are non-negative and monotonously decreasing, we can define
\begin{gather*}
    L_{\theta}(2,2) \coloneqq \frac{4\beta^{4}(\kappa-1)(\kappa-2)}{\beta^{6}}\frac{c_{\nu,\kappa-2}}{c_{\nu,\kappa}}\sigma^{2},
\end{gather*}
and
\begin{gather*}
    L_{\theta}(2,2,2) \coloneqq \frac{8\beta^{6}(\kappa-1)(\kappa-2)(\kappa-3)}{\beta^{9}}\frac{c_{\nu,\kappa-3}}{c_{\nu,\kappa}}\sigma^{2} + \frac{24 \beta^{2}(\kappa-1)}{\beta^{5}}\frac{c_{\nu,\kappa-1}}{c_{\nu,\kappa}}\sigma^{2},
\end{gather*}
as well as $C_{\theta}(2,2) = C_{\theta}(2,2,2)  \coloneqq \beta$, and have that
\begin{gather*}
    \frac{\partial^{2} \phi_{\theta}}{\partial \beta^{2}} \in \mathcal{C}_{\text{C}}(\mathbb{R}_{+};[0,C_{\theta}(2,2)]), \quad \supnorm{\frac{\partial^{2} \phi_{\theta}}{\partial \beta^{2}}} \leq L_{\theta}(2,2),
\end{gather*}
\begin{gather*}
    \frac{\partial^{3} \phi_{\theta}}{\partial \beta^{3}}\in \mathcal{C}_{\text{C}}(\mathbb{R}_{+};[0,C_{\theta}(2,2,2)]), \quad \supnorm{\frac{\partial^{3} \phi_{\theta}}{\partial \beta^{3}}} \leq L_{\theta}(2,2,2).
\end{gather*}
Then, we also find
\begin{gather*}
    \sup_{\theta \in \Theta}L_{\theta}(2,2) \leq \frac{4(\kappa-1)(\kappa-2)}{\beta_{\min}^{2}}\frac{c_{\nu,\kappa-2}}{c_{\nu,\kappa}}\sigma_{\sup}^{2} \eqqcolon L(2,2),
\end{gather*}
as well as
\begin{gather*}
    \sup_{\theta \in \Theta}L_{\theta}(2,2,2) \leq \frac{8(\kappa-1)}{\beta_{\min}^{3}}\frac{\sigma_{\sup}^{2}}{c_{\nu,\kappa}}\left((\kappa-2)(\kappa-3)c_{\nu,\kappa-3} + 3c_{\nu,\kappa-1}\right) \eqqcolon L(2,2,2),
\end{gather*}
and $C(2,2) = C(2,2,2) \coloneqq \beta_{\max} = \sup_{\theta \in \Theta}C_{\theta}(2,2) = \sup_{\theta \in \Theta}C_{\theta}(2,2,2)$, where $L(2,2)$, $L(2,2,2)$ and $C(2,2)$ and $C(2,2,2)$ do not depend on $\theta \in \Theta$. This then shows that the partial derivatives of $\phi_{\theta}$ with respect to the range parameter $\beta$ exist up to order three and are continuous on $\mathbb{R}_{+}$ with uniform bounds that do not depend on $\theta \in \Theta$ and compact supports that are subsets of $\left[0,\beta_{\max}\right]$. Let us now focus on the partial derivatives with respect to $\sigma^{2}$. We can readily see that for $t \in \mathbb{R}_{+}$,
\begin{equation} \label{partialsigma}
    \frac{\partial \phi_{\theta}}{\partial \sigma^{2}}(t) = \phi_{\nu,\kappa}\left(\frac{t}{\beta}\right),
\end{equation}
and thus with $S_{\theta}(1) = \left[0,\beta\right]$ and $L_{\theta}(1) = 1$ we can choose $L(1) = 1$ and $C(1) = \beta_{\max}$ such that (\ref{toshow1}) and (\ref{toshow2}) are satisfied. Notice that for any $t \in \mathbb{R}_{+}$, both $\frac{\partial^{2}\phi_{\theta}}{\partial \left(\sigma^{2}\right)^{2}}(t)$ and $\frac{\partial^{3}\phi_{\theta}}{\partial \left(\sigma^{2}\right)^{3}}(t)$ are zero. Thus, the existence of the desired constants $L_{\theta}(2,2)$, $C_{\theta}(2,2)$ and $L(2,2)$, $C(2,2)$ for $\frac{\partial^{2}\phi_{\theta}}{\partial \left(\sigma^{2}\right)^{2}}(t)$ and $L_{\theta}(2,2,2)$, $C_{\theta}(2,2,2)$ and $L(2,2,2)$, $C(2,2,2)$ for $\frac{\partial^{3}\phi_{\theta}}{\partial \left(\sigma^{2}\right)^{3}}$, such that (\ref{toshow1}) and (\ref{toshow2}) is satisfied, is clear. Let us now consider the mixed partial derivatives. Using (\ref{partialbeta}) and (\ref{partialsigma}), we have
\begin{gather*}
    \frac{\partial^{2}\phi_{\theta}}{\partial \sigma^{2} \partial \beta}(t) = \frac{\partial^{2}\phi_{\theta}}{\partial \beta \partial \sigma^{2}}(t) = \frac{2t^{2}(\kappa-1)}{\beta^{3}}\frac{c_{\nu,\kappa-1}}{c_{\nu,\kappa}}\phi_{\nu, \kappa -1}\left(\frac{t}{\beta}\right),
\end{gather*}
and thus the existence of constants $L_{\theta}(1,2) = L_{\theta}(2,1)$, $C_{\theta}(1,2) = C_{\theta}(2,1)$ and $L(1,2) = L(2,1)$, $C(1,2) = C(2,1)$ for $\frac{\partial^{2}\phi_{\theta}}{\partial \sigma^{2} \partial \beta}(t)$ and $\frac{\partial^{2}\phi_{\theta}}{\partial \beta \partial \sigma^{2}}(t)$ such that (\ref{toshow1}) and (\ref{toshow2}) is satisfied follows with 
\begin{gather*}
    L_{\theta}(1,2) = \frac{2(\kappa-1)}{\beta}\frac{c_{\nu,\kappa-1}}{c_{\nu,\kappa}}, \quad C_{\theta}(1,2) = \beta,
\end{gather*}
and 
\begin{gather*}
    L(1,2) = \frac{2(\kappa-1)}{\beta_{\min}}\frac{c_{\nu,\kappa-1}}{c_{\nu,\kappa}}, \quad C(1,2) = \beta_{\max}.
\end{gather*}
Using (\ref{partialbeta}), (\ref{partialsbeta2}) and (\ref{partialsigma}) we further have that
\begin{gather*}
    \frac{\partial^{3}\phi_{\theta}}{\partial \sigma^{2}\partial \beta^{2}}(t) =  \frac{\partial^{3}\phi_{\theta}}{\partial \beta \partial \sigma^{2} \partial \beta}(t) = \frac{\partial^{3}\phi_{\theta}}{\partial \sigma^{2} \partial \beta \partial \beta}(t) = \frac{1}{\sigma^{2}}q_{1}(t),
\end{gather*}
and thus 
\begin{gather*}
    \frac{\partial^{3}\phi_{\theta}}{\partial \sigma^{2}\partial \beta^{2}} \in \mathcal{C}_{\text{C}}(\mathbb{R}_{+};S_{\theta}(1,2,2))\text{, }\supnorm{\frac{\partial^{3}\phi_{\theta}}{\partial \sigma^{2}\partial \beta^{2}}} \leq L_{\theta}(1,2,2), 
\end{gather*}
with $S_{\theta}(1,2,2) = S_{\theta}(2,2,1) = S_{\theta}(2,1,2) = \left[0,\beta\right]$, and 
\begin{gather*}
    L_{\theta}(1,2,2) = L_{\theta}(2,2,1) = L_{\theta}(2,1,2) = \frac{4(\kappa-1)(\kappa-2)}{\beta^{2}}\frac{c_{\nu,\kappa-2}}{c_{\nu,\kappa}}.
\end{gather*}
Further, (\ref{toshow2}) is satisfied with $C(1,2,2) = C(2,2,1) = C(2,1,2) = \beta_{\max}$, and 
\begin{gather*}
    L(1,2,2) = L(2,2,1) = L(2,1,2) = \frac{4(\kappa-1)(\kappa-2)}{\beta_{\min}^{2}}\frac{c_{\nu,\kappa-2}}{c_{\nu,\kappa}}.
\end{gather*}
Finally we can notice that 
\begin{gather*}
    \frac{\partial^{3}\phi_{\theta}}{\partial \beta \partial\left(\sigma^{2}\right)^{2}}(t) =  \frac{\partial^{3}\phi_{\theta}}{\partial \sigma^{2} \partial \beta \partial \sigma^{2}}(t) = \frac{\partial^{3}\phi_{\theta}}{\partial \sigma^{2} \partial \sigma^{2} \partial \beta}(t) = 0,
\end{gather*}
and hence we can verify the existence of constants 
\begin{gather*}
    L_{\theta}(1,1,2) = L_{\theta}(1,2,1) = L_{\theta}(2,1,1), \quad C_{\theta}(1,1,2) = C_{\theta}(1,2,1) = C_{\theta}(2,1,1),
\end{gather*}
and 
\begin{gather*}
    L(1,1,2) = L(1,2,1) = L(2,1,1), \quad C(1,1,2) = C(1,2,1) = C(2,1,1), 
\end{gather*}
for $\frac{\partial^{3}\phi_{\theta}}{\partial \beta \partial\left(\sigma^{2}\right)^{2}}$, $\frac{\partial^{3}\phi_{\theta}}{\partial \sigma^{2} \partial \beta \partial \sigma^{2}}$, and $\frac{\partial^{3}\phi_{\theta}}{\partial \sigma^{2} \partial \sigma^{2} \partial \beta}$, such that (\ref{toshow1}) and (\ref{toshow2}) is satisfied. Thus, we have shown that for $\kappa > 4$, $\{\phi_{\theta} \colon \theta \in \Theta\}$ satisfies \ref{6.1.2} of Assumption~\ref{cov_a2_iso}, where for any $q = 1,2,3$, $i_{1}, \dotsc, i_{q} \in \{1, \dotsc,p\}$, $\mathcal{B}_{\text{C}}(\mathbb{R}_{+};S_{\theta}(i_{1}, \dotsc, i_{q}))$ can be replaced with $\mathcal{C}_{\text{C}}(\mathbb{R}_{+};S_{\theta}(i_{1}, \dotsc, i_{q}))$. It now remains to show that \ref{6.1.3} of Assumption~\ref{cov_a2_iso} is satisfied. We already know, since $\phi_{\theta} \in \Phi_{d}$, that $w_{\theta}$ is continuous and non-negative definite on $\mathbb{R}^{d}$. We write $L_{1}\left(\mathbb{R}_{+}\right)$ and $L_{1}\left(\mathbb{R}^{d}\right)$ for the spaces of Lebesgue integrable functions on $\mathbb{R}_{+}$ and $\mathbb{R}^{d}$, respectively. Since $t \mapsto t^{d-1}\phi_{\theta}(t) \in L_{1}\left(\mathbb{R}_{+}\right)$ we have that $w_{\theta} \in L_{1}\left(\mathbb{R}^{d}\right)$. Thus we can conclude, using for example Theorems $5.26$ and $6.18$ in \cite{wendland_2004}, that for any $s \in \mathbb{R}^{d}$, $\widehat{w}_{\theta}(s) = \mathcal{F}_{d}\phi_{\theta}\left(\nnorm{s}\right) > 0$, where 
\begin{gather*}
    \mathcal{F}_{d}\phi_{\theta}\left(t\right) = t^{1-(d/2)}\int_{0}^{\infty}\phi_{\theta}(u)u^{d/2}J_{(d/2)-1}(tu)du, \quad t \in \mathbb{R}_{+},
\end{gather*}
with $J_{(d/2)-1}$ the Bessel function of order $(d/2)-1$. This also shows $s \mapsto \widehat{w}_{\theta}(s)$ is uniformly continuous on $\mathbb{R}^{d}$, a member of $L_{1}\left(\mathbb{R}^{d}\right)$ and Fourier inversion holds (see for example Theorem $1.1$ and Corollary $1.26$ in \cite{10.2307/j.ctt1bpm9w6}). It remains to check that $\Theta \times \mathbb{R}^{d} \ni \left(\theta, s\right) \mapsto \widehat{w}_{\theta}(s)$ is continuous. In the present case, where $\kappa > 0$ and $\nu \geq (d+1)/2 + \kappa$, one has actually already established a closed form representation of $\widehat{w}_{\theta}(s)$. We can refer to Theorem $2.1$ in \cite{CHERNIH201417} (see also Theorem $1$ in \cite{bevilacqua2019estimation} for a nice summary and further results) and write for $s \in \mathbb{R}^{d}\setminus\left\{0\right\}$, 
\begin{equation*} \label{closedformwFourier}
    \frac{\widehat{w}_{\theta}(s)}{\left(2 \pi \right)^{d}\sigma^{2}L_{\zeta}\beta^{d}} = {_{1}}F_{2}\bigg(\frac{d+1}{2} + \kappa; \frac{d+1}{2} + \kappa + \frac{\nu}{2},\frac{d+1}{2} + \kappa + \frac{\nu}{2} + \frac{1}{2}; -\frac{\left(\nnorm{s}\beta\right)^{2}}{4} \bigg),
\end{equation*}
where with $\zeta \coloneqq \left(\nu, \kappa, d\right)$, $L_{\zeta} = K^{\zeta}\Gamma(\kappa)/2^{1-\kappa}\operatorname{B}(2\kappa, \nu + 1)$, with
\begin{gather*}
    K^{\zeta} = \frac{2^{-\kappa-d+1}\pi^{-\frac{d}{2}}\Gamma\left(\nu + 1\right)\Gamma\left(2\kappa + d\right)}{\Gamma\left(\kappa + \frac{d}{2}\right)\Gamma\left(\nu + 2\left(\frac{d+1}{2} + \kappa\right)\right)},
\end{gather*}
and for any $z \in \mathbb{R}$,
\begin{gather*}
    {_{1}}F_{2}\left(a; b,c; z\right) = \sum_{k=0}^{\infty}\frac{\left(a\right)_{k}z^{k}}{\left(b\right)_{k}\left(c\right)_{k}k!},
\end{gather*}
a special case of the generalized hypergeometric functions ${_{1}}F_{2}$ (see also \cite{10.5555/1098650}), where for $k \in \natnum$, $\left(q\right)_{k} = \Gamma(q+k)/\Gamma(q)$ denotes the Pochhammer symbol. Note that for $z \in \mathbb{R}$, $\norm{z} \geq 1$ ($z \neq 1$),  ${_{1}}F_{2}\left(a; b,c; z\right)$ is defined via its analytic continuation. Since we know that $s \mapsto \widehat{w}_{\theta}(s)$ is continuous on the entire $\mathbb{R}^{d}$ and 
\begin{gather*}
    z \mapsto {_{1}}F_{2}\left(\frac{d+1}{2} + \kappa; \frac{d+1}{2} + \kappa + \frac{\nu}{2},\frac{d+1}{2} + \kappa + \frac{\nu}{2} + \frac{1}{2}; z \right),
\end{gather*}
is continuous in $0$, we can further note that
\begin{align*}
    \widehat{w}_{\theta}(0) &= \left(2 \pi \right)^{d}\sigma^{2}L_{\zeta}\beta^{d}{_{1}}F_{2}\left(\frac{d+1}{2} + \kappa; \frac{d+1}{2} + \kappa + \frac{\nu}{2},\frac{d+1}{2} + \kappa + \frac{\nu}{2} + \frac{1}{2}; 0 \right)\\
    &= \left(2 \pi \right)^{d}\sigma^{2}L_{\zeta}\beta^{d}.
\end{align*}
This then shows that $\Theta \times \mathbb{R}^{d} \ni \left(\theta, s\right) \mapsto \widehat{w}_{\theta}(s)$ is continuous as a composition of continuous functions and hence the proposition is proven.
\end{proof}

\begin{proof} [Proof of Proposition~\ref{check2}]
We first show that Assumption~\ref{cov_a3} is satisfied. We write $\theta_{1} = \left(\sigma_{1}^{2}, \beta_{1}\right)$ and $\theta_{2} = \left(\sigma_{2}^{2}, \beta_{2}\right)$ and show that $\theta_{1} \neq \theta_{2}$ implies that $\phi_{\theta_{1}}(\nnorm{h}) \neq \phi_{\theta_{2}}(\nnorm{h})$ for all $h \in \mathrm{B}\left(0;\min\left\{\beta_{1}, \beta_{2}\right\}\right)\setminus\{0\}$. Suppose first that $\beta_{1} = \beta_{2}$ but $\sigma_{1}^{2} \neq \sigma_{2}^{2}$ we then have that $\phi_{\theta_{1}}(\nnorm{h}) \neq \phi_{\theta_{2}}(\nnorm{h})$ for all $h \in \mathrm{B}\left(0;\min\left\{\beta_{1}, \beta_{2}\right\}\right)$, since for any $h \in \mathrm{B}\left(0;\min\left\{\beta_{1}, \beta_{2}\right\}\right)$, $\nnorm{h}/\beta_{1} = \nnorm{h}/\beta_{2} < 1$ and thus  $\phi_{\nu,\kappa}\left(\nnorm{h}/\beta_{1}\right) = \phi_{\nu,\kappa}\left(\nnorm{h}/\beta_{2}\right) > 0$. Suppose now that either $\sigma_{1}^{2} \neq \sigma_{2}^{2}$, with $\sigma_{2}^{2} < \sigma_{1}^{2}$ but $\beta_{1} \neq \beta_{2}$, or $\sigma_{1}^{2} = \sigma_{2}^{2} = \sigma^{2}$ but $\beta_{1} \neq \beta_{2}$. Then, let us assume that $\min\left\{\beta_{1}, \beta_{2}\right\} = \beta_{2}$. We have with $\nnorm{h}/\beta_{1} < \nnorm{h}/\beta_{2}$, by monotonicity of $r \mapsto \phi_{\nu, \kappa}(r)$, that
either 
\begin{gather*}
    \phi_{\theta_{1}}(\nnorm{h}) - \phi_{\theta_{2}}(\nnorm{h}) = \left(\phi_{\nu, \kappa}\left(\frac{\nnorm{h}}{\beta_{1}}\right) - \frac{\sigma_{2}^{2}}{\sigma_{1}^{2}}\phi_{\nu, \kappa}\left(\frac{\nnorm{h}}{\beta_{2}}\right)\right) > 0 
\end{gather*}
for all $h \in \mathrm{B}\left(0;\min\left\{\beta_{1}, \beta_{2}\right\}\right)\setminus\{0\}$ or
\begin{gather*}
    \phi_{\theta_{1}}(\nnorm{h}) - \phi_{\theta_{2}}(\nnorm{h}) = \sigma^{2}\left(\phi_{\nu, \kappa}\left(\frac{\nnorm{h}}{\beta_{1}}\right) - \phi_{\nu, \kappa}\left(\frac{\nnorm{h}}{\beta_{2}}\right)\right) > 0
\end{gather*}
for all $h \in \mathrm{B}\left(0;\min\left\{\beta_{1}, \beta_{2}\right\}\right)\setminus\{0\}$. When $\min\left\{\beta_{1}, \beta_{2}\right\} = \beta_{1}$, we will in either of the above cases have $\phi_{\theta_{1}}(\nnorm{h}) - \phi_{\theta_{2}}(\nnorm{h}) < 0$ for all $h \in \mathrm{B}\left(0;\min\left\{\beta_{1}, \beta_{2}\right\}\right)\setminus\{0\}$. Further, we can also use a similar argument for the case where either $\sigma_{1}^{2} \neq \sigma_{2}^{2}$, with $\sigma_{2}^{2} > \sigma_{1}^{2}$ but $\beta_{1} > \beta_{2}$ or $\beta_{1} < \beta_{2}$, or $\sigma_{1}^{2} = \sigma_{2}^{2} = \sigma^{2}$ but $\beta_{1} > \beta_{2}$ or $\beta_{1} < \beta_{2}$. Thus we have shown that $\theta_{1} \neq \theta_{2}$ implies that $\phi_{\theta_{1}}(\nnorm{h}) \neq \phi_{\theta_{2}}(\nnorm{h})$ for all $h \in \mathrm{B}\left(0;\min\left\{\beta_{1}, \beta_{2}\right\}\right)\setminus\{0\}$. Then, for $\tau = 0$, since $\min\left\{\beta_{1}, \beta_{2}\right\} > 1$, $\mathrm{B}\left(0;\min\left\{\beta_{1}, \beta_{2}\right\}\right)\setminus\{0\}$ at least contains integers $z \in \left\{p \in \mathbb{Z}^{d} \colon \nnorm{p} = 1\right\}$. Therefore $\theta_{1} \neq \theta_{2}$ implies $\phi_{\theta_{1}}(z) \neq \phi_{\theta_{2}}(z)$ on $\left\{p \in \mathbb{Z}^{d} \colon \nnorm{p} = 1\right\}$. If $\tau \in \left(0, 1/2\right)$, since $\min\left\{\beta_{1}, \beta_{2}\right\} > 0$, $\mathrm{B}\left(0;\min\left\{\beta_{1}, \beta_{2}\right\}\right) \cap D_{\tau}$ has non zero Lebesgue measure. We have thus shown that Assumption~\ref{cov_a3} is satisfied. Let us now show that $\{w_{\theta} \colon \theta \in \Theta\}$ also satisfies Assumptions~\ref{cov_a4}. To do so, fix some interval $I = (0,b] \subset \mathbb{R}_{+}$, where $1-2\tau < b < \beta$. We will show that for any $\theta \in \Theta$, there exists $t_{0} \in I$ such that
\begin{equation} \label{wronskian}
    W\left(\frac{\partial \phi_{\theta}}{\partial \sigma^{2}}, \frac{\partial \phi_{\theta}}{\partial \beta}\right)(t_{0}) = \ddet{\begin{pmatrix} \frac{\partial \phi_{\theta}}{\partial \sigma^{2}}(t_{0}) & \frac{\partial \phi_{\theta}}{\partial \beta}(t_{0}) \\ \frac{d \frac{\partial \phi_{\theta}}{\partial \sigma^{2}}}{dt}(t_{0}) & \frac{d \frac{\partial \phi_{\theta}}{\partial \beta}}{dt}(t_{0})  \end{pmatrix}} \neq 0,
\end{equation}
where $W\big(\frac{\partial \phi_{\theta}}{\partial \sigma^{2}}, \frac{\partial \phi_{\theta}}{\partial \beta}\big)(t_{0})$ is called the Wronskian of $t \mapsto \frac{\partial \phi_{\theta}}{\partial \sigma^{2}}(t)$ and $t \mapsto \frac{\partial \phi_{\theta}}{\partial \beta}(t)$ at $t_{0} \in I$. This then shows that the functions $t \mapsto \frac{\partial \phi_{\theta}}{\partial \sigma^{2}}(t)$ and $t \mapsto \frac{\partial \phi_{\theta}}{\partial \beta}(t)$ are linearly independent on the entire interval $I$, more explicitly, for any $t \in I$, 
\begin{gather*}
    \alpha_{1}\frac{\partial \phi_{\theta}}{\partial \sigma^{2}}(t) + \alpha_{2}\frac{\partial \phi_{\theta}}{\partial \beta}(t) = 0
\end{gather*}
will imply that $\alpha_{1} = \alpha_{2} = 0$. This then shows that there does not exist $\left(\alpha_{1}, \alpha_{2}\right) \in \mathbb{R}^{2}\setminus\{0\}$, such that for any $\theta \in \Theta$, 
\begin{gather*}
h \mapsto \alpha_{1}\frac{\partial \phi_{\theta}}{\partial \sigma^{2}}(\nnorm{h}) + \alpha_{2}\frac{\partial \phi_{\theta}}{\partial \beta}(\nnorm{h}) = 0, 
\end{gather*}
a.e.\ with respect to the Lebesgue measure on on $\mathrm{B}\left[0;b\right]\setminus\{0\}$. This then justifies, for both cases, either $\tau = 0$, or $\tau > 0$, that also Assumption~\ref{cov_a4} must be satisfied. Hence, let us show (\ref{wronskian}). We can calculate, using arguments from the proof of Proposition~\ref{check1}, that for $t \in I$,
\begin{align*}
    \frac{\partial \phi_{\theta}}{\partial \sigma^{2}}(t) &= \phi_{\nu,\kappa}\left(\frac{t}{\beta}\right), \\
    \frac{\partial \phi_{\theta}}{\partial \beta}(t) &= \frac{2t^{2}(\kappa-1)}{\beta^{3}}\frac{c_{\nu,\kappa-1}}{c_{\nu,\kappa}}\sigma^{2}\phi_{\nu,\kappa-1}\left(\frac{t}{\beta}\right),\\
    \frac{d\frac{\partial \phi_{\theta}}{\partial \sigma^{2}}}{d t}(t) &= -\frac{2t(\kappa-1)}{\beta}\frac{c_{\nu,\kappa-1}}{c_{\nu,\kappa}}\phi_{\nu,\kappa-1}\left(\frac{t}{\beta}\right),
\end{align*}
and
\begin{gather*}
    \frac{d\frac{\partial \phi_{\theta}}{\partial \beta}}{d t}(t) = \frac{4t(\kappa-1)}{\beta^{3}}\frac{\sigma^{2}}{c_{\nu,\kappa}}\left(\ c_{\nu,\kappa-1}\phi_{\nu,\kappa-1}\left(\frac{t}{\beta}\right) - \frac{t^{2}(\kappa-2)}{\beta}c_{\nu,\kappa-2}\phi_{\nu,\kappa-2}\left(\frac{t}{\beta}\right)\right).
\end{gather*}
Therefore we have that for $t \in I$, $W\big(\frac{\partial \phi_{\theta}}{\partial \sigma^{2}}, \frac{\partial \phi_{\theta}}{\partial \beta}\big)(t)$ is given by
\begin{equation*}
    \begin{split}
    \phi_{\nu,\kappa}\left(\frac{t}{\beta}\right)\frac{4t(\kappa-1)}{\beta^{3}}\frac{\sigma^{2}}{c_{\nu,\kappa}}\left(\!c_{\nu,\kappa-1}\phi_{\nu,\kappa-1}\left(\frac{t}{\beta}\right)-\frac{t^{2}(\kappa-2)}{\beta}c_{\nu,\kappa-2}\phi_{\nu,\kappa-2}\left(\frac{t}{\beta}\right)\!\right) \\
    \quad + \frac{4t^{3}(\kappa-1)^{2}}{\beta^{4}}\sigma^{2}\left(\frac{c_{\nu,\kappa-1}}{c_{\nu,\kappa}}\phi_{\nu,\kappa-1}\left(\frac{t}{\beta}\right)\right)^{2}.
\end{split}
\end{equation*}
But the latter expression is not equal to zero on the entire $I$. To see it, assume by contradiction that indeed $ W\big(\frac{\partial \phi_{\theta}}{\partial \sigma^{2}}, \frac{\partial \phi_{\theta}}{\partial \beta}\big)(t) = 0$ for all $t \in I$. Using standard algebraic manipulations one can show that this is equivalent to assume that the function
\begin{gather*}
    g(t) \coloneqq \frac{f_{1}(t)}{f_{2}(t)},
\end{gather*}
with
\begin{gather*}
    f_{1}(t) \coloneqq (\kappa-2)c_{\nu,\kappa}\phi_{\nu,\kappa}\bigg(\frac{t}{\beta}\bigg)c_{\nu,\kappa-2}\phi_{\nu,\kappa-2}\bigg(\frac{t}{\beta}\bigg) - (\kappa-1)\bigg(c_{\nu,\kappa-1}\phi_{\nu,\kappa-1}\bigg(\frac{t}{\beta}\bigg)\bigg)^{2}
\end{gather*}
and
\begin{gather*}
    f_{2}(t) \coloneqq c_{\nu,\kappa}c_{\nu,\kappa-1}\phi_{\nu,\kappa}\bigg(\frac{t}{\beta}\bigg)\phi_{\nu,\kappa-1}\bigg(\frac{t}{\beta}\bigg),
\end{gather*}
is constant equal to $\beta$ on $I$. But this makes no sense and thus we arrive at a contradiction. Hence, there exists $t_{0} \in I$ such that (\ref{wronskian}) is satisfied, which shows that Assumption~\ref{cov_a4} is satisfied and thus concludes the proof of Proposition~\ref{check2}.
\end{proof}

\begin{proof} [Proof of Proposition~\ref{propositionExample0}]
The goal is to check that $\big\{(\mathfrak{T}_{m, \theta})\colon \theta \in \Theta\big\}$ satisfies Assumption~\ref{uniformly_approx_covfun_iso}, then we conclude using Propositions~\ref{check1} and \ref{check2}, as well as Theorems~\ref{consistencyradial} and \ref{asymnormalradial}. We first notice that for any $\theta \in \Theta$, $m \in \natnum$ and any $q = 1,2,3$, $i_{1}, \dotsc, i_{q} \in \{1, \dotsc,p\}$, $\mathfrak{T}_{m, \theta}$ and $\frac{\partial^{q}\mathfrak{T}_{m, \theta}}{\partial \theta_{i_1}\cdots\partial \theta_{i_q}}$ are Borel measurable functions on $\mathbb{R}_{+}$. In addition, for any $\theta \in \Theta$ and $m \in \natnum$, $\mathfrak{T}_{m, \theta}$ has support $[0,U_{\theta,m}]$, with $U_{\theta,m} = \min\left\{C_{m},\beta\right\}$ that satisfies $\sup_{m \in \natnum}\sup_{\theta \in \Theta}U_{\theta,m} = \beta_{\max}$. Further, one can verify that the family $\big\{(\mathfrak{T}_{m, \theta})\colon \theta \in \Theta\big\}$ is also uniformly bounded by $\sigma^{2}_{\max}$ on $\mathbb{R}_{+}$ and it converges uniformly to $\phi_{\theta}$ on $\mathbb{R}_{+}$, independent of $\theta \in \Theta$, that is $\sup_{\theta \in \Theta}\lVert \mathfrak{T}_{m, \theta} - \phi_{\theta}\rVert_{\infty} \tonormalm 0$. Thus \ref{6.2.1}, \ref{6.2.2} and \ref{6.2.3} of Assumption~\ref{uniformly_approx_covfun_iso} are satisfied. To verify the remaining assumptions, we view $\mathfrak{T}_{m, \theta}(t)$ as the result of a truncation operator $g \mapsto \mathfrak{T}_{m}(g) = g\mathbbm{1}_{[0,C_{m}]}$ evaluated at $t$. That is, $\mathfrak{T}_{m, \theta}(t) = \mathfrak{T}_{m}(\phi_{\theta})(t)$. Then, we remark that for any $q = 1,2,3$, $i_{1}, \dotsc, i_{q} \in \{1, \dotsc,p\}$, for any $\theta \in \Theta$,
\begin{gather*}
    \frac{\partial^{q} \mathfrak{T}_{m, \theta}}{\partial \theta_{i_1}\cdots\partial \theta_{i_q}}(t) = \mathfrak{T}_{m}\left(\frac{\partial^{q}\phi_{\theta}}{\partial \theta_{i_1}\cdots\partial \theta_{i_q}}\right)(t).
\end{gather*}    
Thus, by Proposition~\ref{check1}, also \ref{6.2.4} and \ref{6.2.5} of Assumption~\ref{uniformly_approx_covfun_iso} are satisfied.  

\end{proof}

\begin{proof}[Proof of Proposition~\ref{propositionExample1}]
For any $m \in \natnum$, for any $\theta \in \Theta$, $B_{m, \theta}(t;b_{m})$ is continuous on $[0,M]$ and it is also continuous as a function of $\theta \in \Theta$ (see also the proof of Proposition~\ref{check1}). Further, it converges uniformly to $\phi_{\theta}$ on $[0,M]$, independent of $\theta \in \Theta$. That is
\begin{gather*}
    \sup_{\theta \in \Theta}\sup_{t \in [0,M]}\norm{\mathfrak{P}_{m, \theta}(t) - \phi_{\theta}(t)} \tonormalm 0.
\end{gather*} 
To see this we can rely, for example, on the proof of Theorem 2.3.1 in \cite{lorentz1953bernstein}. There, it is shown that for any $t \in [0,M]$ and $\theta \in \Theta$, for any $\varepsilon > 0$, there exists $\delta(t) > 0$ such that 
\begin{gather*}
    \norm{B_{m, \theta}(t;b_{m}) - \phi_{\theta}(t)} \leq \varepsilon + 2\sigma^{2}_{\max}\frac{b_{m}t}{m\delta(t)^{2}},
\end{gather*}
for $m$ large enough. Since $[0,M]$ is compact and $\phi_{\theta}$ is continuous, we can choose $\delta^{*} \equiv \delta(t)$, independent of $t \in [0,M]$ and $\theta \in \Theta$, such that the above inequality is satisfied for arbitrary $\varepsilon > 0$, with $\delta(t)$ replaced with $\delta^{*}$. Then, we conclude by taking the supremum on the left and right over $[0,M]$ and $\theta \in \Theta$. For any $m \in \natnum$, for any $\theta \in \Theta$, we can write 
\begin{gather*}
    \norm{\mathfrak{P}_{m, \theta}(t) - \phi_{\theta}(t)} = \norm{\mathfrak{P}_{m, \theta}(t) - \phi_{\theta}(t)}\mathbbm{1}_{[0,M]}(t) + \norm{\mathfrak{P}_{m, \theta}(t) - \phi_{\theta}(t)}\mathbbm{1}_{(M,\infty)}(t).
\end{gather*}
Notice that because $M \geq \beta_{\max}$, the latter term is actually zero independent of $\theta \in \Theta$ and thus we have that $\left(\mathfrak{P}_{m, \theta}\right)_{m \in \natnum}$ converges uniformly to $\phi_{\theta}$ on the entire $\mathbb{R}_{+}$, independent of $\theta \in \Theta$. Thus, we have that $\sup_{\theta \in \Theta}\lVert \mathfrak{P}_{m, \theta} - \phi_{\theta}\rVert_{\infty} \tonormalm 0$. Note also that the convergence (in the uniform norm) of $\mathfrak{P}_{m, \theta}$ to $\phi_{\theta}$ in particular implies that the sequence of functions $\left(\mathfrak{P}_{m, \theta}\right){}_{m \in \natnum}$ is bounded on $\mathbb{R}_{+}$ for any $\theta \in \Theta$. Therefore we can use that
\begin{gather*}
    \sup_{m \in \natnum}\sup_{t \in \mathbb{R}_{+}}\sup_{\theta \in \Theta}\mathfrak{P}_{m, \theta}(t) = \sup_{m \in \natnum}\sup_{t \in [0,M]}\sup_{\theta \in \Theta}\mathfrak{P}_{m, \theta}(t),
\end{gather*}
to find $\widetilde{C} \coloneqq M$ and $\widetilde{L} \coloneqq \sup_{m \in \natnum}\sup_{t \in [0,M]}\sup_{\theta \in \Theta}\mathfrak{P}_{m, \theta}(t)$, two constants, which are independent of $m \in \natnum$ and $\theta \in \Theta$ (recall that $\Theta$ is compact), such that \ref{6.2.2} of Assumption~\ref{uniformly_approx_covfun_iso} is satisfied. Clearly, for any $\theta \in \Theta$ and for any $m \in \natnum$, the function $\mathfrak{P}_{m, \theta} \colon (\mathbb{R}_{+}, \mathfrak{B}(\mathbb{R}_{+})) \to (\mathbb{R}, \mathfrak{B}(\mathbb{R}))$ is measurable. In conclusion we have shown that \ref{6.2.1}, \ref{6.2.2} and \ref{6.2.3} of Assumption~\ref{uniformly_approx_covfun_iso} are satisfied. In the proof of Proposition~\ref{check1} we have shown that for any $\theta \in \Theta$, for any $q = 1,2,3$, $i_{1}, \dotsc, i_{q} \in \{1, \dotsc,p\}$, there exist constants $C_{\theta}(i_{1}, \dotsc, i_{q})$, $L_{\theta}(i_{1}, \dotsc, i_{q})< \infty$, such that
\begin{gather*} 
\frac{\partial^{q} \phi_{\theta}}{\partial \theta_{i_1}\cdots\partial \theta_{i_q}} \in \mathcal{C}_{\text{C}}(\mathbb{R}_{+};[0,C_{\theta}(i_{1}, \dotsc, i_{q})]), \quad \supnorm{\frac{\partial^{q} \phi_{\theta}}{\partial \theta_{i_1}\cdots\partial \theta_{i_q}}} \leq L_{\theta}(i_{1}, \dotsc, i_{q}),
\end{gather*}  
where for any $q = 1,2,3$, $i_{1}, \dotsc, i_{q} \in \{1, \dotsc,p\}$, $\sup_{\theta \in \Theta}C_{\theta}(i_{1}, \dotsc, i_{q}) \leq \beta_{\max}$. In addition, we notice that for any $q = 1,2,3$, $i_{1}, \dotsc, i_{q} \in \{1, \dotsc,p\}$, for any $\theta \in \Theta$,
\begin{gather*}
    \frac{\partial^{q} \mathfrak{P}_{m, \theta}}{\partial \theta_{i_1}\cdots\partial \theta_{i_q}}(t) = \mathfrak{P}_{m}\left(\frac{\partial^{q}\phi_{\theta}}{\partial \theta_{i_1}\cdots\partial \theta_{i_q}}\right)(t),
\end{gather*}
where $g \mapsto \mathfrak{P}_{m}(g)$ is the Bernstein polynomial operator for a function $g$ with support included $[0,M]$:
\begin{gather*}
    \mathfrak{P}_{m}(g)(t) = \sum_{k=0}^{m}g\bigg(b_{m}\frac{k}{m}\bigg)\binom{m}{k}\bigg(\frac{t}{b_{m}}\bigg)^{k}\bigg(1-\frac{t}{b_{m}}\bigg)^{m-k},
\end{gather*}
for $t \leq M$ and zero otherwise. Therefore, we can rely on the same arguments that we have used to show that \ref{6.2.2} and \ref{6.2.3} of Assumption~\ref{uniformly_approx_covfun_iso} are satisfied, to show that also \ref{6.2.4} and \ref{6.2.5} of Assumption~\ref{uniformly_approx_covfun_iso} must be satisfied. This then concludes the proof of Proposition~\ref{propositionExample1}.   
\end{proof}

\begin{proof}[Proof of Proposition~\ref{propositionExample2}]
The proof follows the same reasoning as the proof of Proposition~\ref{propositionExample1}.
\end{proof}

\begin{proof}[Proof of Proposition~\ref{propositionExample4}]
Since $\left(\delta(m)\right){}_{m \in \natnum}$ does not depend on $\theta \in \Theta$ and $t \in \mathbb{R}_{+}$, and is such that $\delta(m) \xrightarrow{} 0$, as $m \xrightarrow{} \infty$ we can see that $\big\{(\mathfrak{S}_{m, \theta})\colon \theta \in \Theta\big\}$ satisfies Assumption~\ref{uniformly_approx_covfun_iso}. Thus, using Propositions~\ref{check1} and \ref{check2}, under application of Theorems~\ref{consistencyradial} and \ref{asymnormalradial}, the proposition is proven. 
\end{proof}

\subsection{Proof of results in Section~\ref{sec:app2}}

\begin{proof}[Proof of Theorem~\ref{thm:taper}]
Given a collection $S_{(n)}$ of $S$, let ${K_{n,\theta}}_{i,j} = k_{\theta}(S_{i}-S_{j})$, $1 \leq i,j \leq n$, denote the $n \times n$ covariance matrix based on the family $\{k_{\theta} \colon \theta \in \Theta\}$. We first note that under the given assumptions on the family $\{k_{\theta} \colon \theta \in \Theta\}$, we have that  
\begin{gather*}
    \sup_{n \in \natnum}\sup_{\theta \in \Theta} \big \lVert K_{n,\theta} \big \rVert_{2} < \infty \text{ and } \inf_{n \in \natnum} \inf_{\theta \in \Theta}\lambda_{n}(K_{n,\theta}) > 0,
\end{gather*}
with $\mathbb{P}$ probability one. This can be seen from Proposition~D.4 and Lemma~D.5 in \cite{bachoc2014asymptotic}. Using this, the proof of \eqref{eq:taper1} is immediate, it follows from Lemmas~\ref{eigenbehave} and~\ref{lemmalogdet}. 

If we proof
\begin{equation} \label{eq:taper3}
    d_{n,\hat{\theta}_{n}(kt)} = \inf_{\theta \in \Theta}d_{n,\theta} + \delta^{\prime}_{n}, \quad \text{as $n \xrightarrow[]{} \infty$,}
\end{equation}
where $\delta^{\prime}_{n}\toinp 0$, \eqref{eq:taper2} follows from \eqref{eq:taper1}, and we are done. We note that \eqref{eq:taper3} is established if we prove
\begin{equation} \label{eq:taper4}
    \sup_{\theta \in \Theta}\norm{l_{n,\text{\tiny t-ML}}(\theta) - \mathbb{E}\left[l_{n,\text{\tiny t-ML}}(\theta)\;\vert\;S_{(n)}\right]}\toinp 0,
\end{equation}
where
\begin{gather*}
    l_{n,\text{\tiny t-ML}}(\theta) \coloneqq \frac{1}{n}\llog{\operatorname{det}\left(R_{n,\theta}\right)} + \frac{1}{n}\big\langle Z_{(n)}, R_{n,\theta}^{-1}Z_{(n)}\big\rangle,
\end{gather*}
the random version of the modified log-likelihood function based on the tapered covariance function. This is seen from the proof of Theorem~3.3 in \cite{bachoc2018asymptotic}. But under the given assumptions, the family $\{k_{\theta}t_{\beta_{0}} \colon \theta \in \Theta\}$ satisfies Assumption~\ref{cov_a2} (regarding~\ref{3.2.2}, up to $q=1$ and the continuity of first order partial derivatives). Thus \eqref{eq:taper4} can be shown as it was shown (see \eqref{consistency2}) in the proof of Theorem~\ref{thm1}. 
\end{proof}

\end{appendix}

\section*{Acknowledgments}
The authors thank Roman Flury for all the stimulating discussions that were held during the development of this work. This work was supported by the Swiss National Science Foundation SNSF-175529.

\bibliographystyle{plain} 
\bibliography{references}

\end{document}